\documentclass[]{nmeauth}

\usepackage{subfigure}
\usepackage{endfloat}

\newcommand{\real}{\mathbb{R}} 
\newcommand{\integer}{\mathbb{N}}
\newcommand{\vectornorm}[1]{\left\|#1\right\|}

\newtheorem{theorem}{Theorem}[section]

\usepackage[numbers,sort&compress]{natbib}
\usepackage[lined,boxed]{algorithm2e}

\begin{document}
\NME{0}{0}{00}{00}{00}

\runningheads{M. Arnst, R. Ghanem, E. Phipps and J. Red-Horse}{Stochastic modeling of coupled problems}

\received{4 January 2010}
\norevised{}
\noaccepted{}

\title{Measure transformation and efficient quadrature\\in reduced-dimensional stochastic modeling of coupled problems}

\author{M. Arnst\affil{1}$^{,}$\affil{2}\corrauth, R. Ghanem\affil{2}, E. Phipps\affil{3} and J. Red-Horse\affil{3}}

\address{\affilnum{1} B52/3, Universit\'{e} de Li\`{e}ge, Chemin des Chevreuils 1, B-4000 Li\`{e}ge, Belgium.\\
                 \affilnum{2} 210 KAP Hall, University of Southern California, Los Angeles, CA 90089, USA.\\
                 \affilnum{3} Sandia National Laboratories\footnotemark[2], P.O. Box 5800, Albuquerque, NM 87185, USA.}

\corraddr{B52/3, Universit\'{e} de Li\`{e}ge, Chemin des Chevreuils 1, B-4000 Li\`{e}ge, Belgium.}
\footnotetext[2]{Sandia National Laboratories is a multi-program laboratory managed and operated by Sandia Corporation, a wholly owned subsidiary of Lockheed Martin Corporation, for the U.S. Department of Energy's National Nuclear Security Administration under contract DE-AC04-94AL85000.}

\begin{abstract}
Coupled problems with various combinations of multiple physics, scales, and domains are found in numerous areas of science and engineering.
A key challenge in the formulation and implementation of corresponding coupled numerical models is to facilitate the communication of information across physics, scale, and domain interfaces, as well as between the iterations of solvers used for response computations.
In a probabilistic context, any information that is to be communicated between subproblems or iterations should be characterized by an appropriate probabilistic representation. 
Although the number of sources of uncertainty can be expected to be large in most coupled problems, our contention is that exchanged probabilistic information often resides in a considerably lower dimensional space than the sources themselves.  
In this work, we thus use a dimension-reduction technique for obtaining the representation of the exchanged information. 
The main subject of this work is the investigation of a measure-transformation technique that allows implementations to exploit this dimension reduction to achieve computational gains.
The effectiveness of the proposed dimension-reduction and measure-transformation methodology is demonstrated through a multiphysics problem relevant to nuclear engineering.
\end{abstract}


\section{Introduction}
The modeling and simulation of coupled systems governed by multiple physical processes that may exist simultaneously across multiple scales and domains are critical tools for addressing numerous challenges encountered in many areas of science and engineering. 
However, models are, by definition, only approximations of their target scenarios and are thus prone to modeling errors.
Additionally, parametric uncertainties may exist owing to various limitations in manufacturing and experimental methods. 
Uncertainty quantification~(UQ) thus constitutes a key requirement for achieving realistic predictive simulations.

Probability theory provides a rigorous mathematical framework for UQ, which permits a unified treatment of modeling errors and parametric uncertainties.
The first step in a probabilistic UQ analysis typically involves using methods from mathematical statistics~\citep{cramer1946,kullback1968} to characterize the uncertain features associated with a model as one or more random variables, random fields, random matrices, or random operators.  
The second step is to map this probabilistic representation of inputs through the system model into a probabilistic representation of responses.
This can be achieved in several ways, which include Monte Carlo sampling techniques~\citep{robert2005} and stochastic expansion methods.
The latter typically involve the computation of a representation of the predictions as a polynomial chaos~(PC) expansion.
Several approaches are available to calculate the coefficients in this expansion, such as embedded projection~\citep{ghanem2003,soize2004}, nonintrusive projection~\citep{soize2004}, and collocation~\citep{ghanem1998,ghanem1998b,ghiocel2002,xiu2005,babuska2007}.

A key challenge in the formulation and implementation of a coupled model is to facilitate the communication of information across physics, scale, and domain interfaces, as well as between the iterations of solvers used for response computations.
This information can comprise physical properties, energetic quantities, or solution patches, among other quantities.
Although the number of sources of uncertainty can be expected to be large in most coupled problems, we believe that the exchanged information often resides in a considerably lower dimensional space than the sources themselves.
Exchanged information can be expected to have a \textit{low effective stochastic dimension} in multiphysics problems when this information consists of a solution field that has been smoothed by a forward operator and in multiscale problems when this information is obtained by summarizing fine-scale quantities into a coarse-scale representation.

In a previous paper~\citep{arnst2011a}, we had proposed the use of a \textit{dimension-reduction} technique to represent the exchanged information: we proposed to represent the exchanged information by an adaptation of the Karhunen-Lo\`{e}ve~(KL) decomposition as this information passes from subproblem to subproblem and from iteration to iteration.
The main objective of this paper is to complement this dimension-reduction technique by a \textit{measure-transformation} technique that allows implementations to exploit the dimension reduction to achieve computational gains.
Here, we specifically consider implementations that use stochastic expansion methods.
The proposed measure-transformation technique allows such implementations to carry out key algorithmic operations, such as the construction of PC expansions, with respect to the reduced stochastic degrees of freedom of the exchanged uncertainty representations, thus leading to a solution in a reduced-dimensional space, which in turn reduces the computational cost.

The organization of this paper is as follows.
First, in Sec.~\ref{sec:sec2}, we outline the proposed dimension-reduction and measure-transformation methodology.
Then, in Secs.~\ref{sec:sec3a}--\ref{sec:sec4}, we present the dimension-reduction and measure-transformation techniques.
In Sec.~\ref{sec:sec5}, we provide details on the implementation of these techniques.
Finally, in Secs.~\ref{sec:sec6} and~\ref{sec:sec7}, we demonstrate the proposed methodology through an illustration problem.

\section{Dimension-reduction and measure-transformation methodology}\label{sec:sec2}

\subsection{Model problem}
This paper is devoted to the solution of a stochastic coupled model of the following form:
\begin{equation}
\label{eq:coupling12}\begin{aligned}
\boldsymbol{f}(\boldsymbol{u},\boldsymbol{x},\boldsymbol{\xi})=\boldsymbol{0},&&\boldsymbol{y}=\boldsymbol{h}(\boldsymbol{u},\boldsymbol{\xi}),&&\boldsymbol{f}:\real^{r}\times\real^{s_{0}}\times\real^{m}\rightarrow\real^{r},&&\boldsymbol{h}:\real^{r}\times\real^{m}\rightarrow\real^{r_{0}},\\
\boldsymbol{g}(\boldsymbol{y},\boldsymbol{v},\boldsymbol{\zeta})=\boldsymbol{0},&&\boldsymbol{x}=\boldsymbol{k}(\boldsymbol{v},\boldsymbol{\zeta}),&&\boldsymbol{g}:\real^{r_{0}}\times\real^{s}\times\real^{n}\rightarrow\real^{s},&&\boldsymbol{k}:\real^{s}\times\real^{n}\rightarrow\real^{s_{0}}.
\end{aligned}
\end{equation}
To avoid certain technicalities involved in infinite-dimensional representations, we assume that these equations are discretized representations of a stochastic model that couples two physics, two scales, two domains, or a combination of these subproblems.
For instance, these equations may be obtained from the spatial discretization of a steady-state problem, or they may be the equations obtained at a single time step after the spatial and temporal discretization of an evolution problem.
Further, we assume that the data of the first subproblem, which enter this subproblem as coefficients or loadings or both, depend on a finite number of uncertain real parameters denoted as $\xi_{1},\ldots,\xi_{m}$, and that the data of the second subproblem depend on a finite number of uncertain real parameters denoted as~$\zeta_{1},\ldots,\zeta_{n}$.
Lastly, we model these sources of uncertainty as random variables and collect them into vectors, $\boldsymbol{\xi}=(\xi_{1},\ldots,\xi_{m})$ and~$\boldsymbol{\zeta}=(\zeta_{1},\ldots,\zeta_{n})$, which are assumed to be defined on a probability triple $(\Theta,\mathcal{T},P)$ and considered to have values in $\real^{m}$ and $\real^{n}$, respectively.

The stochastic coupled model~(\ref{eq:coupling12}) is a general bidirectionally coupled model. 
The \textit{solution variables}, $\boldsymbol{u}$, of the first subproblem, $\boldsymbol{f}$, depend on the solution variables, $\boldsymbol{v}$, of the second subproblem, $\boldsymbol{g}$, through \textit{coupling variables}, $\boldsymbol{x}$; likewise, $\boldsymbol{v}$ depends on~$\boldsymbol{u}$ through~$\boldsymbol{y}$.

Thus, to solve this stochastic coupled model, we require to find the random variables~$\boldsymbol{u}$ and~$\boldsymbol{v}$ defined on~$(\Theta,\mathcal{T},P)$ with values in~$\real^{r}$ and~$\real^{s}$ such that~(\ref{eq:coupling12}) is satisfied
under the assumption that the stochastic coupled model is well-posed in that it admits a unique and stable solution.

\subsection{Partitioned iterative solution}
Because a coupled model usually characterizes its response only in an implicit manner, the numerical solution of a coupled model typically requires an iterative method. 
Here, we assume that iterative methods and associated solvers already exist for each subproblem, and we therefore consider a \textit{partitioned} iterative method that reuses the aforementioned separate solvers as steps in a global iterative method built around them to obtain a solution to the coupled model.
Let us assume that each of the aforementioned separate iterative methods is based on the reformulation of the associated subproblem as a fixed-point problem:
\begin{equation}
\label{eq:couplingA12}\begin{aligned}
&\boldsymbol{u}=\boldsymbol{a}(\boldsymbol{u},\boldsymbol{x},\boldsymbol{\xi}),&&\boldsymbol{y}=\boldsymbol{h}(\boldsymbol{u},\boldsymbol{\xi}),&&\boldsymbol{a}:\real^{r}\times\real^{s_{0}}\times\real^{m}\rightarrow\real^{r},&&\boldsymbol{h}:\real^{r}\times\real^{m}\rightarrow\real^{r_{0}},\\
&\boldsymbol{v}=\boldsymbol{b}(\boldsymbol{y},\boldsymbol{v},\boldsymbol{\zeta}),&&\boldsymbol{x}=\boldsymbol{k}(\boldsymbol{v},\boldsymbol{\zeta}),&&\boldsymbol{b}:\real^{r_{0}}\times\real^{s}\times\real^{n}\rightarrow\real^{s},&&\boldsymbol{k}:\real^{s}\times\real^{n}\rightarrow\real^{s_{0}}.
\end{aligned}
\end{equation}
It should be noted that these equations can be obtained by setting~$\boldsymbol{a}(\boldsymbol{u},\boldsymbol{v},\boldsymbol{\xi})=\boldsymbol{u}-\boldsymbol{f}(\boldsymbol{u},\boldsymbol{v},\boldsymbol{\xi})$ and~$\boldsymbol{b}(\boldsymbol{u},\boldsymbol{v},\boldsymbol{\zeta})=\boldsymbol{v}-\boldsymbol{g}(\boldsymbol{u},\boldsymbol{v},\boldsymbol{\zeta})$, but that alternative reformulations, such as those involving a direct solution of the subproblems or of their linear approximations, are often better adapted.
We then consider the solution of the stochastic coupled model by a \textit{Gauss-Seidel} iterative method using suitable initial values~$\boldsymbol{u}^{0}$, $\boldsymbol{v}^{0}$, and~$\boldsymbol{x}^{0}=\boldsymbol{k}(\boldsymbol{v}^{0},\boldsymbol{\zeta})$ as follows:
\begin{equation}
\label{eq:couplingAGS12}\begin{aligned}
&\boldsymbol{u}^{\ell}=\boldsymbol{a}\big(\boldsymbol{u}^{\ell-1},\boldsymbol{x}^{\ell-1},\boldsymbol{\xi}\big),&&\quad\quad\quad\boldsymbol{y}^{\ell}=\boldsymbol{h}(\boldsymbol{u}^{\ell},\boldsymbol{\xi}),\\
&\boldsymbol{v}^{\ell}=\boldsymbol{b}\big(\boldsymbol{y}^{\ell},\boldsymbol{v}^{\ell-1},\boldsymbol{\zeta}\big),&&\quad\quad\quad\boldsymbol{x}^{\ell}=\boldsymbol{k}(\boldsymbol{v}^{\ell},\boldsymbol{\zeta}).
\end{aligned}
\end{equation}
This is not the only partitioned iterative method available; however, for simplicity, we employ only this method in this work.
It should be noted that although we implement the proposed methodology using the Gauss-Seidel iterative method, one can readily use the proposed methodology with other iterative methods such as Jacobi, relaxation, and Newton methods.

\subsection{Dimension reduction}\label{sec:dimensionreduction}
We believe that exchanged information often resides in a considerably lower dimensional space than the input sources of uncertainty themselves.
Therefore, in~\citep{arnst2011a}, we investigated the effectiveness of dimension-reduction techniques for the representation of the exchanged information.
Rather than exchanging the coupling variables $\boldsymbol{x}^{\ell}$ and $\boldsymbol{y}^{\ell}$ and the solution variables $\boldsymbol{u}^{\ell}$ and $\boldsymbol{v}^{\ell}$ in their original form, we proposed to approximate these random variables by a truncated KL decomposition as they pass from subproblem to subproblem and from iteration to iteration.
The use of this dimension-reduction technique leads to the solution of the stochastic coupled model by a Gauss-Seidel iterative method as follows:
\begin{equation}
\label{eq:couplingAGSred12}\begin{aligned}
&\hat{\boldsymbol{u}}{}^{\ell}=\boldsymbol{a}\big(\hat{\boldsymbol{u}}{}^{\ell-1,e},\hat{\boldsymbol{x}}{}^{\ell-1,e},\boldsymbol{\xi}\big),&&\quad\quad\quad\hat{\boldsymbol{y}}{}^{\ell}=\boldsymbol{h}(\hat{\boldsymbol{u}}{}^{\ell},\boldsymbol{\xi}),\\
&\hat{\boldsymbol{v}}{}^{\ell}=\boldsymbol{b}\big(\hat{\boldsymbol{y}}{}^{\ell,d},\hat{\boldsymbol{v}}^{\ell-1,d},\boldsymbol{\zeta}\big),&&\quad\quad\quad\hat{\boldsymbol{x}}{}^{\ell}=\boldsymbol{k}(\hat{\boldsymbol{v}}{}^{\ell},\boldsymbol{\zeta}),
\end{aligned}
\end{equation}
where $\boldsymbol{q}^{\ell,d}=[\hat{\boldsymbol{y}}{}^{\ell,d};\hat{\boldsymbol{v}}{}^{\ell-1,d}]$ and $\boldsymbol{r}^{\ell,e}=[\hat{\boldsymbol{u}}{}^{\ell-1,e};\hat{\boldsymbol{x}}{}^{\ell-1,e}]$ are truncated KL decompositions of $\boldsymbol{q}^{\ell}=[\hat{\boldsymbol{y}}{}^{\ell};\hat{\boldsymbol{v}}{}^{\ell-1}]$ and $\boldsymbol{r}^{\ell}=[\hat{\boldsymbol{u}}{}^{\ell-1};\hat{\boldsymbol{x}}{}^{\ell-1}]$, respectively, which are written as follows:
\begin{equation}
\label{eq:KLqqqq}\begin{aligned}
&\boldsymbol{q}^{\ell,d}=\overline{\boldsymbol{q}}{}^{\ell}+\sum_{j=1}^{d}\sqrt{\lambda_{j}^{\ell}}\eta_{j}^{\ell}\boldsymbol{\phi}^{j,\ell},\\
&\boldsymbol{r}^{\ell,e}=\overline{\boldsymbol{r}}{}^{\ell}+\sum_{j=1}^{e}\sqrt{\kappa_{j}^{\ell}}\iota_{j}^{\ell}\boldsymbol{\psi}^{j,\ell}.
\end{aligned}
\end{equation}
These decompositions are described in detail in later sections.  
Truncation of KL decompositions most often results in approximation errors, and we thus use a hat superscript to distinguish the successive approximations determined by~(\ref{eq:couplingAGSred12}) from those determined by~(\ref{eq:couplingAGS12}).

It should be noted that~(\ref{eq:KLqqqq}) provides a combined reduced-dimensional representation of~$\hat{\boldsymbol{y}}{}^{\ell}$ and~$\hat{\boldsymbol{v}}{}^{\ell-1}$ in terms of a single set of reduced random variables~$\boldsymbol{\eta}^{\ell}=(\eta_{1}^{\ell},\ldots,\eta_{d}^{\ell})$ and a combined reduced-dimensional representation of~$\hat{\boldsymbol{u}}{}^{\ell-1}$ and~$\hat{\boldsymbol{x}}{}^{\ell-1}$ in terms of a single set of reduced random variables~$\boldsymbol{\iota}^{\ell}=(\iota_{1}^{\ell},\ldots,\iota_{e}^{\ell})$.
However, this is not the only construction of a reduced-dimensional representation that could be considered.  
One can readily use the proposed methodology with other dimension-reduction techniques such as those involving the construction of a separate reduced-dimensional representation of the coupling and solution variables, with each representation having its own reduced random variables.  

Finally, although our notations do not express a potential dependence of $d$ and $e$ on $\ell$, it should be noted that the reduced dimensions can be allowed to depend on the iteration. 

\subsection{Measure transformation}\label{sec:measuretransformation}
The successive approximations determined by the iterative method~(\ref{eq:couplingAGS12}) that does not involve dimension reduction can be constructed as random variables of the following form:
\begin{equation}
\label{eq:usefulnessbefore12}
\begin{aligned}
&{\boldsymbol{u}}{}^{\ell}(\theta)\equiv{\boldsymbol{u}}{}^{\ell}\big(\boldsymbol{\xi}(\theta),\boldsymbol{\zeta}(\theta)\big),\\
&{\boldsymbol{v}}{}^{\ell}(\theta)\equiv{\boldsymbol{v}}{}^{\ell}\big(\boldsymbol{\xi}(\theta),\boldsymbol{\zeta}(\theta)\big);
\end{aligned}
\end{equation}
i.e., $\boldsymbol{u}{}^{\ell}$ and~$\boldsymbol{v}{}^{\ell}$ can be constructed as transformations of the input random variables~$\boldsymbol{\xi}=(\xi_{1},\ldots,\xi_{m})$ and~$\boldsymbol{\zeta}=(\zeta_{1},\ldots,\zeta_{n})$.
The random variables~$\boldsymbol{u}{}^{\ell}$ and~$\boldsymbol{v}{}^{\ell}$ thus exist in a solution space of stochastic dimension $m+n$.
Consequently, implementations of~(\ref{eq:couplingAGS12}) using stochastic expansion methods would typically involve the approximation of $\boldsymbol{u}{}^{\ell}$ and~$\boldsymbol{v}{}^{\ell}$ by finite-dimensional representations as follows:
\begin{equation}
\label{eq:usefulnessbefore12b}
\begin{aligned}
&{\boldsymbol{u}}{}^{\ell,p}=\sum_{|\boldsymbol{\alpha}|=0}^{p}\boldsymbol{u}{}^{\ell}_{\boldsymbol{\alpha}}\psi_{\boldsymbol{\alpha}}\big(\boldsymbol{\xi},\boldsymbol{\zeta}\big),&&\quad\quad\quad\boldsymbol{u}{}^{\ell}_{\boldsymbol{\alpha}}\in\real^{r},\\
&{\boldsymbol{v}}{}^{\ell,p}=\sum_{|\boldsymbol{\alpha}|=0}^{p}\boldsymbol{v}{}^{\ell}_{\boldsymbol{\alpha}} \psi_{\boldsymbol{\alpha}}\big(\boldsymbol{\xi},\boldsymbol{\zeta}\big),&&\quad\quad\quad\boldsymbol{v}{}^{\ell}_{\boldsymbol{\alpha}}\in\real^{s},
\end{aligned}
\end{equation}
where~$\{\psi_{\boldsymbol{\alpha}},\boldsymbol{\alpha}\in\integer^{m+n}\}$ is a Hilbertian basis for the Hilbert space of~$P_{(\boldsymbol{\xi},\boldsymbol{\zeta})}$-square-integrable functions from~$\real^{m+n}$ into~$\real$, in which~$P_{(\boldsymbol{\xi},\boldsymbol{\zeta})}$ denotes the joint probability distribution of the input random variables~$\boldsymbol{\xi}$ and~$\boldsymbol{\zeta}$, and~$|\boldsymbol{\alpha}|=\alpha_{1}+\ldots+\alpha_{m+n}$. 
Then, the task of the solution algorithm would be to compute the coordinates~$\boldsymbol{u}{}^{\ell}_{\boldsymbol{\alpha}}$ and~$\boldsymbol{v}{}^{\ell}_{\boldsymbol{\alpha}}$ in these expansions. 
Depending on the specific nonintrusive projection, embedded projection, or collocation method that is chosen, this task would require the construction of basis functions, quadrature rules, moment tensors, or a combination of these with respect to the joint probability distribution~$P_{(\boldsymbol{\xi},\boldsymbol{\zeta})}$.

In contrast, owing to the dimension reduction, the successive approximations determined by the iterative method~(\ref{eq:couplingAGSred12}) can be constructed as random variables of the following form:
\begin{equation}
\label{eq:usefulness12}\begin{aligned}
&\hat{\boldsymbol{u}}{}^{\ell}(\theta)\equiv\hat{\boldsymbol{u}}{}^{\ell}\big(\boldsymbol{\xi}(\theta),\boldsymbol{\iota}^{\ell}(\theta)\big),\\
&\hat{\boldsymbol{v}}{}^{\ell}(\theta)\equiv\hat{\boldsymbol{v}}{}^{\ell}\big(\boldsymbol{\eta}^{\ell}(\theta),\boldsymbol{\zeta}(\theta)\big);
\end{aligned}
\end{equation}
i.e., $\hat{\boldsymbol{u}}{}^{\ell}$ can be constructed as a transformation of the input random variables~$\boldsymbol{\xi}=(\xi_{1},\ldots,\xi_{m})$ and the reduced random variables~$\boldsymbol{\iota}^{\ell}=(\iota_{1}^{\ell},\ldots,\iota_{e}^{\ell})$ and~$\hat{\boldsymbol{v}}{}^{\ell}$ can be constructed as a transformation of~$\boldsymbol{\eta}^{\ell}=(\eta_{1}^{\ell},\ldots,\eta_{d}^{\ell})$ and~$\boldsymbol{\zeta}=(\zeta_{1},\ldots,\zeta_{n})$.
Thus, the random variable~$\hat{\boldsymbol{u}}{}^{\ell}$ exists in a space of stochastic dimension~$m+e$ and the random variable $\hat{\boldsymbol{v}}{}^{\ell}$ exists in a space of stochastic dimension~$d+n$.
Consequently, implementations can exploit the dimension reduction to approximate $\boldsymbol{u}{}^{\ell}$ and~$\boldsymbol{v}{}^{\ell}$ by finite-dimensional representations as follows:
\begin{equation}
\label{eq:usefulness12b}\begin{aligned}
&\hat{\boldsymbol{u}}{}^{\ell,q}=\sum_{|\boldsymbol{\beta}|=0}^{q}\hat{\boldsymbol{u}}{}^{\ell}_{\boldsymbol{\beta}}\Gamma{}_{\boldsymbol{\beta}}^{\ell}\big(\boldsymbol{\xi},\boldsymbol{\iota}^{\ell}\big),&&\quad\quad\quad\hat{\boldsymbol{u}}{}^{\ell}_{\boldsymbol{\beta}}\in\real^{r},\\
&\hat{\boldsymbol{v}}{}^{\ell,q}=\sum_{|\boldsymbol{\gamma}|=0}^{q}\hat{\boldsymbol{v}}{}^{\ell}_{\boldsymbol{\gamma}}\Upsilon{}_{\boldsymbol{\gamma}}^{\ell}\big(\boldsymbol{\eta}^{\ell},\boldsymbol{\zeta}\big),&&\quad\quad\quad\hat{\boldsymbol{v}}{}^{\ell}_{\boldsymbol{\gamma}}\in\real^{s},
\end{aligned}
\end{equation}
where~$\{\Gamma{}_{\boldsymbol{\beta}}^{\ell},\boldsymbol{\beta}\in\integer^{m+e}\}$ and~$\{\Upsilon{}_{\boldsymbol{\gamma}}^{\ell},\boldsymbol{\gamma}\in\integer^{d+n}\}$ are Hilbertian bases for the Hilbert spaces of~$P_{(\boldsymbol{\xi},\boldsymbol{\iota}^{\ell})}$-square-integrable functions from~$\real^{m+e}$ and of~$P_{(\boldsymbol{\eta}^{\ell},\boldsymbol{\zeta})}$-square integrable functions from~$\real^{d+n}$, respectively, into~$\real$.
Here, $P_{(\boldsymbol{\xi},\boldsymbol{\iota}^{\ell})}$ and~$P_{(\boldsymbol{\eta}^{\ell},\boldsymbol{\zeta})}$ are the joint probability distributions of~$\boldsymbol{\xi}$ and~$\boldsymbol{\iota}^{\ell}$ and of~$\boldsymbol{\eta}^{\ell}$ and~$\boldsymbol{\zeta}$, respectively.
Now, depending on the specific solution method that is chosen, the construction of these expansions requires the construction of basis functions, quadrature rules, moment tensors, or a combination of these with respect to the joint probability distributions~$P_{(\boldsymbol{\xi},\boldsymbol{\iota}^{\ell})}$ and~$P_{(\boldsymbol{\eta}^{\ell},\boldsymbol{\zeta})}$.
This approach involves a \textit{measure transformation} because from the dimension reduction given by~(\ref{eq:KLqqqq}), it follows that the probability distributions of the reduced random variables $\boldsymbol{\iota}^{\ell}$ and~$\boldsymbol{\eta}^{\ell}$ are transformations of the probability distribution of the input random variables~$\boldsymbol{\xi}$ and~$\boldsymbol{\zeta}$, as described in detail in a later section. 

Throughout this work, we have employed Hilbertian bases that are constituted of polynomials and we thus refer to them as polynomial chaos (PC) bases. 

Although our notations do not express a potential dependence of $p$ and $q$ on the subproblem or $\ell$, it should be noted that the subsets of basis functions used to construct the finite-dimensional representations can be allowed to depend on the subproblem and the iteration.

\subsection{Effectiveness of the dimension-reduction and measure-transformation methodology}
The key feature of the proposed methodology is that it enables a solution of the subproblems in a reduced-dimensional space when the exchanged information has a low effective stochastic dimension.
Specifically, a solution in a reduced-dimensional space is enabled when the reduced dimensions can be selected such that $d<m$ and $e<n$ while sufficient accuracy is maintained; refer to~(\ref{eq:usefulnessbefore12}) and~(\ref{eq:usefulness12}).
This benefit is of particular significance for implementations of stochastic coupled models that use stochastic expansion methods.
These methods suffer from a curse of dimensionality in that their computational cost increases quickly with an increase in the stochastic dimension.
The proposed methodology addresses the curse of dimensionality by mitigating the increase in stochastic dimension when information is exchanged.

Finally, it should be noted that the proposed methodology can readily be adapted to meet various requirements of specific applications. 
One could, for instance, use the KL decomposition to represent only those exchanged random variables that are of low effective stochastic dimension and use another representation for the remaining variables; further, one could implement a measure transformation to solve only those subproblems whose computational cost would thus be lowered and implement another approach for the remaining subproblems. 

\section{Karhunen-Loeve decomposition}\label{sec:sec3a}

Here, we concisely recall the KL decomposition used in~\citep{arnst2011a} to construct a reduced-dimensional representation of a random variable~$\boldsymbol{q}$ that is defined on a probability triple~$(\Theta,\mathcal{T},P)$, takes values in a finite-dimensional Euclidean space~$\real^{w}$, and is of the second order:
\begin{equation}
\int_{\Theta}\vectornorm{\boldsymbol{q}}^{2}dP<+\infty,
\end{equation}
where~$\vectornorm{\cdot}$ denotes the Euclidean norm.
Refer to~\citep{arnst2011a} and the references therein for more details on this construction. 

\subsection{Second-order descriptors}
The mean vector~$\overline{\boldsymbol{q}}$ and covariance matrix~$\boldsymbol{C}_{\boldsymbol{q}}$ of the second-order random variable~$\boldsymbol{q}$ are defined as the $w$-dimensional vector and square matrix, respectively, such that
\begin{align}
\label{eq:meanu}\overline{\boldsymbol{q}}&=\int_{\Theta}\boldsymbol{q}dP,\\
\label{eq:covaru}\boldsymbol{C}_{\boldsymbol{q}}&=\int_{\Theta}(\boldsymbol{q}-\overline{\boldsymbol{q}})(\boldsymbol{q}-\overline{\boldsymbol{q}})^{\mathrm{T}}dP.
\end{align}

\subsection{Reduced-dimensional representation}
Let $\boldsymbol{W}$ be a $w$-dimensional, square, symmetric, and positive definite matrix; we refer to this matrix as the weighting matrix and comment on its usefulness later.
Then, because~$\boldsymbol{C}_{\boldsymbol{q}}$ is symmetric and positive semidefinite, the solution of the generalized eigenproblem
\begin{equation}
\boldsymbol{W}^{\mathrm{T}}\boldsymbol{C}_{\boldsymbol{q}}\boldsymbol{W}\boldsymbol{\phi}^{j}=\lambda_{j}\boldsymbol{W}\boldsymbol{\phi}^{j}\label{eq:eigenproblemeee}
\end{equation}
provides a set of~$w$ eigenvalues~$\lambda_{1}\geq\lambda_{2}\geq\ldots\geq\lambda_{w}\geq 0$ in addition to~$w$ eigenvectors~$\boldsymbol{\phi}^{1},\ldots,\boldsymbol{\phi}^{w}$, which constitute a $\boldsymbol{W}$-weighted orthonormal basis of $\real^{w}$ such that
\begin{equation}
(\boldsymbol{\phi}^{i})^{\mathrm{T}}\boldsymbol{W}\boldsymbol{\phi}^{j}=\delta_{ij},\label{eq:orthonromalkfa3}
\end{equation}
where $\delta_{ij}$ is the Kronecker delta, which is equal to 1 if $i=j$ and 0 otherwise.
The KL-type decomposition of~$\boldsymbol{q}$ is then given by
\begin{equation}
\boldsymbol{q}=\overline{\boldsymbol{q}}+\sum_{j=1}^{w}\sqrt{\lambda_{j}}\eta_{j}\boldsymbol{\phi}^{j},\label{eq:KLvectorrr}
\end{equation}
where the~$\eta_{j}$ are random variables defined on~$(\Theta,\mathcal{T},P)$, with values in~$\real$, such that
\begin{equation}
\eta_{j}=\frac{1}{\sqrt{\lambda_{j}}}(\boldsymbol{q}-\overline{\boldsymbol{q}})^{\mathrm{T}}\boldsymbol{W}\boldsymbol{\phi}^{j}.
\end{equation}
These reduced random variables~$\eta_{j}$ are zero-mean and uncorrelated:
\begin{align}
\label{eq:eta1}&\int_{\Theta}\eta_{j}dP=0,\\
\label{eq:eta2}&\int_{\Theta}\eta_{i}\eta_{j}dP=\delta_{ij}.
\end{align}
The truncation of~(\ref{eq:KLvectorrr}) after~$d$ terms provides a reduced-dimensional representation as follows:
\begin{equation}
\boldsymbol{q}^{d}=\overline{\boldsymbol{q}}+\sum_{j=1}^{d}\sqrt{\lambda_{j}}\eta_{j}\boldsymbol{\phi}^{j}.\label{eq:KLred}
\end{equation}
Although the reduced random variables of the KL decomposition are uncorrelated, it should be noted that they are in general statistically dependent.
Further, although the joint probability distribution of the reduced random variables is usually complicated, it is determined by the KL decomposition.
Finally, if the random variable to be reduced is Gaussian, then the reduced random variables are also Gaussian; however, the joint probability distribution of the reduced random variables is most often not a ``labeled" probability distribution.

\subsection{Concluding remarks}
In this section, we presented a version of the KL decomposition that features a weighting matrix.
We had demonstrated in~\citep{arnst2011a} that this weighting matrix is particularly useful for the construction of a reduced-dimensional representation of a random variable that solves a space-time discretized stochastic model: we showed that by appropriately choosing the weighting matrix as the Gram matrix of the discretization basis, a reduced-dimensional representation is obtained which is consistent with the function-analytic structure that the stochastic model exhibited prior to its discretization.  
Nevertheless, it should be noted that the standard KL decomposition is recovered by simply setting the weighting matrix equal to the identity matrix.

\section{Polynomial chaos with respect to arbitrary probability distributions}\label{sec:sec3}

Here, we consider the construction of a PC basis with respect to an arbitrary, although fully specified, probability distribution~$P_{\boldsymbol{\chi}}$ defined on a finite-dimensional Euclidean space~$\real^{z}$.

\subsection{Objective}
In the context of the proposed methodology, we consider $P_{\boldsymbol{\chi}}$ as the probability distribution of a random variable~$\boldsymbol{\chi}$ that collects the sources of uncertainty that enter a subproblem of a stochastic coupled model.
With reference to~(\ref{eq:usefulness12b}), at a specific iteration, $\boldsymbol{\chi}$ could collect the components of $\boldsymbol{\xi}$ and $\boldsymbol{\iota}$ with~$z=m+e$ or the components of $\boldsymbol{\eta}$ and $\boldsymbol{\zeta}$ with~$z=d+n$, in which case the sought PC basis would be needed to build a PC expansion of $\boldsymbol{u}$ or $\boldsymbol{v}$, respectively. 

\subsection{Methodology}
If the probability distribution~$P_{\boldsymbol{\chi}}$ does not exhibit statistical dependence and its marginal probability distributions are ``labeled" univariate probability distributions, 
then a PC basis can readily be constructed by tensorization of the classical univariate orthonormal polynomials that are associated with these univariate probability distributions in the Askey table~\citep{wiener1938,cameron1947,ghanem2003,xiu2003,soize2004}.
Well-known examples include the use of PC bases of Hermite and Legendre polynomials with respect to Gaussian and uniform probability distributions, respectively. 

However, here, we are interested in a general case wherein the probability distribution~$P_{\boldsymbol{\chi}}$ may exhibit statistical dependence and may not be ``labeled."
This case is considered in the context of the proposed methodology because the reduced random variables of a KL decomposition are in general statistically dependent and not ``labeled."

If~$P_{\boldsymbol{\chi}}$ exhibits statistical dependence or is not ``labeled," the Hilbert space of~$P_{\boldsymbol{\chi}}$-square-integrable functions from~$\real^{z}$ into~$\real$ does not necessarily admit a PC basis, as discussed in the following sections.
Further, polynomials that are orthonormal with respect to an arbitrary probability distribution usually cannot be read from tables in the literature, thus requiring a computational construction.
Therefore, the approach adopted in this section is as follows.
First, we describe a computational construction of a set of orthonormal polynomials.
Then, we discuss conditions under which this set is also a PC basis.

\subsection{Set constituted of~$P_{\boldsymbol{\chi}}$-orthonormal polynomials}
We use standard notations involving multi-indices to work with multivariate polynomials.
We refer to elements~$\boldsymbol{\beta}=(\beta_{1},\ldots,\beta_{z})$ of~$\integer^{z}$ as multi-indices.  
A (multivariate) monomial~$\boldsymbol{\chi}^{\boldsymbol{\beta}}$ associated with a multi-index~$\boldsymbol{\beta}$ is then a function from~$\real^{z}$ into~$\real$ defined by
\begin{equation}
\boldsymbol{\chi}^{\boldsymbol{\beta}}=\chi{}_{1}^{\beta_{1}}\times\ldots\times\chi{}_{z}^{\beta_{z}}.
\end{equation}
Further, we refer to the number~$|\boldsymbol{\beta}|=\beta_{1}+\ldots+\beta_{z}$ as the modulus of the multi-index~$\boldsymbol{\beta}$ and also as the total degree of the corresponding monomial~$\boldsymbol{\chi}^{\boldsymbol{\beta}}$. 
A (multivariate) polynomial is then a function from~$\real^{z}$ into~$\real$ that maps any~$\boldsymbol{\chi}$ to a finite sum $\sum_{\boldsymbol{\beta}}c_{\boldsymbol{\beta}}\boldsymbol{\chi}^{\boldsymbol{\beta}}$ with real coefficients~$c_{\boldsymbol{\beta}}$.
Let~$\mathcal{P}_{z}^{q}$ be the space of all polynomials in~$z$ variables with a total degree of at most~$q$:
\begin{equation}
\mathcal{P}_{z}^{q}=\bigg\{\pi\;:\;\boldsymbol{\chi}=(\beta_{1},\ldots,\beta_{z})\mapsto\pi(\boldsymbol{\chi})=\sum_{|\boldsymbol{\beta}|=0}^{q}c_{\boldsymbol{\beta}}\boldsymbol{\chi}^{\boldsymbol{\beta}}\bigg\};\label{eq:Pnp} 
\end{equation} 
the dimension of~$\mathcal{P}_{z}^{q}$ is commonly expressed as
\begin{equation}
\text{dim}\big(\mathcal{P}^{q}_{z}\big)=\begin{pmatrix}z+q\\z\end{pmatrix}=\frac{(z+q)!}{z!q!}.
\end{equation}

In order to construct a set of~$P_{\boldsymbol{\chi}}$-orthonormal polynomials in~$\mathcal{P}^{q}_{z}$ for a given total degree~$q$, we propose a procedure wherein we first arrange the monomials~$\{\boldsymbol{\chi}^{\boldsymbol{\beta}},\;0\leq|\boldsymbol{\beta}|\leq q\}$ spanning~$\mathcal{P}^{q}_{z}$ in a sequence and then orthonormalize this sequence by the Gram-Schmidt method as follows:  
\begin{align}
\label{eq:GS1}\widetilde{\Gamma}_{\boldsymbol{\beta}}(\boldsymbol{\chi})&=\boldsymbol{\chi}^{\boldsymbol{\beta}}-\sum_{\boldsymbol{\gamma}<\boldsymbol{\beta}}c_{\boldsymbol{\beta\gamma}}\widetilde{\Gamma}_{\boldsymbol{\gamma}}(\boldsymbol{\beta}),\quad\quad\quad\text{with}\quad c_{\boldsymbol{\beta\gamma}}=\frac{\int_{\real^{z}}\boldsymbol{\chi}^{\boldsymbol{\beta}}\;\widetilde{\Gamma}_{\boldsymbol{\gamma}}(\boldsymbol{\chi})dP_{\boldsymbol{\chi}}}{\int_{\real^{z}}\widetilde{\Gamma}_{\boldsymbol{\gamma}}(\boldsymbol{\chi})\widetilde{\Gamma}_{\boldsymbol{\gamma}}(\boldsymbol{\chi})dP_{\boldsymbol{\chi}}},\\
\label{eq:GS2}\Gamma_{\boldsymbol{\beta}}(\boldsymbol{\chi})&=\widetilde{\Gamma}_{\boldsymbol{\beta}}(\boldsymbol{\chi})\Big/\int_{\real^{z}}\widetilde{\Gamma}_{\boldsymbol{\beta}}(\boldsymbol{\chi})\widetilde{\Gamma}_{\boldsymbol{\beta}}(\boldsymbol{\chi})dP_{\boldsymbol{\chi}}.
\end{align}
Clearly, the Gram-Schmidt method cannot always be executed without difficulty: it will fail when the integral in the numerator of~(\ref{eq:GS1}) does not exist or when the integral in the denominator of~(\ref{eq:GS2}) vanishes.
The former difficulty may occur when a few or all the moments of~$P_{\boldsymbol{\chi}}$ are unbounded.
The latter difficulty corresponds to the case wherein~$\widetilde{\Gamma}_{\boldsymbol{\beta}}$ has a vanishing weight with respect to~$P_{\boldsymbol{\chi}}$ and may occur either when~$P_{\boldsymbol{\chi}}$ is a discrete probability distribution or when~$P_{\boldsymbol{\chi}}$ is degenerate in that it is concentrated on a hypersurface in $\real^{z}$.

One can readily impose suitable conditions on~$P_{\boldsymbol{\chi}}$ to ensure proper execution of the Gram-Schmidt method.
Indeed, if~$P_{\boldsymbol{\chi}}$ has finite moments up to a total degree of~$2q$, i.e.,
\begin{equation}
\int_{\real^{z}}\big|\chi_{1}^{\beta_{1}}\times\ldots\times\chi_{z}^{\beta_{z}}\big|dP_{\boldsymbol{\chi}}<+\infty,\quad0\leq|\boldsymbol{\beta}|\leq 2q\label{eq:finitemoments}
\end{equation}
and if~$P_{\boldsymbol{\chi}}$ is nondegenerate, i.e.,
\begin{equation}
\int_{\real^{z}}\pi(\boldsymbol{\chi})^{2}dP_{\boldsymbol{\chi}}>0,\quad\forall\pi\in\mathcal{P}^{q}_{z},\quad\pi\neq 0,\label{eq:nondegen}
\end{equation}
then the aforementioned difficulties cannot occur; thus, it is then ensured that the Gram-Schmidt method provides a unique set of polynomials~$\Gamma_{\boldsymbol{\beta}}$ that are mutually orthogonal in that~$\int_{\real^{z}}\Gamma_{\boldsymbol{\beta}}(\boldsymbol{\chi})\Gamma_{\boldsymbol{\gamma}}(\boldsymbol{\chi})dP_{\boldsymbol{\chi}}=0$ if~$\boldsymbol{\beta}\neq\boldsymbol{\gamma}$ and normalized so that~$\int_{\real^{z}}|\Gamma_{\boldsymbol{\beta}}(\boldsymbol{\chi})|^{2}dP_{\boldsymbol{\chi}}=1$.

A frequently encountered example of a case wherein~$P_{\boldsymbol{\chi}}$ has finite moments and is non-degenerate is the one in which~$P_{\boldsymbol{\chi}}$ admits a probability density function that has a closed and bounded support with a nonempty interior.

It should be noted that cases wherein~$P_{\boldsymbol{\chi}}$ is a discrete or degenerate probability distribution can easily be handled by simply discarding polynomials with a vanishing $P_{\boldsymbol{\chi}}$-weighted norm.

\subsection{PC basis constituted of~$P_{\boldsymbol{\chi}}$-orthonormal polynomials}\label{sec:hilbert}
The probability distribution~$P_{\boldsymbol{\chi}}$ determines a space~$L^{2}_{P_{\boldsymbol{\chi}}}(\real^{z},\real)$ of functions from~$\real^{z}$ into~$\real$ that are square-integrable with respect to~$P_{\boldsymbol{\chi}}$:
\begin{equation}
\int_{\real^{z}}|f(\boldsymbol{\chi})|^{2}dP_{\boldsymbol{\chi}}<+\infty,\quad\forall f\in L^{2}_{P_{\boldsymbol{\chi}}}(\real^{z},\real).
\end{equation}
The function space~$L^{2}_{P_{\boldsymbol{\chi}}}(\real^{z},\real)$ is a Hilbert space for the inner product
\begin{equation}
\langle f,g\rangle=\int_{\real^{z}}f(\boldsymbol{\chi})g(\boldsymbol{\chi})dP_{\boldsymbol{\chi}},\quad\forall f,g\in L^{2}_{P_{\boldsymbol{\chi}}}(\real^{z},\real).
\end{equation}
If the conditions~(\ref{eq:finitemoments}) and~(\ref{eq:nondegen}) are fulfilled for a specific total degree $q$, the Gram-Schmidt method given by~(\ref{eq:GS1}) and~(\ref{eq:GS2}) provides a unique set~$\{\Gamma_{\boldsymbol{\beta}},0\leq|\boldsymbol{\beta}|\leq q\}$ of~$P_{\boldsymbol{\chi}}$-orthonormal polynomials spanning~$\mathcal{P}^{q}_{z}$.
This set of polynomials allows one to associate to each function~$f$ in~$L^{2}_{P_{\boldsymbol{\chi}}}(\real^{z},\real)$ a corresponding set~$\{f_{\boldsymbol{\beta}},0\leq|\boldsymbol{\beta}|\leq q\}$ of coordinates~$f_{\boldsymbol{\beta}}$ in~$\real$ defined by
\begin{equation}
f_{\boldsymbol{\beta}}=\langle f,\Gamma_{\boldsymbol{\beta}}\rangle=\int_{\real^{z}}f(\boldsymbol{\chi})\Gamma_{\boldsymbol{\beta}}(\boldsymbol{\chi})dP_{\boldsymbol{\chi}},\quad 0\leq|\boldsymbol{\beta}|\leq q.
\end{equation}
This set of coordinates in turn determines an approximation $f^{q}$ of~$f$ in~$\mathcal{P}^{q}_{z}$ as follows:
\begin{equation}
f\approx f^{q}=\sum_{|\boldsymbol{\beta}|=0}^{q}f_{\boldsymbol{\beta}}\Gamma_{\boldsymbol{\beta}}.\label{eq:pce}
\end{equation}
Clearly, if the conditions~(\ref{eq:finitemoments}) and~(\ref{eq:nondegen}) are fulfilled for any total degree $q$, then the above construction provides a $q$-indexed family of the approximations $f^{q}$ of~$f$.
Then , the convergence of these approximations with respect to the total degree~$q$ is a critical issue. 
Ideally, one would like the approximation~$f^{q}$ to converge to~$f$ as~$q$ is increased, i.e.,
\begin{equation}
\lim_{q\rightarrow\infty}\int_{\real^{z}}\bigg|f(\boldsymbol{\chi})-\sum_{|\boldsymbol{\beta}|=0}^{q}f_{\boldsymbol{\beta}}\Gamma_{\boldsymbol{\beta}}(\boldsymbol{\chi})\bigg|^{2}dP_{\boldsymbol{\chi}}=0,\quad\forall f\in L^{2}_{P_{\boldsymbol{\chi}}}(\real^{z},\real).\label{eq:completeness1}
\end{equation}
Now, a set~$\{\Gamma_{\boldsymbol{\beta}},\;\boldsymbol{\beta}\in\integer^{z}\}$ of~$P_{\boldsymbol{\chi}}$-orthonormal functions in~$L^{2}_{P_{\boldsymbol{\chi}}}(\real^{z},\real)$ (which need not necessarily be polynomials) 
is called~\citep{reed1980} a \textit{Hilbertian basis} for~$L^{2}_{P_{\boldsymbol{\chi}}}(\real^{z},\real)$ if it satisfies
\begin{equation}
f=\sum_{\boldsymbol{\beta}\in\integer^{z}}f_{\boldsymbol{\beta}}\Gamma_{\boldsymbol{\beta}},\quad\forall f\in L^{2}_{P_{\boldsymbol{\chi}}}(\real^{z},\real)\label{eq:completeness122}
\end{equation}
in the sense that the right-hand side should converge to the left-hand side in the $L^{2}$-norm without depending on the ordering of the multi-indices.
Clearly, the fulfillment of~(\ref{eq:completeness122}) for any ordering implies the fulfillment of~(\ref{eq:completeness1}) for that specific ordering by the total degree.
As mentioned previously, we refer to a Hilbertian basis constituted of polynomials as a \textit{PC basis}.

The family of Hermite polynomials is well known to provide a PC basis for a Hilbert space of functions that are square-integrable with respect to a Gaussian probability distribution~\citep{cameron1947}.
In general, the families of classical orthonormal polynomials are well known to provide PC bases for Hilbert spaces of functions that are square-integrable with respect to the corresponding ``labeled" probability distributions in the Askey table~\citep{cameron1947,ghanem2003,xiu2003,soize2004}. 
However, it is not obvious for a set of orthonormal polynomials, such as the one determined by the Gram-Schmidt method given by~(\ref{eq:GS1}) and~(\ref{eq:GS2}), to be a PC basis; this property generally depends on the probability distribution under consideration.

If the support of the probability distribution~$P_{\boldsymbol{\chi}}$ is a nonempty, closed, and bounded subset of the $z$-dimensional Euclidean space~$\real^{z}$, then by the Stone-Weierstrass theorem~\citep{aliprantis2006}, any set of~$P_{\boldsymbol{\chi}}$-orthonormal polynomials that spans the space of all the $z$-variate polynomials, including the set of orthonormal polynomials determined by the Gram-Schmidt method given by~(\ref{eq:GS1}) and~(\ref{eq:GS2}), 
forms a PC basis for the Hilbert space of functions that are square-integrable with respect to this probability distribution.
The illustration problem considered in the last two sections of this paper falls within this case. 

On the other hand, cases wherein~$P_{\boldsymbol{\chi}}$ is a probability distribution with an unbounded support are more complex.
For brevity, here, we only give one example to illustrate that sets of orthonormal polynomials defined with respect to a probability distribution with unbounded support need not form a PC basis:
consider a Gaussian random variable~$\xi$; it can then be shown~\citep{riesz1926,stoyanov2000} that the Hilbert space of functions from~$\real$ into~$\real$ that are square-integrable with respect to the probability distribution of~$\chi=\xi^{2k+1}$ with $k=1,2\ldots$ does not admit a PC basis.
Refer to~\citep{ernst2011} and the references therein for further details on this issue.

\subsection{Bibliographical comments}
Orthogonal polynomials have been a subject of intensive research in the literature; for instance, refer to the recent references~\citep{dunkl2001,gautschi2004,golub2009}.
In the research area of stochastic modeling and analysis, the construction of PC bases with respect to arbitrary probability distributions has already been considered in~\citep{wan2006,witteveen2006,soize2010}: reference~\citep{wan2006} presents a multielement PC basis synthesized from univariate orthonormal polynomials generated by the Stieltjes method; the method given in~\citep{witteveen2006} involves a multidimensional PC basis deduced by tensorization of univariate orthonormal polynomials generated by the Gram-Schmidt method; and reference~\citep{soize2010} proposes the use of a singular value decomposition to obtain a multidimensional PC basis that is orthonormal with respect to a discrete approximation of the probability distribution.

\section{Quadrature rules with respect to arbitrary probability distributions}\label{sec:sec4}
In this section, we consider the construction of a family of quadrature rules for the evaluation of the integrals of continuous functions from a finite-dimensional Euclidean space~$\real^{z}$ into~$\real$ with respect to an arbitrary, although fully specified, probability distribution~$P_{\boldsymbol{\chi}}$ on~$\real^{z}$.

\subsection{Objective}
Again, we consider $P_{\boldsymbol{\chi}}$ as the probability distribution of a random variable~$\boldsymbol{\chi}$ that collects the sources of uncertainty that enter a subproblem of a stochastic coupled model.
With reference to~(\ref{eq:usefulness12b}), at a specific iteration, $\boldsymbol{\chi}$ could collect the components of $\boldsymbol{\xi}$ and $\boldsymbol{\iota}$ with~$z=m+e$ or the components of $\boldsymbol{\eta}$ and $\boldsymbol{\zeta}$ with~$z=d+n$.
The purpose of the family of quadrature rules depends on the specific stochastic expansion method that is selected for discretization.
For instance, if a nonintrusive projection method were adopted, a family of quadrature rules would be required to evaluate the projection of the solution to the subproblem onto the employed PC basis.

\subsection{Methodology}
If the probability distribution~$P_{\boldsymbol{\chi}}$ does not exhibit statistical dependence, 
then an appropriate family of quadrature rules can often be constructed through a full or sparse-grid tensorization of the families of Gaussian quadrature rules associated with these univariate probability distributions.  
Gaussian quadrature rules can be read from tables in the literature for ``labeled" univariate probability distributions and can be computed easily for arbitrary univariate probability distributions~\citep{golub2009}.
Well-known examples include the use of Gauss-Hermite and Gauss-Legendre quadrature rules with respect to Gaussian and uniform probability distributions, respectively.
The families of quadrature rules thus obtained are well known to be efficient in that they provide fast convergence rates for smooth integrands.

However, here, we are interested in a general case wherein the probability distribution~$P_{\boldsymbol{\chi}}$ may exhibit statistical dependence.
This case is considered here because the reduced random variables of a KL decomposition are generally statistically dependent.

General multivariate probability distributions do not necessarily admit Gaussian quadrature rules, as discussed in the following sections.
Further, quadrature rules for integration with respect to arbitrary probability distributions usually cannot be read from tables in the literature, thus requiring a computational construction.
Thus, the approach adopted in this section is as follows.
First, we propose a method for the computational construction of a family of quadrature rules for the evaluation of the integrals of continuous functions from~$\real^{z}$ into~$\real$ with respect to the probability distribution~$P_{\boldsymbol{\chi}}$.
Then, we discuss the efficiency and convergence properties of this method in addition to its relationship to other available methods.

\subsection{Proposed computational construction of quadrature rules}\label{sec:sec41}
Let~$\lambda\geq 1$ be a fixed finite integer that determines the \textit{level} of the quadrature rule to be constructed by requiring its polynomial degree of exactness to be $2\lambda-1$. 
Further, let~$P_{\boldsymbol{\chi}}$ be a probability distribution on~$\real^{z}$ that satisfies
\begin{equation}
\int_{\real^{z}}\big|\chi_{1}^{\beta_{1}}\times\ldots\times\chi_{z}^{\beta_{z}}\big|dP_{\boldsymbol{\chi}}<+\infty,\quad0\leq|\boldsymbol{\beta}|\leq 2\lambda-1.\label{eq:finitemomentsb}
\end{equation}
We then propose the following two-step strategy for the construction of a quadrature rule:
\begin{itemize}
\item First, we construct a quadrature rule of the following form: 
\begin{equation}
\int_{\real^{z}}f(\boldsymbol{\chi})dP_{\boldsymbol{\chi}}\approx\sum_{k=1}^{\tilde{\nu}}f\big(\tilde{\boldsymbol{\chi}}_{k}\big)\tilde{w}_{k},\label{eq:originalcubature}
\end{equation}
which allows integrals relative to~$P_{\boldsymbol{\chi}}$ to be approximated accurately, but may have a very large number~$\tilde{\nu}$ of nodes~$\tilde{\boldsymbol{\chi}}_{1},\ldots,\tilde{\boldsymbol{\chi}}_{\tilde{\nu}}$ and associated weights~$\tilde{w}_{1},\ldots,\tilde{w}_{\tilde{\nu}}$.
\item Then, we construct an \textit{embedded quadrature rule} of the following form:
\begin{equation}
\int_{\real^{z}}f(\boldsymbol{\chi})dP_{\boldsymbol{\chi}}\approx
\sum_{k=1}^{\nu}f\big(\boldsymbol{\chi}_{k}\big)w_{k},\label{eq:cubature}
\end{equation}
which has \textit{level} $\lambda$ in that it integrates all polynomials up to total degree~$2\lambda-1$ exactly:
\begin{equation}
\int_{\real^{z}}\pi(\boldsymbol{\chi})dP_{\boldsymbol{\chi}}=\sum_{k=1}^{\nu}\pi\big(\boldsymbol{\chi}_{k}\big)w_{k},\quad\quad\quad\forall\pi\in\mathcal{P}_{z}^{2\lambda-1}, 
\end{equation}
but uses only a small subset~$\boldsymbol{\chi}_{1}=\tilde{\boldsymbol{\chi}}_{k_{1}},\ldots,\boldsymbol{\chi}_{\nu}=\tilde{\boldsymbol{\chi}}_{k_{\nu}}$ of $\nu\ll\tilde{\nu}$ nodes of the original rule. 
\end{itemize}
The nodes and weights of the original quadrature rule~(\ref{eq:originalcubature}) could, for instance, be constructed by a Monte Carlo integration approach.
Here, the key challenge is rather in the selection of the small subset of nodes to be retained by the embedded quadrature rule.
Clearly, in order to solve this subset selection problem, we require to determine a set of weights $\varpi_{1},\ldots,\varpi_{\tilde{\nu}}$ for the candidate nodes~$\tilde{\boldsymbol{\chi}}_{1},\ldots,\tilde{\boldsymbol{\chi}}_{\tilde{\nu}}$, respectively, such that these weights are mostly zero while they still allow integrals of polynomials up to the prescribed total degree of~$2\lambda-1$ to be evaluated exactly. 
We propose to construct such a \textit{sparse} set of weights through the solution of the following~$L^{1}$-\textit{minimization} problem:
\begin{equation}
\min_{\boldsymbol{\varpi}\in\real^{\tilde{\nu}}}\vectornorm{\boldsymbol{\varpi}}_{L^{1}}
,\quad\quad\text{subject to $\boldsymbol{A}\boldsymbol{\varpi}=\boldsymbol{b}$},\label{eq:L1}
\end{equation}
where $\vectornorm{\boldsymbol{\varpi}}_{L^{1}}=|\varpi_{1}|+|\varpi_{2}|+\ldots+|\varpi_{\tilde{\nu}-1}|+|\varpi_{\tilde{\nu}}|$ and the equality constraints~$\boldsymbol{A}\boldsymbol{\varpi}=\boldsymbol{b}$ serve to ensure that polynomials up to the prescribed total degree of $2\lambda-1$ are integrated exactly:
\begin{align}
\label{eq:A}\boldsymbol{A}&=\begin{bmatrix}
\pi_{1}\big(\tilde{\boldsymbol{\chi}}_{1}\big) & \pi_{1}\big(\tilde{\boldsymbol{\chi}}_{2}\big) & \ldots & \pi_{1}\big(\tilde{\boldsymbol{\chi}}_{\tilde{\nu}-1}\big) & \pi_{1}\big(\tilde{\boldsymbol{\chi}}_{\tilde{\nu}}\big)\\
\vdots & \vdots & & \vdots & \vdots \\
\pi_{\mu}\big(\tilde{\boldsymbol{\chi}}_{1}\big) & \pi_{\mu}\big(\tilde{\boldsymbol{\chi}}_{2}\big) & \ldots & \pi_{\mu}\big(\tilde{\boldsymbol{\chi}}_{\tilde{\nu}-1}\big) & \pi_{\mu}\big(\tilde{\boldsymbol{\chi}}_{\tilde{\nu}}\big)
\end{bmatrix},\\
\label{eq:b}\boldsymbol{b}&=\begin{bmatrix} \int_{\real^{z}}\pi_{1}(\boldsymbol{\chi})dP_{\boldsymbol{\chi}} & \ldots & \int_{\real^{z}}\pi_{\mu}(\boldsymbol{\chi})dP_{\boldsymbol{\chi}}\end{bmatrix}^{\mathrm{T}}.
\end{align}
Here, $\pi_{1},\ldots,\pi_{\mu}$ is a basis for~$\mathcal{P}^{2\lambda-1}_{z}$ with~$\mu=\text{dim}(\mathcal{P}^{2\lambda-1}_{z})$.
Let~$\mathcal{B}\cup\mathcal{N}$ then be a partitioning of the index set~$\{1,\ldots,\tilde{\nu}\}$ associated with a solution~$\boldsymbol{\varpi}$ to~(\ref{eq:L1}) such that the~$k$-th component~$\varpi_{k}$ of that solution is nonzero when~$k$ is included in~$\mathcal{B}$ and vanishes when~$k$ is included in~$\mathcal{N}$.
The nodes of the original quadrature rule labeled by the indices in~$\mathcal{B}=\{k_{1},\ldots,k_{\nu}\}$ are then the nodes~$\boldsymbol{\chi}_{1}=\tilde{\boldsymbol{\chi}}_{k_{1}},\ldots,\boldsymbol{\chi}_{\nu}=\tilde{\boldsymbol{\chi}}_{k_{\nu}}$ to be retained by the embedded quadrature rule with the associated weights~$w_{1}=\varpi_{k_{1}},\ldots,w_{\nu}=\varpi_{k_{\nu}}$, while the nodes labeled by the indices in~$\mathcal{N}$ are to be discarded.
In Sec.~\ref{sec:sec5}, we have demonstrated how such a sparse solution can be obtained using either a simplex ~\citep{nocedal2006} or an interior-point optimization algorithm~\citep{wright1997,nocedal2006,megiddo1991}.
Both of these approaches yield a sparse solution that has at most $\text{dim}(\mathcal{P}^{2\lambda-1}_{z})$ nonzero components.
Thus, the proposed construction yields a quadrature rule that uses at most as many nodes as there are equality constraints imposed in~(\ref{eq:L1}) to ensure the prescribed polynomial exactness.

\subsection{Efficiency}
The proposed construction yields a quadrature rule that uses at most~$\text{dim}(\mathcal{P}^{2\lambda-1}_{z})$ nodes to achieve exactness for all~$z$-variate polynomials up to a total degree of $2\lambda-1$.
An important question thus pertains to the optimality of the proposed construction with respect to the number of nodes used to achieve this degree of exactness.
For quadrature rules used to evaluate one-dimensional integrals, it is well known that a $\lambda$-point Gaussian quadrature rule has a degree of exactness of~$2\lambda-1$ and that there exists no $\lambda$-point quadrature rule that is exact for all polynomials up to a degree of~$2\lambda$.
In contrast, for quadrature rules used to evaluate multidimensional integrals, the minimum number of nodes required for the exact integration of a given number of polynomials is yet unknown in current state-of-the-art mathematics.  

Nevertheless, if the probability distribution~$P_{\boldsymbol{\chi}}$ has a closed and bounded support, it is known that the minimum number of nodes required for a quadrature rule to have a degree of exactness of~$2\lambda-1$ is~\citep{cools2002} greater than or equal to~$\text{dim}(\mathcal{P}^{\lambda-1}_{z})$ and less than or equal to~$\text{dim}(\mathcal{P}^{2\lambda-1}_{z})$.
This upper bound corresponds to the Tchakaloff theorem:
\begin{theorem}\label{theorem:tchakaloff}
(Tchakaloff~\citep{tchakaloff1957}) Let~$P_{\boldsymbol{\chi}}$ be a probability distribution on~$\real^{z}$ that admits a probability density function that is supported by a closed and bounded subset~$\mathcal{K}$ of $\real^{z}$.
Then, for any integer~$2\lambda-1\geq 0$, there exists $\mu=\text{dim}(\mathcal{P}_{z}^{2\lambda-1})$ nodes~$\boldsymbol{\chi}_{1},\ldots,\boldsymbol{\chi}_{\mu}$ in $\mathcal{K}$ and positive weights $w_{1},\ldots,w_{\mu}$ such that the resulting quadrature rule has a degree of exactness of~$2\lambda-1$. 
\end{theorem}
Hence, our construction provides a quadrature rule that uses at most as many nodes as the Tchakaloff theorem indicates are at most required to achieve the prescribed degree of exactness.
It should be noted that this theorem not only ensures the existence of a~$\text{dim}(\mathcal{P}_{z}^{2\lambda-1})$-point quadrature rule with a degree of exactness of~$2\lambda-1$, but also guarantees the existence of such a rule with only positive weights.
In contrast, our construction does not necessarily yield a quadrature rule with only positive weights. 
Further, it should be noted that extensions of this theorem to cases involving probability distributions with unbounded support exist~\citep{putinar1997}.  

\subsection{Convergence}\label{sec:convergenceee}
When~$P_{\boldsymbol{\chi}}$ has finite moments of any order and the assumption~(\ref{eq:finitemoments}) is thus fulfilled for any level~$\lambda$, the proposed construction provides a $\lambda$-indexed family of quadrature rules, each member of which is required to achieve a corresponding degree of exactness of $2\lambda-1$.
An important question then pertains to the convergence of the approximation of an integral of a given function using a sequence of quadrature rules with an increasing degree of exactness.

If the probability distribution~$P_{\boldsymbol{\chi}}$ has a closed and bounded support, an insightful, although not comprehensive, answer is contained in the following theorem:
\begin{theorem}\label{eq:cubatureerror}
(Convergence of polynomial-based quadrature rules~\citep{engels1980}) Let~$P_{\boldsymbol{\chi}}$ be a probability distribution on~$\real^{z}$ that admits a probability density function that is supported by a closed and bounded subset~$\mathcal{K}$ of $\real^{z}$.
Then, a quadrature rule with degree of exactness of~$2\lambda-1$ satisfies for every continuous function~$f$ from~$\mathcal{K}\subset\real^{z}$ into~$\real$ the following inequality:
\begin{equation}
\bigg|\int_{\mathcal{K}}f(\boldsymbol{\chi})dP_{\boldsymbol{\chi}}-\sum_{k=1}^{\nu}f\big(\boldsymbol{\chi}_{k}\big)w_{k}\bigg|\leq\bigg(1+\sum_{k=1}^{\nu}|w_{k}|\bigg)\min_{\pi\in\mathcal{P}^{2\lambda-1}_{z}}\max_{\boldsymbol{\chi}\in\mathcal{K}}|f(\boldsymbol{\chi})-\pi(\boldsymbol{\chi})|.\label{eq:theoremcubaturerror}
\end{equation}
\end{theorem} 
The proof follows from the triangle inequality; in fact, the left-hand side satisfies the inequality
\small{
\begin{equation}
\text{l.h.s.}\leq\underbrace{\bigg|\int_{\mathcal{K}}\hspace{-2mm}f(\boldsymbol{\chi})dP_{\boldsymbol{\chi}}-\hspace{-1mm}\int_{\mathcal{K}}\hspace{-2mm}\pi(\boldsymbol{\chi})dP_{\boldsymbol{\chi}}\bigg|}_{\leq\max_{\boldsymbol{\chi}\in\mathcal{K}}|f(\boldsymbol{\chi})-\pi(\boldsymbol{\chi})|}+\underbrace{\bigg|\int_{\mathcal{K}}\hspace{-2mm}\pi(\boldsymbol{\chi})dP_{\boldsymbol{\chi}}-\hspace{-1mm}\sum_{k=1}^{\nu}\pi(\boldsymbol{\chi}_{k})w_{k}\bigg|}_{=0}+\underbrace{\bigg|\sum_{k=1}^{\nu}\pi(\boldsymbol{\chi}_{k})w_{k}\hspace{-1mm}-\hspace{-1mm}\sum_{k=1}^{\nu}f(\boldsymbol{\chi}_{k})w_{k}\bigg|}_{\leq\sum_{k=1}^{\nu}\hspace{-1mm}|w_{k}|\hspace{-0.5mm}\max_{\boldsymbol{\chi}\in\mathcal{K}}\hspace{-0.5mm}|f(\boldsymbol{\chi})-\pi(\boldsymbol{\chi})|}.
\end{equation}}\normalsize
The result~(\ref{eq:theoremcubaturerror}) shows that the quadrature error is bounded by a product of two factors; the first grows with the sum of the absolute values of the weights, and the second is a bound on the approximation error introduced by the best approximation of the integrand by a polynomial with a total degree of at most~$2\lambda-1$. 
This result justifies our construction~(\ref{eq:L1}) of the embedded quadrature rule.
Our choice of the objective function as the sum of the absolute values of the weights results in the minimization of the first factor mentioned above, contributing to the upper bound on the quadrature error.
Our requirement for the embedded quadrature rule to have a polynomial degree of exactness of~$2\lambda-1$ precisely results in the second factor mentioned above.
This requirement allows the embedded quadrature rule to achieve a fast convergence rate for smooth integrands as the polynomial degree of exactness increases because the approximation error introduced by the best approximation of a smooth function by a polynomial of a specified degree on a closed and bounded set is well known to decrease at a fast rate with an increase in this degree.

\subsection{Bibliographical comments}
Three classes of approaches for multivariate integration have received most attention in the literature, namely, probabilistic and number-theoretic methods, polynomial-based methods, and adaptive techniques.
Probabilistic and number-theoretic methods include Monte Carlo and quasi-Monte Carlo integration; for instance, refer to~\citep{robert2005,calfish1998}.
These methods are well known to exhibit a rather slow rate of convergence as a function of the number of integrand evaluations, but Monte Carlo methods have the advantage that their rate of convergence is independent of the stochastic dimension. 
Polynomial-based integration rules are designed to be exact for a prescribed collection of polynomials; for instance, refer to~\citep{engels1980,dunkl2001,cools2002,davis2007}. 
Polynomial-based methods have the advantage that their rate of convergence increases rapidly with the smoothness of the integrand, but these methods suffer from a curse of dimensionality in that their rate of convergence decreases with an increase in the dimension.
Adaptive methods involve the choice of the nodes and weights in a manner that is dependent on the integrand to achieve fast convergence rates, while still limiting the increase in computational cost as the dimension increases; for instance, refer to~\citep{holtz2010}.

Among the polynomial-based methods, four classes of approaches have been investigated extensively in the literature.
The first class includes methods that rely on the direct numerical solution of the system of nonlinear equations that express the polynomial exactness in order to find the nodes and weights; for instance, refer to~\citep{cools1997}. 
Such methods have been found to yield good results for rather low-dimensional problems that feature symmetries and other invariance properties that can be exploited for simplification. 
The second class includes methods that search for polynomials that vanish at the nodes of the quadrature rule; for instance, refer to~\citep{dunkl2001}.
These methods have facilitated the study of important mathematical properties of quadrature rules using ideal and other theories.
The third class consists of methods that involve the construction of quadrature rules by full or sparse-grid tensorization of suitable univariate rules; for instance, refer to~\citep{holtz2010}.
These methods have found many successful applications in the context of stochastic expansion methods, but are limited in application to probability distributions that do not exhibit statistical dependence. 
Finally, methods belonging to the fourth class involve the selection of a subset of nodes from a given quadrature rule in order to construct a more efficient embedded rule.
Many aspects of polynomial-based methods, including Theorems~\ref{theorem:tchakaloff} and~\ref{eq:cubatureerror}, are studied in~\citep{engels1980,dunkl2001,cools2002,davis2007} and the references therein.  

Finally, it should be noted that references~\citep{davis1967,wilson1970,sommariva2009,xiao2010} propose alternative methods for the construction of embedded quadrature rules.
The method given in~\citep{davis1967,wilson1970} relies on the computation of a basic optimal solution~\citep{nocedal2006} to the linear program $\min_{\boldsymbol{\varpi}}\sum_{k=1}^{\nu}\varpi_{k}$, subject to $\boldsymbol{A}\boldsymbol{\varpi}=\boldsymbol{b}$ and $\boldsymbol{\varpi}\geq\boldsymbol{0}$, where~$\boldsymbol{A}$ and~$\boldsymbol{b}$ are still defined by~(\ref{eq:A}) and~(\ref{eq:b}), respectively.
A drawback of this method is that the existence of a basic optimal solution is generally not guaranteed.
The existence of this solution is ensured~\citep{davis1967} when the original quadrature rule has only positive weights and has the targeted degree of exactness, but a basic optimal solution may fail to exist in other cases because the imposed constraints may then be overly restrictive.
The method given in~\citep{sommariva2009,xiao2010} relies on the solution of~$\boldsymbol{A}\boldsymbol{\varpi}=\boldsymbol{b}$ for a basic solution using a QR factorization with column pivoting.
This method is computationally less costly than the proposed method, but it can generally be expected to provide a quadrature rule for which the sum of the absolute values of the weights is larger and for which the error bound~(\ref{eq:theoremcubaturerror}) is thus less favorable than the one for the quadrature rule obtained using the proposed method. 

\section{Implementation}\label{sec:sec5}
In this section, we provide details on the implementation of the proposed dimension-reduction and measure-transformation methodology using stochastic expansion methods.

\subsection{Karhunen-Loeve decomposition}\label{sec:sec5a}
Here, we describe the implementation of the KL decomposition of a random variable represented by a PC expansion and we show how this implementation in turn naturally provides a PC expansion of the reduced random variables.
Adopting the notations used in Secs.~\ref{sec:sec2} and~\ref{sec:sec3a}, we consider the construction of a reduced-dimensional representation of a second-order random variable~$\boldsymbol{q}^{p}$ with values in~$\real^{w}$ that is represented by a PC expansion of the following form:
\begin{equation}
\boldsymbol{q}^{p}=\sum_{|\boldsymbol{\alpha}|=0}^{p}\boldsymbol{q}_{\boldsymbol{\alpha}}\psi_{\boldsymbol{\alpha}}(\boldsymbol{\xi},\boldsymbol{\zeta}),\quad\boldsymbol{q}_{\boldsymbol{\alpha}}\in \real^{w},\label{eq:upcekl}
\end{equation}
where~$\{\psi_{\boldsymbol{\alpha}},\;\boldsymbol{\alpha}\in\integer^{m+n}\}$ denotes a PC basis for the Hilbert space of~$P_{(\boldsymbol{\xi},\boldsymbol{\zeta})}$-square-integrable functions from~$\real^{m+n}$ into~$\real$ with~$\psi_{\boldsymbol{0}}=1$, as mentioned previously.
We specifically consider the construction of the KL decomposition described in Sec.~\ref{sec:sec3a} that involves a weighting matrix~$\boldsymbol{W}$; the standard KL decomposition can be recovered easily by setting the weighting matrix equal to the identity matrix.
Because of the orthonormality of the~$\psi_{\boldsymbol{\alpha}}$, the mean vector $\overline{\boldsymbol{q}}$ and the covariance matrix $\boldsymbol{C}_{\boldsymbol{q}}$ of~$\boldsymbol{q}^{p}$ can be deduced immediately from the PC coordinates as follows:
\begin{align}
\overline{\boldsymbol{q}}&=\boldsymbol{q}_{\boldsymbol{0}},\\
\boldsymbol{C}_{\boldsymbol{q}}&=\sum_{|\boldsymbol{\alpha}|=1}^{p}\boldsymbol{q}_{\boldsymbol{\alpha}}\boldsymbol{q}_{\boldsymbol{\alpha}}^{\mathrm{T}}.
\end{align}   
Then, the solution of the generalized eigenproblem~$\boldsymbol{W}^{\mathrm{T}}\boldsymbol{C}_{\boldsymbol{q}}\boldsymbol{W}\boldsymbol{\Phi}^{j}=\lambda_{j}\boldsymbol{W}\boldsymbol{\Phi}^{j}$ provides the eigenvalues~$\lambda_{1}\geq\lambda_{2}\geq\ldots\lambda_{w}\geq 0$ and the associated eigenmodes~$\boldsymbol{\phi}^{1},\ldots,\boldsymbol{\phi}^{w}$ required to construct a reduced-dimensional representation~$\boldsymbol{q}^{p,d}$ of $\boldsymbol{q}^{p}$ as follows:
\begin{equation}
\boldsymbol{q}^{p,d}=\overline{\boldsymbol{q}}+\sum_{j=1}^{d}\sqrt{\lambda_{j}}\eta^{p}_{j}\boldsymbol{\phi}^{j},
\end{equation}
where the~$\eta_{j}^{p}$ are random variables with values in~$\real$ such that
\begin{equation}
\eta_{j}^{p}=\frac{1}{\sqrt{\lambda_{j}}}\big(\boldsymbol{q}^{p}-\overline{\boldsymbol{q}}\big)^{\mathrm{T}}\boldsymbol{W}\boldsymbol{\phi}^{j}\label{eq:upcekleta}
\end{equation}
and are zero-mean and uncorrelated.
By substituting~(\ref{eq:upcekl}) in~(\ref{eq:upcekleta}), a representation of each reduced random variable as a PC expansion is immediately obtained: 
\begin{equation}
\eta^{p}_{j}=\sum_{|\boldsymbol{\alpha}|=1}^{p}\eta_{j,\boldsymbol{\alpha}}\psi_{\boldsymbol{\alpha}}(\boldsymbol{\xi},\boldsymbol{\zeta})\quad\text{with}\quad\eta_{j,\boldsymbol{\alpha}}=\frac{1}{\sqrt{\lambda_{j}}}\boldsymbol{q}_{\boldsymbol{\alpha}}^{\mathrm{T}}\boldsymbol{W}\boldsymbol{\phi}^{j},\label{eq:upcerveta}
\end{equation}
thus indicating that the KL decomposition of a PC expansion naturally provides a complete probabilistic characterization of the reduced random variables as a PC expansion.

It should be noted that the proposed methodology provides a representation of the solution and coupling variables associated with the subproblems as PC expansions in a combination of input random variables and reduced random variables of reduced-dimensional representations of exchanged information --- refer to~(\ref{eq:usefulness12}).
Thus, when these reduced random variables are represented by PC expansions in the input random variables themselves --- refer to~(\ref{eq:upcerveta}), the proposed methodology requires the construction of KL decompositions of random variables that are represented by compositions of PC expansions.
The implementation of such KL decompositions falls within the scope of the implementation mentioned above because the composition of two PC expansions can always be written equivalently as a PC expansion of the form~(\ref{eq:upcekl}), although as one that is truncated at a higher total degree.

\subsection{Polynomial chaos with respect to arbitrary probability distributions}\label{sec:orthopol}
Here, we provide details on the implementation of the Gram-Schmidt method described in Sec.~\ref{sec:sec3}.
Adopting the notations used in Sec.~\ref{sec:sec3}, we consider the construction of a set of $z$-variate orthonormal polynomials~$\{\Gamma_{\boldsymbol{\beta}},0\leq|\boldsymbol{\beta}|\leq q\}$ up to a specific total degree~$q$ with respect to the probability distribution~$P_{\boldsymbol{\chi}}$ on the finite-dimensional Euclidean space~$\real^{z}$.
To ensure that the Gram-Schmidt method can be executed properly, we assume that~$P_{\boldsymbol{\chi}}$ satisfies the conditions given by~(\ref{eq:finitemoments}) and~(\ref{eq:nondegen}).
Let the $z$-variate monomials~$\{\boldsymbol{\chi}^{\boldsymbol{\beta}},\;0\leq|\boldsymbol{\beta}|\leq q\}$ then be ordered in a sequence~$\boldsymbol{\chi}^{\boldsymbol{\beta}_{1}},\ldots,\boldsymbol{\chi}^{\boldsymbol{\beta}_{\mu}}$, where~$\mu=\text{dim}(\mathcal{P}_{z}^{q})$ denotes the number of monomials in this sequence, as mentioned previously.
Following the approach given in~\citep{golub1969,golub2009}, we propose to implement the Gram-Schmidt  orthonormalization of this sequence using a Cholesky factorization of its Gram matrix.
The Gram matrix~$\boldsymbol{G}$ of the sequence of monomials~$\boldsymbol{\chi}^{\boldsymbol{\beta}_{1}},\ldots,\boldsymbol{\chi}^{\boldsymbol{\beta}_{\mu}}$ is the~$\mu$-dimensional, square, and symmetric matrix that collects the inner products of these monomials as follows:
\begin{equation}
\boldsymbol{G}=
\begin{bmatrix}
\int_{\real^{z}}\boldsymbol{\chi}^{\boldsymbol{\beta}_{1}}\,\boldsymbol{\chi}^{\boldsymbol{\beta}_{1}}dP_{\boldsymbol{\chi}} & \ldots & \int_{\real^{z}}\boldsymbol{\chi}^{\boldsymbol{\beta}_{1}}\,\boldsymbol{\chi}^{\boldsymbol{\beta}_{\mu}}dP_{\boldsymbol{\chi}}\\
\vdots & & \vdots \\
\int_{\real^{z}}\boldsymbol{\chi}^{\boldsymbol{\beta}_{\mu}}\,\boldsymbol{\chi}^{\boldsymbol{\beta}_{1}}dP_{\boldsymbol \chi} & \ldots & \int_{\real^{z}}\boldsymbol{\chi}^{\boldsymbol{\beta}_{\mu}}\,\boldsymbol{\chi}^{\boldsymbol{\beta}_{\mu}}dP_{\boldsymbol{\chi}}\\
\end{bmatrix}.\label{eq:Gnu}
\end{equation}
Owing to the assumption that the conditions given by~(\ref{eq:finitemoments}) and~(\ref{eq:nondegen}) are fulfilled, the Gram matrix~$\boldsymbol{G}$ is bounded and positive definite.
The Cholesky factorization~$\boldsymbol{G}=\boldsymbol{R}^{\mathrm{T}}\boldsymbol{R}$ of the Gram matrix thus provides a~$\mu$-dimensional, square, and upper triangular matrix $\boldsymbol{R}$ with strictly positive diagonal entries.
The inversion of this upper triangular matrix $\boldsymbol{R}$ in turn provides an upper triangular matrix~$\boldsymbol{R}^{-1}$.
Let~$s_{ij}$  denote the entries of~$\boldsymbol{R}^{-1}$; then, the orthonormal polynomials~$\{\Gamma_{\boldsymbol{\beta}},0\leq|\boldsymbol{\beta}|\leq q\}$ to be determined are obtained as follows:
\begin{equation}
\Gamma_{\boldsymbol{\beta}_{j}}(\boldsymbol{\chi})=s_{1j}\boldsymbol{\chi}^{\boldsymbol{\beta}_{1}}+s_{2j}\boldsymbol{\chi}^{\boldsymbol{\beta}_{2}}+\ldots+s_{jj}\boldsymbol{\chi}^{\boldsymbol{\beta}_{j}},\quad j=1,\ldots,\mu;\label{eq:chol}
\end{equation}
in fact, with~$\boldsymbol{\chi}_{\mu}=[\boldsymbol{\chi}^{\boldsymbol{\beta}_{1}} \ldots \boldsymbol{\chi}^{\boldsymbol{\beta}_{\mu}}]^{\mathrm{T}}$ and~$\boldsymbol{\Gamma}_{\mu}=[\Gamma_{\boldsymbol{\beta}_{1}} \ldots \Gamma_{\boldsymbol{\beta}_{\mu}}]^{\mathrm{T}}$, we immediately obtain
\begin{equation}
\int_{\real^{z}}\boldsymbol{\Gamma}_{\mu}(\boldsymbol{\chi})\boldsymbol{\Gamma}_{\mu}(\boldsymbol{\chi})^{\mathrm{T}}dP_{\boldsymbol{\chi}}=\boldsymbol{R}^{-\mathrm{T}}\int_{\real^{z}}\boldsymbol{\chi}_{\mu}\boldsymbol{\chi}_{\mu}^{\mathrm{T}}dP_{\boldsymbol{\chi}}\;\boldsymbol{R}^{-1}=\boldsymbol{R}^{-\mathrm{T}}\boldsymbol{G}\boldsymbol{R}^{-1}=\boldsymbol{I}.
\end{equation}
Here,~$\boldsymbol{I}$ denotes the~$\mu$-dimensional identity matrix.

It should be noted that the computation of the basis of multivariate polynomials by orthonormalization of the basis of monomials may be overly sensitive to numerical approximation and roundoff errors for large values of the dimension~$z$ and total degree~$q$.   
This issue can be addressed, for instance, by following the approaches given in~\citep{golub2009,giraud2005}.

\subsection{Quadrature rules with respect to arbitrary probability distributions}\label{sec:l1minimization}
Here, we provide details on the method for the computational construction of embedded quadrature rules described in Sec.~\ref{sec:sec4}.
At the core of this method is the solution of the $L^{1}$-minimization problem~(\ref{eq:L1}) for a sparse optimal solution.
Optimization theory provides two types of approaches for solving $L^{1}$-minimization problems.
Both types of approaches first involve the reformulation of the $L^{1}$-minimization problem as an equivalent linear-programming problem; then, the approaches of the first type apply a simplex algorithm, whereas the approaches of the second type apply an interior-point algorithm~\citep{wright1997,nocedal2006}.
The application of a simplex algorithm directly provides a sparse optimal solution because a simplex algorithm computes an optimal solution by performing a sequence of pivoting operations on basic feasible solutions.
In contrast, the application of an interior-point algorithm generally provides a fully populated optimal solution; thus, the application of an interior-point algorithm requires an additional computational procedure to extract a sparse optimal solution~\citep{wright1997,megiddo1991}.
In this section, we provide details on the solution of~(\ref{eq:L1}) using an interior-point algorithm.
It should be noted that our discussion can readily be extended to obtain a solution to~(\ref{eq:L1}) by a simplex algorithm; however, for brevity, we do not explicitly demonstrate this extension.

Adopting the notations used in Sec.~\ref{sec:sec4}, we consider the construction of a quadrature rule for integration with respect to the probability distribution~$P_{\boldsymbol{\chi}}$ on the finite-dimensional Euclidean space~$\real^{z}$ with accuracy level~$\lambda$. 
To ensure that the proposed method can be executed properly, we assume that the probability distribution~$P_{\boldsymbol{\chi}}$ satisfies the condition given by~(\ref{eq:finitemomentsb}).

Let us assume that a quadrature rule for integration with respect to~$P_{\boldsymbol{\chi}}$ is already available and this quadrature rule has a very large number (denoted by~$\tilde{\nu}$) of nodes and associated weights.
Then, we focus on the construction of an embedded quadrature rule with fewer nodes and weights using the solution of the~$L^{1}$-minimization problem~$\min_{\boldsymbol{\varpi}\in\real^{\tilde{\nu}}}\vectornorm{\boldsymbol{\varpi}}_{L^{1}}$ subject to~$\boldsymbol{A}\boldsymbol{\varpi}=\boldsymbol{b}$, in which~$\boldsymbol{A}$ is the~$\mu\times\tilde{\nu}$-dimensional matrix and~$\boldsymbol{b}$ the~$\mu$-dimensional vector defined by~(\ref{eq:A}) and~(\ref{eq:b}), respectively, introduced to impose the prescribed accuracy level~$\lambda$; it should be noted that~$\mu=\text{dim}(\mathcal{P}_{z}^{2\lambda-1})$ denotes the number of imposed equality constraints, as mentioned previously.
Following the approach given in~\citep{wright1997,nocedal2006}, we convert the~$L^{1}$-minimization problem into a linear program as follows:
\begin{equation}
\min_{\boldsymbol{\varpi},\boldsymbol{t}\in\real^{\tilde{\nu}}}\;\boldsymbol{e}^{\mathrm{T}}\boldsymbol{t},\quad\text{subject to~$\boldsymbol{A}\boldsymbol{\varpi}=\boldsymbol{b}$, $\boldsymbol{\varpi}-\boldsymbol{t}\leq\boldsymbol{0}$, and~$\boldsymbol{\varpi}+\boldsymbol{t}\geq\boldsymbol{0}$},\label{eq:Lp}
\end{equation}
where $\boldsymbol{t}$ is a~$\tilde{\nu}$-dimensional vector of \textit{slack variables} and~$\boldsymbol{e}$ is a~$\tilde{\nu}$-dimensional vector defined by~$\boldsymbol{e}=[1,\ldots,1]^{\mathrm{T}}$. 
The Karush-Kuhn-Tucker optimality conditions for the Lagrangian associated with the constrained optimization problem~(\ref{eq:Lp}) are then expressed as follows: 
\begin{align}
&\boldsymbol{A}^{\mathrm{T}}\boldsymbol{\lambda}-\boldsymbol{\mu}+\tilde{\boldsymbol{\mu}}=\boldsymbol{0}\quad\text{and}\quad\boldsymbol{e}-\boldsymbol{\mu}-\tilde{\boldsymbol{\mu}}=\boldsymbol{0}&&\text{(stationarity)},\\
&\boldsymbol{A}\boldsymbol{\varpi}-\boldsymbol{b}=\boldsymbol{0},\quad-\boldsymbol{\varpi}+\boldsymbol{t}\geq\boldsymbol{0},\quad\text{and}\quad\boldsymbol{\varpi}+\boldsymbol{t}\geq\boldsymbol{0}&&\text{(primal feasibility)},\\
&\boldsymbol{\mu}\geq\boldsymbol{0}\quad\text{and}\quad\tilde{\boldsymbol{\mu}}\geq\boldsymbol{0}&&\text{(dual feasibility)},\\
&\mu_{k}(-\varpi_{k}+t_{k})=0\quad\text{and}\quad\tilde{\mu}_{k}(\varpi_{k}+t_{k})=0,\quad1\leq k\leq\tilde{\nu}&&\text{(complementarity)},
\end{align}
where the~$\mu$-dimensional vector~$\boldsymbol{\lambda}$ and the~$\tilde{\nu}$-dimensional vectors~$\boldsymbol{\mu}$ and~$\tilde{\boldsymbol{\mu}}$ collect the Lagrange multipliers.
After eliminating~$\tilde{\boldsymbol{\mu}}=\boldsymbol{e}-\boldsymbol{\mu}$, these optimality conditions are expressed as follows:
\begin{align}
\label{eq:lagrange1}&\boldsymbol{A}^{\mathrm{T}}\boldsymbol{\lambda}-2\boldsymbol{\mu}+\boldsymbol{e}=\boldsymbol{0}&&\text{(stationarity)},\\
\label{eq:lagrange2}&\boldsymbol{A}\boldsymbol{\varpi}-\boldsymbol{b}=\boldsymbol{0},\quad-\boldsymbol{\varpi}+\boldsymbol{t}\geq\boldsymbol{0},\quad\text{and}\quad\boldsymbol{\varpi}+\boldsymbol{t}\geq\boldsymbol{0}&&\text{(primal feasibility)},\\
&\boldsymbol{\mu}\geq\boldsymbol{0}\quad\text{and}\quad\boldsymbol{e}-\boldsymbol{\mu}\geq\boldsymbol{0}&&\text{(dual feasibility)},\\
\label{eq:lagrange4}&\mu_{k}(-\varpi_{k}+t_{k})=0\quad\text{and}\quad(1-\mu_{k})(\varpi_{k}+t_{k})=0,\quad1\leq k\leq\tilde{\nu}&&\text{(complementarity)}.
\end{align}
The \textit{primal linear program}~(\ref{eq:Lp}) is associated with the \textit{dual linear program}
\begin{equation}
\min_{\boldsymbol{\lambda}\in\real^{\mu},\boldsymbol{\mu}\in\real^{\tilde{\nu}}}\;\boldsymbol{b}^{\mathrm{T}}\boldsymbol{\lambda},\quad\text{subject to~$\boldsymbol{A}^{\mathrm{T}}\boldsymbol{\lambda}-2\boldsymbol{\mu}+\boldsymbol{e}=\boldsymbol{0}$, $\boldsymbol{\mu}\geq\boldsymbol{0}$, and~$\boldsymbol{e}-\boldsymbol{\mu}\geq\boldsymbol{0}$};\label{eq:Lpdual}
\end{equation}
the relationship between the primal linear program~(\ref{eq:Lp}) and the dual linear program~(\ref{eq:Lpdual}) is that their Karush-Kuhn-Tucker optimality conditions are identical.
Consequently, using the complementarity~(\ref{eq:lagrange4}), which implies that $\boldsymbol{\mu}^{\mathrm{T}}(-\boldsymbol{\varpi}+\boldsymbol{t})+(\boldsymbol{e}-\boldsymbol{\mu})^{\mathrm{T}}(\boldsymbol{\varpi}+\boldsymbol{t})=0$, the stationarity~(\ref{eq:lagrange1}), which implies that $-\boldsymbol{\lambda}^{\mathrm{T}}\boldsymbol{A}\boldsymbol{\varpi}+\boldsymbol{e}^{\mathrm{T}}\boldsymbol{t}=0$, and the primal feasibility~(\ref{eq:lagrange2}), it can be shown that the optimal values for the primal and dual linear programs are the same:
\begin{equation}
\boldsymbol{\lambda}^{\mathrm{T}}\boldsymbol{b}=\boldsymbol{e}^{\mathrm{T}}\boldsymbol{t}.\label{eq:L1primaldual}
\end{equation}

Further, it should be noted that the complementarity~(\ref{eq:lagrange4}) implies that~$\mu_{k}=0$ when~$\varpi_{k}<0$ and that~$\mu_{k}=1$ when~$\varpi_{k}>0$; the values of the Lagrange multipliers~$\mu_{k}$ associated with the nonzero components~$\varpi_{k}$ are thus determined solely by the sign of these components.

Now, the solution of the linear program using a primal-dual interior-point algorithm yields a \textit{primal-dual optimal solution}~$(\boldsymbol{\varpi},\boldsymbol{t},\boldsymbol{\lambda},\boldsymbol{\mu})$ that solves the Karush-Kuhn-Tucker optimality conditions~(\ref{eq:lagrange1})--(\ref{eq:lagrange4}).
However, the vector~$\boldsymbol{\varpi}$ obtained as a part of such a quadruple, i.e.~$(\boldsymbol{\varpi},\boldsymbol{t},\boldsymbol{\lambda},\boldsymbol{\mu})$, is in general fully populated and therefore not sparse.

Following the approach given in~\citep{wright1997,nocedal2006}, we propose to use the procedure described next to extract a sparse optimal solution.   
Let~$\mathcal{B}\cup\mathcal{N}$ be a partitioning of the index set~$\{1,\ldots,\tilde{\nu}\}$ such that~$\mathcal{B}$ collects the indices associated with the nonzero components of~$\boldsymbol{\varpi}$ and~$\mathcal{N}$ collects the remaining indices associated with the vanishing components of~$\boldsymbol{\varpi}$.
Let the number of nonzero components be greater than the number (denoted by~$\mu$) of imposed equality constraints.
Then, there exists necessarily a vector~$\boldsymbol{z}\neq\boldsymbol{0}$ such that $\boldsymbol{A}\boldsymbol{z}=\boldsymbol{0}$ and the components of~$\boldsymbol{z}$ labeled by the indices in~$\mathcal{N}$ vanish; and there exists necessarily a scalar~$\alpha$ such that~$\boldsymbol{\varpi}+\alpha\boldsymbol{z}$ has more vanishing components than~$\boldsymbol{\varpi}$ and each component~$\varpi_{k}+\alpha z_{k}$ either vanishes or has the same sign as~$\varpi_{k}$.  
Because the components~$\varpi_{k}+\alpha z_{k}$ either vanish or have the same sign as the components~$\varpi_{k}$, the Lagrange multipliers~$\mu_{k}$ satisfy
\begin{equation}
\mu_{k}\Big(\hspace{-1mm}-(\varpi_{k}+\alpha z_{k})+| \varpi_{k}+\alpha z_{k}|\Big)=0\;\text{and}\;(1-\mu_{k})\Big((\varpi_{k}+\alpha z_{k})+| \varpi_{k}+\alpha z_{k}|\Big)=0,\;1\leq k\leq\tilde{\nu};
\end{equation}
hence, the quadruple~$(\boldsymbol{w}+\alpha\boldsymbol{z},|\boldsymbol{w}+\alpha\boldsymbol{z}|,\boldsymbol{\lambda},\boldsymbol{\mu})$ is also a primal-dual optimal solution that satisfies~(\ref{eq:lagrange1})--(\ref{eq:lagrange4}).
Consequently, the quadruple~$(\boldsymbol{w}+\alpha\boldsymbol{z},|\boldsymbol{w}+\alpha\boldsymbol{z}|,\boldsymbol{\lambda},\boldsymbol{\mu})$ satisfies $-\boldsymbol{\lambda}^{\mathrm{T}}\boldsymbol{A}(\boldsymbol{\varpi}+\alpha\boldsymbol{z})+\boldsymbol{e}^{\mathrm{T}}|\boldsymbol{\varpi}+\alpha\boldsymbol{z}|=0$ such that, with~$\boldsymbol{A}\boldsymbol{z}=\boldsymbol{0}$, the value taken by the objective function at~$(\boldsymbol{\varpi}+\alpha\boldsymbol{z},|\boldsymbol{\varpi}+\alpha\boldsymbol{z}|,\boldsymbol{\lambda},\boldsymbol{\mu})$ is equal to the value taken by the objective function at~$(\boldsymbol{\varpi},\boldsymbol{t},\boldsymbol{\lambda},\boldsymbol{\mu})$:
\begin{equation}
\boldsymbol{e}^{\mathrm{T}}|\boldsymbol{\varpi}+\alpha\boldsymbol{z}|=\boldsymbol{\lambda}^{\mathrm{T}}\boldsymbol{A}(\boldsymbol{\varpi}+\alpha\boldsymbol{z})=\boldsymbol{\lambda}^{\mathrm{T}}\boldsymbol{A}\boldsymbol{\varpi}=\boldsymbol{\lambda}^{\mathrm{T}}\boldsymbol{b}=\boldsymbol{e}^{\mathrm{T}}\boldsymbol{t}=\boldsymbol{e}^{\mathrm{T}}|\boldsymbol{\varpi}|.
\end{equation}
In conclusion, this procedure yields a vector~$\boldsymbol{\varpi}+\alpha\boldsymbol{z}$ that at least has one nonzero component less than~$\boldsymbol{\varpi}$ but is still optimal in that the value of the objective function remains unchanged.
Clearly, the repeated application of this procedure will yield an optimal solution that has at most as many nonzero components as there are imposed equality constraints. 
Algorithm~\ref{algo:algo4} outlines the method thus obtained for the computation of an embedded quadrature rule.

\begin{algorithm}[htp]
\SetKwInOut{Input}{Input}\SetKwInOut{Output}{Output}
\SetKw{KwAnd}{and}
\SetKwBlock{transformation}{Transformation of variables}{end}
\SetKwBlock{minimization}{$L^{1}$-minimization}{end}
\SetKwBlock{rounding}{Extraction of sparse optimal solution}{end}
\SetKwBlock{embedded}{Construction of embedded quadrature rule}{end}
\Input{Quadrature level $\lambda;$\\
Quadrature rule~$\big\{(\tilde{\boldsymbol{\chi}}_{k},\tilde{w}_{k}),\;1\leq k\leq\tilde{\nu}\big\}$ for integration with respect to~$P_{\boldsymbol{\chi}}$\;}
\BlankLine
\minimization{Select polynomial basis~$\pi_{1},\ldots,\pi_{\mu}$ for~$\mathcal{P}_{z}^{2\lambda-1}$ with~$\mu=\text{dim}(\mathcal{P}_{z}^{2\lambda-1})$\;
Construct matrix~$\boldsymbol{A}$ using~(\ref{eq:A}) as follows:
\begin{equation*}
\boldsymbol{A}=\begin{bmatrix}
\pi_{1}\big(\tilde{\boldsymbol{\chi}}_{1}\big) & \pi_{1}\big(\tilde{\boldsymbol{\chi}}_{2}\big) & \ldots & \pi_{1}\big(\tilde{\boldsymbol{\chi}}_{\tilde{\nu}-1}\big) & \pi_{1}\big(\tilde{\boldsymbol{\chi}}_{\tilde{\nu}}\big)\\
\vdots & \vdots & & \vdots & \vdots \\
\pi_{\mu}\big(\tilde{\boldsymbol{\chi}}_{1}\big) & \pi_{\mu}\big(\tilde{\boldsymbol{\chi}}_{2}\big) & \ldots & \pi_{\mu}\big(\tilde{\boldsymbol{\chi}}_{\tilde{\nu}-1}\big) & \pi_{\mu}\big(\tilde{\boldsymbol{\chi}}_{\tilde{\nu}}\big)
\end{bmatrix};
\end{equation*}
Construct vector~$\boldsymbol{b}$ using~(\ref{eq:b}) as follows:
\begin{equation*}
\boldsymbol{b}=\begin{bmatrix} \int_{\real^{z}}\pi_{1}(\boldsymbol{\chi})dP_{\boldsymbol{\chi}} & \ldots & \int_{\real^{z}}\pi_{\mu}(\boldsymbol{\chi})dP_{\boldsymbol{\chi}}\end{bmatrix}^{\mathrm{T}};
\end{equation*}
Use interior-point algorithm to solve linear program
\begin{equation*}
\min_{\boldsymbol{\lambda}\in\real^{\mu},\boldsymbol{\mu}\in\real^{\tilde{\nu}}}\begin{bmatrix}-\boldsymbol{b}&\boldsymbol{0}\end{bmatrix}\begin{bmatrix}\boldsymbol{\lambda}\\\boldsymbol{\mu}\end{bmatrix},\quad\text{subject to}\;\begin{bmatrix}\boldsymbol{A^{\mathrm{T}}} & -2\boldsymbol{I}\end{bmatrix}\begin{bmatrix}\boldsymbol{\lambda}\\\boldsymbol{e}\end{bmatrix}=-\boldsymbol{e}\;\text{and}\;\boldsymbol{0}\leq\boldsymbol{\mu}\leq\boldsymbol{e},
\end{equation*}
for primal-dual optimal solution~$(\boldsymbol{\varpi},\boldsymbol{t},\boldsymbol{\lambda},\boldsymbol{\mu})$\;}
\BlankLine
\rounding{\Repeat{$\alpha=0$}{Partition~$\{1,\ldots,\tilde{\nu}\}$ as~$\mathcal{B}\cup\mathcal{N}$ such that~$\varpi_{k}=0$ if~$k$ is in $\mathcal{N}$\;
Find vector $\boldsymbol{z}$ such that $\boldsymbol{A}\boldsymbol{z}=\boldsymbol{0}$ and $z_{k}=0$ if~$k$ is in~$\mathcal{N}$\;
Set~$\alpha=\min\{|z_{k}/\varpi_{k}|:k\in\mathcal{B}\;\text{, $z_{k}\neq 0$, and sign}(z_{k})\neq\text{sign}(\varpi_{k})\}$\;
Update $\boldsymbol{\varpi}$ to $\boldsymbol{\varpi}+\alpha\boldsymbol{z}$\;}}
\BlankLine
\embedded{Synthesize quadrature rule~$\big\{\big(\boldsymbol{\chi}_{\ell}=\tilde{\boldsymbol{\chi}}_{k_{\ell}},w_{\ell}=\varpi_{k_{\ell}}\big),1\leq\ell\leq\nu\big\}$,\\ in which $\mathcal{B}=\{k_{1},\ldots,k_{\nu}\}$ with necessarily $\nu\leq\mu=\text{dim}(\mathcal{P}_{z}^{2\lambda-1})$\;}
\caption{Computation of an embedded quadrature rule.}\label{algo:algo4}
\end{algorithm}

\section{Realization for a multiphysics problem}\label{sec:sec6}

\subsection{Problem formulation}
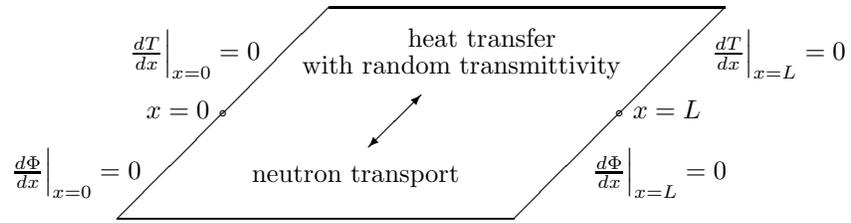
\begin{figure}[htp]
  \begin{center}
    \begin{picture}(250,100)(0,0)
      \put(50,50){\line(1,1){40}}
      \put(50,50){\line(-1,-1){40}}
      \put(200,50){\line(1,1){40}}
      \put(200,50){\line(-1,-1){40}}
      \put(10,10){\line(1,0){150}}
      \put(90,90){\line(1,0){150}}
      \put(120,75){\makebox{heat transfer}}
      \put(50,24){\makebox{$\;\;\;\;$neutron transport}}
      \put(80,65){\makebox{with random transmittivity}}
      \put(115,47){\vector(1,1){10}}
      \put(115,47){\vector(-1,-1){10}}
      \put(50,50){\circle{2}}
      \put(200,50){\circle{2}}
      \put(21,48){\makebox{$x=0$}}
      \put(205,48){\makebox{$x=L$}}
      \put(15,70){\makebox{$\frac{dT}{dx}\Big|_{x=0}=0$}}
      \put(235,70){\makebox{$\frac{dT}{dx}\Big|_{x=L}=0$}}
      \put(-30,25){\makebox{$\frac{d\Phi}{dx}\Big|_{x=0}=0$}}
      \put(190,25){\makebox{$\frac{d\Phi}{dx}\Big|_{x=L}=0$}}
    \end{picture}
  \end{center}
  \caption{Schematic representation of the problem.}\label{fig:figure0}
\end{figure}
We consider the stationary transport of neutrons in a one-dimensional reactor with temperature feedback~\citep{lamarsh2002}.
Let the reactor occupy an open interval~$]0,L[$ (Fig.~\ref{fig:figure0}).
The problem then involves finding the temperature~$T$ and neutron flux~$\Phi$ such that
\begin{equation}
\label{eq:neutron1c}
\begin{aligned}
&\frac{d}{dx}\left(k\frac{dT}{dx}\right)-h(T-T_{\infty})=-E_{\text{f}}\Sigma_{\text{f}}(T)\Phi,\\
&\frac{d}{dx}\left(D(T)\frac{d\Phi}{dx}\right)-\Big(\Sigma_{\text{a}}(T)-\nu\Sigma_{\text{f}}(T)\Big)\Phi=-s,
\end{aligned}
\end{equation}
under homogeneous Neumann boundary conditions.
The first term on the left-hand side of the heat subproblem represents heat conduction, and the second term represents the transmission of heat to the surroundings; further, the right-hand side represents a distributed heat source proportional to the neutron flux.  
The first term on the left-hand side of the neutronics subproblem represents neutron diffusion, and the second term represents the net effect of the absorption and generation of neutrons; further, the right-hand side represents a distributed neutron source.  
The coefficients~$k$ and~$h$ are the heat conductivity and heat transmittivity, respectively; the temperature $T_{\infty}$ is the ambient temperature; and $\nu$ and~$E_{\text{f}}$ are the number of neutrons and the energy released per fission reaction, respectively.
The coefficients~$D$, $\Sigma_{\text{a}}$, and~$\Sigma_{\text{f}}$ are the neutron diffusion constant, fission cross section, and absorption cross section, respectively; these coefficients depend on the reactor temperature as follows:
\begin{equation}
D\big(T(x)\big)=D_{\text{ref}}\sqrt{\frac{T(x)}{T_{\text{ref}}}},\quad\Sigma_{\text{a}}\big(T(x)\big)=\Sigma_{\text{a,ref}}\sqrt{\frac{T_{\text{ref}}}{T(x)}},\quad\Sigma_{\text{f}}\big(T(x)\big)=\Sigma_{\text{f,ref}}\sqrt{\frac{T_{\text{ref}}}{T(x)}}.\label{eq:couplingmechanism}
\end{equation}

\subsection{Deterministic weak formulation}
Let~$H=H^{1}(]0,L[)$ be the space of functions that are sufficiently regular to describe the solutions to the heat and neutronics subproblems.
The weak formulation then involves finding~$T$ and~$\Phi$ in~$H$ such that
\begin{equation}
\label{eq:neutron2}\begin{aligned}
&\int_{0}^{L}k\frac{dT}{dx}\frac{dS}{dx}dx+\int_{0}^{L}h(T-T_{\infty})Sdx=\int_{0}^{L}E_{\text{f}}\Sigma_{\text{f}}(T)\Phi Sdx,&&\forall S\in H,\\
&\int_{0}^{L}D(T)\frac{d\Phi}{dx}\frac{d\Psi}{dx}dx+\int_{0}^{L}\Big(\Sigma_{\text{a}}(T)-\nu\Sigma_{\text{f}}(T)\Big)\Phi\Psi dx=\int_{0}^{L}s\Psi dx,&&\forall \Psi\in H.\\
\end{aligned}
\end{equation}

\subsection{Random thermal transmittivity}
We incorporate uncertainties by modeling the thermal transmittivity as a random field~$\{h(x,\cdot),1\leq x\leq L\}$ such that
\begin{equation}
h(x,\boldsymbol{\xi})=\overline{h}\bigg(1+\delta\sum_{j=1}^{m}\sqrt{\lambda_{j}}\sqrt{3}\xi_{j}\phi^{j}(x)\bigg),\label{eq:hN}
\end{equation}
where the~$\xi_{j}$ are statistically independent uniform random variables defined on a probability triple~$(\Theta,\mathcal{T},P)$ with values in $[-1,1]$ and the~$\sqrt{3}\xi_{j}$ are thus uniform random variables with unit standard deviation; further, the~$\lambda_{j}$ and~$\phi^{j}$ are the eigenvalues and eigenmodes, respectively, of the eigenproblem~$\mathcal{C}(\phi^{j})=\lambda_{j}\phi^{j}$, where~$\mathcal{C}$ is the covariance integral operator with
\begin{equation}
C(x,y)=\frac{4a^{2}}{\pi^{2}(x-y)^{2}}\sin^{2}\left(\frac{\pi(x-y)}{2a}\right)\label{eq:covar}
\end{equation}
as the kernel; here, the parameter~$a$ is the spatial correlation length of~$\{h(x,\cdot),1\leq x\leq L\}$.
Clearly, the random field~$\{h(x,\cdot),1\leq x\leq L\}$ thus obtained is such that the random variable~$h(x,\cdot)$ has the mean~$\overline{h}$ and coefficient of variation~$\delta$ at every position~$x$, at least when the approximation error introduced owing to the truncation of the expansion after $m$ terms is not taken into account.  

\subsection{Stochastic weak formulation} 
The weak formulation of the stochastic problem involves finding random variables~$T$ and~$\Phi$ defined on~$(\Theta,\mathcal{T},P)$, with values in~$H$, such that
\begin{equation}
\label{eq:neutron2b}
\begin{aligned}
&\int_{0}^{L}k\frac{dT}{dx}\frac{dS}{dx}dx+\int_{0}^{L}h(\boldsymbol{\xi})(T-T_{\infty})Sdx=\int_{0}^{L}E_{\text{f}}\Sigma_{\text{f}}(T)\Phi Sdx,&&\forall S\in H,\\
&\int_{0}^{L}D(T)\frac{d\Phi}{dx}\frac{d\Psi}{dx}dx+\int_{0}^{L}\Big(\Sigma_{\text{a}}(T)-\nu\Sigma_{\text{f}}(T)\Big)\Phi\Psi dx=\int_{0}^{L}s\Psi dx,&&\forall \Psi\in H.
\end{aligned}
\end{equation}
 
\subsection{Discretization of space}
The finite element~(FE) method is used for the discretization of space.
The domain~$[0,L]$ is meshed using $r-1$ elements of equal length.  
Let~$N_{1},\ldots,N_{r}$ then be a basis of element-wise linear shape functions such that~$N_{j}$ takes value~$1$ at the~$j$-th node and~0 at other nodes. 
Using this basis, the random temperature~$T$ and neutron flux~$\Phi$ are approximated as follows:
\begin{equation}
\label{eq:Vh1}
\begin{aligned}
T^{r}(x)=\sum_{j=1}^{r}T_{j}N_{j}(x),\quad\quad\quad T_{j}\in\real,\\
\Phi^{r}(x)=\sum_{j=1}^{r}\Phi_{j}N_{j}(x),\quad\quad\quad \Phi_{j}\in\real.\\
\end{aligned}
\end{equation} 
The FE discretization of the stochastic weak formulation~(\ref{eq:neutron2b}) then involves finding random vectors~$\boldsymbol{T}=\{T_{1},\ldots,T_{r}\}$ and~$\boldsymbol{\Phi}=\{\Phi_{1},\ldots,\Phi_{r}\}$ defined on~$(\Theta,\mathcal{T},P)$, with values in~$\real^{r}$, which collect the nodal values of the random temperature and neutron flux such that
\begin{equation}
\label{eq:discrfE3b}
\begin{aligned}
&[\boldsymbol{K}+\boldsymbol{H}(\boldsymbol{\xi})]\boldsymbol{T}=\boldsymbol{q}(\boldsymbol{\Phi},\boldsymbol{T}),\\
&[\boldsymbol{D}(\boldsymbol{T})+\boldsymbol{M}(\boldsymbol{T})]\boldsymbol{\Phi}=\boldsymbol{s}.
\end{aligned}
\end{equation}
Here, $\boldsymbol{K}$, $\boldsymbol{H}$, $\boldsymbol{D}(\boldsymbol{T})$, and $\boldsymbol{M}(\boldsymbol{T})$ are $r$-dimensional matrices, and $\boldsymbol{q}(\boldsymbol{\Phi},\boldsymbol{T})$ and $\boldsymbol{s}$ are $r$-dimensional vectors such that
\begin{align}
\boldsymbol{S}_{1}^{\mathrm{T}}\boldsymbol{K}\boldsymbol{S}_{2}&=\int_{0}^{L}k\frac{dS^{r}_{1}}{dx}\frac{dS^{r}_{2}}{dx}dx,\\
\boldsymbol{S}_{1}^{\mathrm{T}}\boldsymbol{H}\boldsymbol{S}_{2}&=\int_{0}^{L}hS^{r}_{1}S^{r}_{2}dx,\\
\boldsymbol{\Psi}_{1}^{\mathrm{T}}\boldsymbol{D}(\boldsymbol{T})\boldsymbol{\Psi}_{2}&=\int_{0}^{L}D\big(T^{r}\big)\frac{d\Psi^{r}_{1}}{dx}\frac{d\Psi^{r}_{2}}{dx}dx,\\
\boldsymbol{\Psi}_{1}^{\mathrm{T}}\boldsymbol{M}(\boldsymbol{T})\boldsymbol{\Psi}_{2}&=\int_{0}^{L}\Big(\Sigma_{\text{a}}\big(T^{r}\big)-\nu\Sigma_{\text{f}}\big(T^{r}\big)\Big)\Psi^{r}_{1}\Psi^{r}_{2}dx,\\
\boldsymbol{S}^{\mathrm{T}}\boldsymbol{q}(\boldsymbol{T},\boldsymbol{\Phi})&=\int_{0}^{L}E_{\text{f}}\Sigma_{\text{f}}\big(T^{r}\big)\Phi^{r}S^{r}dx+\int_{0}^{L}hT_{\infty}S^{r}dx,\\
\boldsymbol{S}^{\mathrm{T}}\boldsymbol{s}&=\int_{0}^{L}s\Psi^{r}dx.
\end{align}

\subsection{Reformulation as a realization of the model problem}
The aforementioned illustration problem can be reformulated as a particular realization of the general model problem introduced in Sec.~\ref{sec:sec2} as follows:
\begin{equation}
\begin{aligned}
&\boldsymbol{T}=\boldsymbol{a}(\boldsymbol{T},\boldsymbol{\Phi},\boldsymbol{\xi}),&&\quad\quad\quad\boldsymbol{a}:\real^{r}\times\real^{r}\times\real^{m}\rightarrow\real^{r},\\
&\boldsymbol{\Phi}=\boldsymbol{b}(\boldsymbol{T}),&&\quad\quad\quad\boldsymbol{b}:\real^{r}\rightarrow\real^{r},
\end{aligned}
\end{equation}
where $\boldsymbol{a}(\boldsymbol{T},\boldsymbol{\Phi},\boldsymbol{\xi})=[\boldsymbol{K}+\boldsymbol{H}(\boldsymbol{\xi})]^{-1}\boldsymbol{q}(\boldsymbol{\Phi},\boldsymbol{T})$ and $\boldsymbol{b}(\boldsymbol{T})=[\boldsymbol{D}(\boldsymbol{T})+\boldsymbol{M}(\boldsymbol{T})]^{-1}\boldsymbol{s}$.
This reformulation indicates that the illustration problem is a simplified realization of the model problem, for three reasons.
First, the data of the neutronics subproblem are not affected by their own sources of uncertainty~$\boldsymbol{\zeta}$.
Second, the neutronics subproblem admits a direct solution that does not require iteration.
Lastly, the neutronics and heat subproblems are coupled directly through their solution variables rather than through intermediate coupling variables.

\subsection{Dimension reduction} \label{sec:sec56}
Now, we will demonstrate the proposed methodology by approximating the random temperature by a truncated KL decomposition as it is communicated from the heat to the neutronics subproblem.
At iteration~$\ell$, let the random temperature be represented by a PC expansion as follows:
\begin{equation}
\widehat{\boldsymbol{T}}{}^{\ell,p}=\sum_{|\boldsymbol{\alpha}|=0}^{p}\widehat{\boldsymbol{T}}{}^{\ell}_{\boldsymbol{\alpha}}\psi_{\boldsymbol{\alpha}}(\boldsymbol{\xi}),\quad\quad\quad\widehat{\boldsymbol{T}}{}^{\ell}_{\boldsymbol{\alpha}}\in\real^{r}.\label{eq:upceklb}
\end{equation}
The mean and covariance of~$\widehat{\boldsymbol{T}}{}^{\ell,p}$ are then given by
\begin{align}
\overline{\boldsymbol{T}}{}^{\ell}&=\widehat{\boldsymbol{T}}{}^{\ell}_{\boldsymbol{0}},\\
\boldsymbol{C}_{\boldsymbol{T}}^{\ell}&=\sum_{|\boldsymbol{\alpha}|=1}^{p}\widehat{\boldsymbol{T}}{}^{\ell}_{\boldsymbol{\alpha}}\big(\widehat{\boldsymbol{T}}{}_{\boldsymbol{\alpha}}^{\ell}\big)^{\mathrm{T}}.
\end{align}   
Further, let the $r$-dimensional square matrix~$\boldsymbol{W}$ be the Gram matrix of the FE basis, i.e.,
\begin{equation}
\boldsymbol{W}=\begin{bmatrix}
\langle N_{1},N_{1}\rangle_{H} & \ldots & \langle N_{1},N_{r}\rangle_{H}\\
\vdots & & \vdots\\
\langle N_{r},N_{1}\rangle_{H} & \ldots & \langle N_{r},N_{r}\rangle_{H}
\end{bmatrix},
\end{equation}
where the inner product~$\langle\cdot,\cdot\rangle_{H}$ is such that~$\langle S_{1},S_{2}\rangle_{H}=\int_{0}^{L}S_{1}S_{2}dx+\int_{0}^{L}(dS_{1}/dx)(dS_{2}/dx)dx$ for any pair~$S_{1}$ and~$S_{2}$ of functions in~$H$. 
The solution of the generalized eigenproblem $\boldsymbol{W}^{\mathrm{T}}\boldsymbol{C}_{\boldsymbol{T}}^{\ell}\boldsymbol{W}\boldsymbol{\phi}^{j,\ell}=\lambda_{j}^{\ell}\boldsymbol{W}\boldsymbol{\phi}^{j,\ell}$ then provides the eigenvalues~$\lambda_{j}^{\ell}$ and the associated eigenmodes~$\boldsymbol{\phi}^{j,\ell}$ required to construct a reduced-dimensional representation~$\widehat{\boldsymbol{T}}{}^{\ell,p,d}$ of $\widehat{\boldsymbol{T}}{}^{\ell,p}$ as follows:
\begin{equation}
\widehat{\boldsymbol{T}}{}^{\ell,p,d}=\overline{\boldsymbol{T}}{}^{\ell}+\sum_{j=1}^{d}\sqrt{\lambda_{j}^{\ell}}\eta{}_{j}^{\ell,p}\boldsymbol{\phi}^{j,\ell},\label{eq:Tpdkl}
\end{equation}
where the~$\eta_{j}^{\ell,p}$ are random variables defined on~$(\Theta,\mathcal{T},P)$, with values in~$\real$, such that
\begin{equation}
\eta{}_{j}^{\ell,p}=\frac{1}{\sqrt{\lambda_{j}^{\ell}}}\big(\widehat{\boldsymbol{T}}{}^{\ell,p}-\overline{\boldsymbol{T}}{}^{\ell}\big)^{\mathrm{T}}\boldsymbol{W}\boldsymbol{\phi}^{j,\ell}.\label{eq:upcekletab}
\end{equation}
By substituting~(\ref{eq:upceklb}) in~(\ref{eq:upcekletab}), a representation of the~$\eta^{\ell,p}_{i}$ as a PC expansion is obtained: 
\begin{equation}
\boldsymbol{\eta}^{\ell,p}=\sum_{|\boldsymbol{\alpha}|=1}^{p}\boldsymbol{\eta}^{\ell}_{\boldsymbol{\alpha}}\psi_{\boldsymbol{\alpha}}(\boldsymbol{\xi})\quad\text{with}\quad\boldsymbol{\eta}_{\boldsymbol{\alpha}}^{\ell}=\begin{bmatrix}\frac{1}{\sqrt{\lambda_{1}^{\ell}}}\big(\widehat{\boldsymbol{T}}{}_{\boldsymbol{\alpha}}^{\ell}\big)^{\mathrm{T}}\boldsymbol{W}\boldsymbol{\phi}^{1,\ell} & \ldots & \frac{1}{\sqrt{\lambda_{d}^{\ell}}}\big(\widehat{\boldsymbol{T}}{}_{\boldsymbol{\alpha}}^{\ell}\big)^{\mathrm{T}}\boldsymbol{W}\boldsymbol{\phi}^{d,\ell}\end{bmatrix}^{\mathrm{T}},\label{eq:upcervetaaaa}
\end{equation}
thus completely characterizing the reduced random variables as a PC expansion.

It should be noted that the random neutron flux, in principle, could also be reduced as it passes from the neutronics subproblem to the heat subproblem.
However, because the data of the neutronics subproblem are not affected by their own sources of uncertainty, a reduction of the random neutron flux would not lower the number of sources of uncertainty that enter the heat subproblem and thus would not lead to a solution of the heat subproblem in a reduced-dimensional space.
This extension is therefore not demonstrated.

\subsection{Measure transformation}
Whereas the random variables~$\boldsymbol{\xi}=(\xi_{1},\ldots,\xi_{m})$ necessarily constitute the sources of uncertainty that enter the heat subproblem, the reduced random variables~$\boldsymbol{\eta}^{\ell}=(\eta_{1}^{\ell},\ldots,\eta_{d}^{\ell})$ of the KL decomposition of the random temperature can be construed as the sources of uncertainty that enter the neutronics subproblem.
Then, the proposed methodology leads to the approximation of the random temperature and neutron flux by PC expansions as follows:
\begin{equation}
\label{eq:PCTphi}\begin{aligned}
\widehat{\boldsymbol{T}}{}^{\ell,p}=\sum_{|\boldsymbol{\alpha}|=0}^{p}\widehat{\boldsymbol{T}}{}^{\ell}_{\boldsymbol{\alpha}}\psi_{\boldsymbol{\alpha}}(\boldsymbol{\xi}),&&\quad\quad\quad\widehat{\boldsymbol{T}}{}^{\ell}_{\boldsymbol{\alpha}}\in\real^{r},\\
\widehat{\boldsymbol{\Phi}}{}^{\ell,q}=\sum_{|\boldsymbol{\beta}|=0}^{q}\widehat{\boldsymbol{\Phi}}{}^{\ell}_{\boldsymbol{\beta}}\Gamma_{\boldsymbol{\beta}}^{\ell}(\boldsymbol{\eta}^{\ell,p}),&&\quad\quad\quad\widehat{\boldsymbol{\Phi}}{}^{\ell}_{\boldsymbol{\beta}}\in\real^{r};\\
\end{aligned}
\end{equation}
i.e., we obtain the approximation of the random temperature by a PC expansion in the input random variables and the approximation of the random neutron flux by a PC expansion in the reduced random variables.

We select the polynomials $\{\psi_{\boldsymbol{\alpha}},0\leq|\boldsymbol{\alpha}|\leq p\}$ as normalized Legendre polynomials, and we construct the polynomials $\{\Gamma^{\ell}_{\boldsymbol{\beta}},0\leq|\boldsymbol{\beta}|\leq q\}$ at each iteration using the method given in Secs.~\ref{sec:sec3} and~\ref{sec:orthopol}.
Further, we select the quadrature rule  for integration with respect to the probability distribution of the input random variables, denoted as $\{(\boldsymbol{\xi}_{k},v_{k}),\;1\leq k\leq N\}$, as a sparse-grid Gauss-Legendre quadrature rule of level $p+1$, and we construct the quadrature rule for integration with respect to the probability distribution of the reduced random variables, denoted as $\{(\boldsymbol{\eta}_{k}^{\ell},w_{k}^{\ell}),\;1\leq k\leq \nu^{\ell}\}$, using the method given in Secs.~\ref{sec:sec4} and~\ref{sec:l1minimization} for level $\lambda=q+1$.   

It should be noted that because the polynomial transformation~(\ref{eq:upcervetaaaa}) is continuous and the probability distribution of the input random variables has a closed and bounded support~$[-1,1]^{m}$, the probability distribution of the reduced random variables also has a closed and bounded support.
With reference to Secs.~\ref{sec:hilbert} and~\ref{sec:convergenceee}, this property suffices to ensure the convergence of the PC expansion of the random neutron flux as the total degree~$q$ increases.

\subsection{Implementation}\label{sec:illusimp}
The abovementioned computational construction of the polynomials and quadrature rules requires that integrals be computed with respect to the probability distribution of the reduced random variables to obtain the requisite Gram matrix and moment vector.
Because the probability distribution of the reduced random variables is characterized by the KL decomposition~(\ref{eq:Tpdkl}) as a transformation of the probability distribution of the input random variables through the PC expansion~(\ref{eq:upcervetaaaa}), 
these integrals can be evaluated easily using a quadrature rule for integration with respect to the probability distribution of the input random variables after a ``change of variables."

Owing to the PC expansion~(\ref{eq:upcervetaaaa}), the integral of a function~$f$ from~$\real^{d}$ into~$\real$ with respect to the probability distribution~$P_{\boldsymbol{\eta}}^{\ell}$ of the reduced random variables can be reformulated as an integral with respect to the probability distribution~$P_{\boldsymbol{\xi}}$ by a ``change of variables" as follows:
\begin{equation}
\int_{\real^{d}}f(\boldsymbol{\eta})dP_{\boldsymbol{\eta}}^{\ell}=\int_{\real^{m}}f\Big(\boldsymbol{\eta}^{\ell,p}(\boldsymbol{\xi})\Big)dP_{\boldsymbol{\xi}},\label{eq:changeofvariables}
\end{equation}
provided that either integral is well defined~\citep{dudley2004}; thus, the numerical evaluation of this integral can be performed using a quadrature rule for integration with respect to~$P_{\boldsymbol{\xi}}$. 

In this work, we use the Monte Carlo method to numerically evaluate the entries of the Gram matrix and the components of the moment vector as follows:
\begin{align}
\label{eq:Gij}G_{ij}^{\ell}&=\int_{\Theta}\big(\boldsymbol{\eta}^{\ell,p}\big)^{\boldsymbol{\beta}_{i}}\big(\boldsymbol{\eta}^{\ell,p}\big)^{\boldsymbol{\beta}_{j}}dP\approx\frac{1}{MC}\sum_{k=1}^{MC}\big(\boldsymbol{\eta}^{\ell,p}(\boldsymbol{\xi}_{k})\big)^{\boldsymbol{\beta}_{i}}\big(\boldsymbol{\eta}^{\ell,p}(\boldsymbol{\xi}_{k})\big)^{\boldsymbol{\beta}_{j}},\\
\label{eq:bj}b_{j}^{\ell}&=\int_{\Theta}\big(\boldsymbol{\eta}^{\ell,p}\big)^{\boldsymbol{\beta}_{j}}dP\approx\frac{1}{MC}\sum_{k=1}^{MC}\big(\boldsymbol{\eta}^{\ell,p}(\boldsymbol{\xi}_{k})\big)^{\boldsymbol{\beta}_{j}},
\end{align}
where~$\{\boldsymbol{\xi}_{k},1\leq k\leq MC\}$ collects $MC$ independent samples of the input random variables. 

In addition, the method given in Secs.~\ref{sec:sec4} and~\ref{sec:l1minimization} requires a quadrature rule to provide a collection of candidate nodes from which the nodes of the embedded quadrature rule can be selected.
We adopt the following construction in this work.
First, we construct a sparse-grid Gauss-Legendre quadrature rule, denoted as~$\{(\tilde{\boldsymbol{\xi}}_{k},\tilde{w}_{k}),1\leq k\leq\tilde{\nu}\}$, for integration with respect to the probability distribution of the input random variables.
Next, we carry out a ``change of variables" to transform this sparse-grid Gauss-Legendre quadrature rule to a corresponding quadrature rule of the form~$\{(\tilde{\boldsymbol{\eta}}{}^{\ell}_{k},\tilde{w}_{k}),1\leq k\leq\tilde{\nu}\}$ for integration with respect to the probability distribution of the reduced random variables by choosing the nodes as $\tilde{\boldsymbol{\eta}}{}_{k}^{\ell}=\boldsymbol{\eta}^{\ell,p}(\tilde{\boldsymbol{\xi}}_{k})$ and keeping the weights unchanged.
Finally, from this quadrature rule, we construct an embedded quadrature rule~$\{(\boldsymbol{\eta}_{k}^{\ell},w_{k}^{\ell}),\;1\leq k\leq \nu\}$ of the desired level $\lambda$ by applying Algorithm~1.

The aforementioned sparse-grid Gauss-Legendre quadrature rule is parameterized itself by its own level, denoted here as~$\tilde{\lambda}$.
As a larger~$\tilde{\lambda}$ is chosen, the number of nodes that the sparse-grid Gauss-Legendre quadrature rule has increases; 
thus, an increase in $\tilde{\lambda}$ provides Algorithm~1 with a greater choice of candidate nodes to select from, and therefore, Algorithm~1 can be expected to yield a better embedded quadrature rule that provides a smaller value for the sum of the absolute values of the weights.
Because this sum has a lower bound of 1, the adequate level $\tilde{\lambda}$ can readily be selected by monitoring the convergence of this sum with respect to $\tilde{\lambda}$.

This is not the only implementation available.
The Monte Carlo method can be replaced by a fully tensorized or sparse-grid quadrature rule to evaluate~(\ref{eq:Gij}) and~(\ref{eq:bj}).
Further, the Gram matrix and moment vector can be deduced using polynomial algebra from the PC expansion~(\ref{eq:upcervetaaaa}) and the moments of $P_{\boldsymbol{\xi}}$. 
Finally, the sparse-grid Gauss-Legendre quadrature rule can be replaced by Clenshaw-Curtis, Monte Carlo, or other rules.  

\subsection{Selection of the reduced dimension and the polynomial degree}
At each iteration, we select  the number of terms retained in~(\ref{eq:Tpdkl}) by the KL decomposition~$\widehat{\boldsymbol{T}}{}^{\ell,p,d}$ of~$\widehat{\boldsymbol{T}}{}^{\ell,p}$ as the smallest dimension $d$ that satisfies the following condition:
\begin{equation}
\sqrt{\int_{\Theta}\vectornorm{\widehat{\boldsymbol{T}}{}^{\ell,p}-\widehat{\boldsymbol{T}}{}^{\ell,p,d}}^{2}_{\boldsymbol{W}}dP}\leq\epsilon_{1}\sqrt{\int_{\Theta}\vectornorm{\widehat{\boldsymbol{T}}{}^{\ell,p}}^{2}_{\boldsymbol{W}}dP},\quad\quad\quad\forall\ell\in\integer,\label{eq:criterion1}
\end{equation}
where~$\epsilon_{1}$ is a prescribed tolerance level.
Further, at each iteration, we truncate the PC expansion~$\widehat{\boldsymbol{\Phi}}{}^{\ell,p}$ in (\ref{eq:PCTphi}) at the smallest total degree $q$ that satisfies the following condition:
\begin{equation}
\sqrt{\int_{\Theta}\vectornorm{\widehat{\boldsymbol{\Phi}}{}^{\ell,q}-\widehat{\boldsymbol{\Phi}}{}^{\ell,q-1}}^{2}_{\boldsymbol{W}}dP}\leq\epsilon_{2}\sqrt{\int_{\Theta}\vectornorm{\widehat{\boldsymbol{\Phi}}{}^{\ell,q}}^{2}_{\boldsymbol{W}}dP},\quad\quad\quad\forall\ell\in\integer,\label{eq:criterion2}
\end{equation}
where~$\epsilon_{2}$ is a prescribed tolerance level.
Clearly, these criteria may result in the dependence of the reduced dimension~$d$ and the total degree~$q$ on the iteration~$\ell$.

\subsection{Concluding remarks}
Algorithm~\ref{algo:algo6} summarizes the implementation of the illustration problem presented in this section.
The key feature of this implementation is that it enables a solution of the neutronics subproblem in a reduced-dimensional space when the KL decomposition can extract a low-dimensional representation of the random temperature ($d<m$),
while maintaining accuracy.

The solution of the neutronics subproblem in a reduced-dimensional space can be expected to reduce the number of terms that are required in the PC expansion of the random neutron flux to achieve sufficient accuracy.
Further, the solution of the neutronics subproblem in a reduced-dimensional space can be expected to reduce the number of quadrature nodes that are required for the nonintrusive projection method to achieve sufficient accuracy in the PC coordinates of the random neutron flux and therefore to reduce the number of times a sample of the neutronics subproblem must be solved, thus lowering the computational cost.

\begin{algorithm}[htp]
\SetKwInOut{Input}{Input}\SetKwInOut{Output}{Output}
\SetKw{KwAnd}{and}
\SetKwBlock{FirstMonoProblem}{neutronics subproblem}{end}
\SetKwBlock{SecondMonoProblem}{heat subproblem}{end}
\SetKwBlock{NexusRegion}{dimension reduction}{end}
\SetKwBlock{NexusRegionn}{measure transformation}{end}
\SetKwBlock{Initialization}{initialization}{end}
\Input{Error tolerance levels $\small{\epsilon_{1}}\normalsize$ and $\small{\epsilon_{2}}\normalsize$;\\
PC basis $\small{\big\{\psi_{\boldsymbol{\alpha}},\;0\leq|\boldsymbol{\alpha}|\leq p\big\}}\normalsize$ up to total degree $\small{p}\normalsize$ w.r.t. $\small{P_{\boldsymbol{\xi}}}\normalsize$;\\
Quadrature rule $\small{\{(\boldsymbol{\xi}_{k},v_{k}),\;1\leq k\leq N\}}\normalsize$ of level $\small{p+1}\normalsize$ w.r.t. $\small{P_{\boldsymbol{\xi}}}\normalsize$\;}
\vspace{-1mm}
$\small{\ell=1}\normalsize$\;
\Repeat{$($convergence$)$}{
\SecondMonoProblem{
\For{$k=1$ \KwTo $\,N_{p+1}$}{
Solve $\small{\big[\boldsymbol{K}+\boldsymbol{H}(\boldsymbol{\xi}_{k})\big]\widehat{\boldsymbol{T}}{}^{\ell}\big(\boldsymbol{\xi}_{k}\big)=\boldsymbol{q}\big(\widehat{\boldsymbol{T}}{}^{\ell-1,p}(\boldsymbol{\xi}_{k}),\widehat{\boldsymbol{\Phi}}{}^{\ell-1,q}(\boldsymbol{\xi}_{k})\big),}\normalsize$\\
with $\small{\widehat{\boldsymbol{\Phi}}{}^{\ell-1,q}(\boldsymbol{\xi}_{k})=\sum_{|\boldsymbol{\beta}|=0}^{q}\widehat{\boldsymbol{\Phi}}{}^{\ell-1}_{\boldsymbol{\beta}}\Gamma_{\boldsymbol{\beta}}^{\ell-1}\Big(\sum_{|\boldsymbol{\alpha}|=1}^{p}\boldsymbol{\eta}^{\ell-1}_{\boldsymbol{\alpha}}\psi_{\boldsymbol{\alpha}}(\boldsymbol{\xi}_{k})\Big)}\normalsize$\;}
Compute PC coordinates of $\small{\widehat{\boldsymbol{T}}{}^{\ell,p}}\normalsize$ using $\small{\widehat{\boldsymbol{T}}{}^{\ell}_{\boldsymbol{\alpha}}=\sum_{k=1}^{N_{p=1}}\widehat{\boldsymbol{T}}{}^{\ell}\big(\boldsymbol{\xi}_{k}\big)\psi_{\boldsymbol{\alpha}}\big(\boldsymbol{\xi}_{k}\big)v_{k}}\normalsize$\;
}
\NexusRegion{  
Compute mean $\small{\overline{\boldsymbol{T}}{}^{\ell}=\widehat{\boldsymbol{T}}{}_{\boldsymbol{0}}^{\ell}}\normalsize$ and covariance matrix $\small{\boldsymbol{C}{}_{\widehat{\boldsymbol{T}}}^{\ell}=\sum_{|\boldsymbol{\alpha}|=1}^{p}\widehat{\boldsymbol{T}}{}_{\boldsymbol{\alpha}}^{\ell}(\widehat{\boldsymbol{T}}{}_{\boldsymbol{\alpha}}^{\ell})^{\mathrm{T}}}\normalsize$\;
Solve eigenproblem $\small{\boldsymbol{W}^{\mathrm{T}}\boldsymbol{C}_{\widehat{\boldsymbol{T}}}^{\ell}\boldsymbol{W}\boldsymbol{\phi}^{j,\ell}=\lambda_{j}^{\ell}\boldsymbol{W}\boldsymbol{\phi}^{j,\ell}}\normalsize$\;
Choose $\small{d}\normalsize$ such that $\small{\sqrt{\sum_{j=d+1}^{r}\lambda_{j}^{\ell}}\leq\epsilon_{1}\sqrt{\sum_{|\boldsymbol{\alpha}|=1}^{p}(\widehat{\boldsymbol{T}}{}_{\boldsymbol{\alpha}}^{\ell})^{\mathrm{T}}\boldsymbol{W}\widehat{\boldsymbol{T}}{}_{\boldsymbol{\alpha}}^{\ell}}}\normalsize$\;
Compute PC coordinates of $\small{\eta_{j}^{\ell,p}}\normalsize$ by $\small{\eta_{j,\boldsymbol{\alpha}}^{\ell}=(\widehat{\boldsymbol{T}}{}_{\boldsymbol{\alpha}}^{\ell})^{\mathrm{T}}\boldsymbol{W}\boldsymbol{\phi}^{j,\ell}}\normalsize$ for $\small{j=1}\normalsize$ to $\small{d}\normalsize$\;
}
\FirstMonoProblem{
$\small{q=0}\normalsize$\;
\Repeat{$\small{\Big(\small{\sqrt{\sum_{|\boldsymbol{\beta}|=q}\|\widehat{\boldsymbol{\Phi}}{}^{\ell}_{\boldsymbol{\beta}}\|_{\boldsymbol{W}}^{2}}\leq\epsilon_{2}\sqrt{\sum_{|\boldsymbol{\beta}|=0}^{q}\|\widehat{\boldsymbol{\Phi}}{}^{\ell}_{\boldsymbol{\beta}}\|_{\boldsymbol{W}}^{2}}}\normalsize\;\Big)}\normalsize$}{
\NexusRegionn{
Compute PC basis $\small{\{\Gamma_{\boldsymbol{\beta}}^{\ell},\;0\leq|\boldsymbol{\beta}|\leq q\}}\normalsize$ up to total degree $\small{q}\normalsize$ w.r.t. $\small{P_{\boldsymbol{\eta}}^{\ell}}\normalsize$\;
Compute quadrature rule $\small{\{(\boldsymbol{\eta}_{k}^{\ell},w_{k}^{\ell}),\;1\leq k\leq \nu^{\ell}\}}\normalsize$ of level $\small{q+1}\normalsize$ w.r.t. $\small{P_{\boldsymbol{\eta}}^{\ell}}\normalsize$\;
}  
\For{$k=1$ \KwTo $\,\nu^{\ell}$}{
Solve $\small{\big[\boldsymbol{D}\big(\widehat{\boldsymbol{T}}{}^{\ell,p,d}(\boldsymbol{\eta}^{\ell}_{k})\big)+\boldsymbol{M}\big(\widehat{\boldsymbol{T}}{}^{\ell,p,d}(\boldsymbol{\eta}^{\ell}_{k})\big)\big]\widehat{\boldsymbol{\Phi}}{}^{\ell}(\boldsymbol{\eta}^{\ell}_{k})=\boldsymbol{s},}\normalsize$\\
with $\small{\widehat{\boldsymbol{T}}{}^{\ell,p,d}(\boldsymbol{\eta}^{\ell}_{k})=\overline{\boldsymbol{T}}{}^{\ell}+\sum_{j=1}^{d}\sqrt{\lambda_{j}^{\ell}}\eta_{j,k}^{\ell}\boldsymbol{\phi}^{j,\ell}}\normalsize$\;}
Compute PC coordinates of $\small{\widehat{\boldsymbol{\Phi}}{}^{\ell,q}}\normalsize$ using $\small{\widehat{\boldsymbol{\Phi}}{}^{\ell}_{\boldsymbol{\beta}}=\sum_{k=1}^{\nu^{\ell}}\widehat{\boldsymbol{\Phi}}{}^{\ell}(\boldsymbol{\eta}^{\ell}_{k})\Gamma_{\boldsymbol{\beta}}^{\ell}(\boldsymbol{\eta}^{\ell}_{k})w^{\ell}_{k}}\normalsize$\;
$\small{q=q+1}\normalsize$\;
\vspace{1mm}
}
}
$\small{\ell=\ell+1}\normalsize$\;}
\caption{Implementation of the illustration problem.}\label{algo:algo6}
\end{algorithm}

\section{Numerical results}\label{sec:sec7}
We obtained numerical results using the following properties.  
We assumed the reactor to have a length of~$L=100\,[\text{cm}]$.
Further, we assumed a deterministic and position-independent heat conductivity~$k=100\,\text{$[\text{J/K/cm/s}]$}$; ambient temperature~$T_{\infty}=390\,[\text{K}]$; fission energy~$E_{\text{f}}=3.0E\text{-}11\,[\text{J/neutrons}]$; fission cross section~$\Sigma_{\text{a,ref}}=0.0075\,[\text{cm}^{-1}]$; neutron-diffusion constant~$D_{\text{ref}}=2.2\,[\text{cm}]$; absorption cross section $\Sigma_{\text{a,ref}}=0.0195\,[\text{cm}^{-1}]$; multiplication factor~$\nu=2.2$; neutron source~$s=5.0E11\,[\text{neutrons/s/cm}^{3}$]; and temperatures~$T_{\text{ref}}=390\,[\text{K}]$, $T_{\text{min}}=390\,[\text{K}]$, and~$T_{\text{max}}=1000\,[\text{K}]$.

\begin{figure}[htp]
  \begin{center}
    \subfigure[Samples.]{\includegraphics[width=0.8\textwidth]{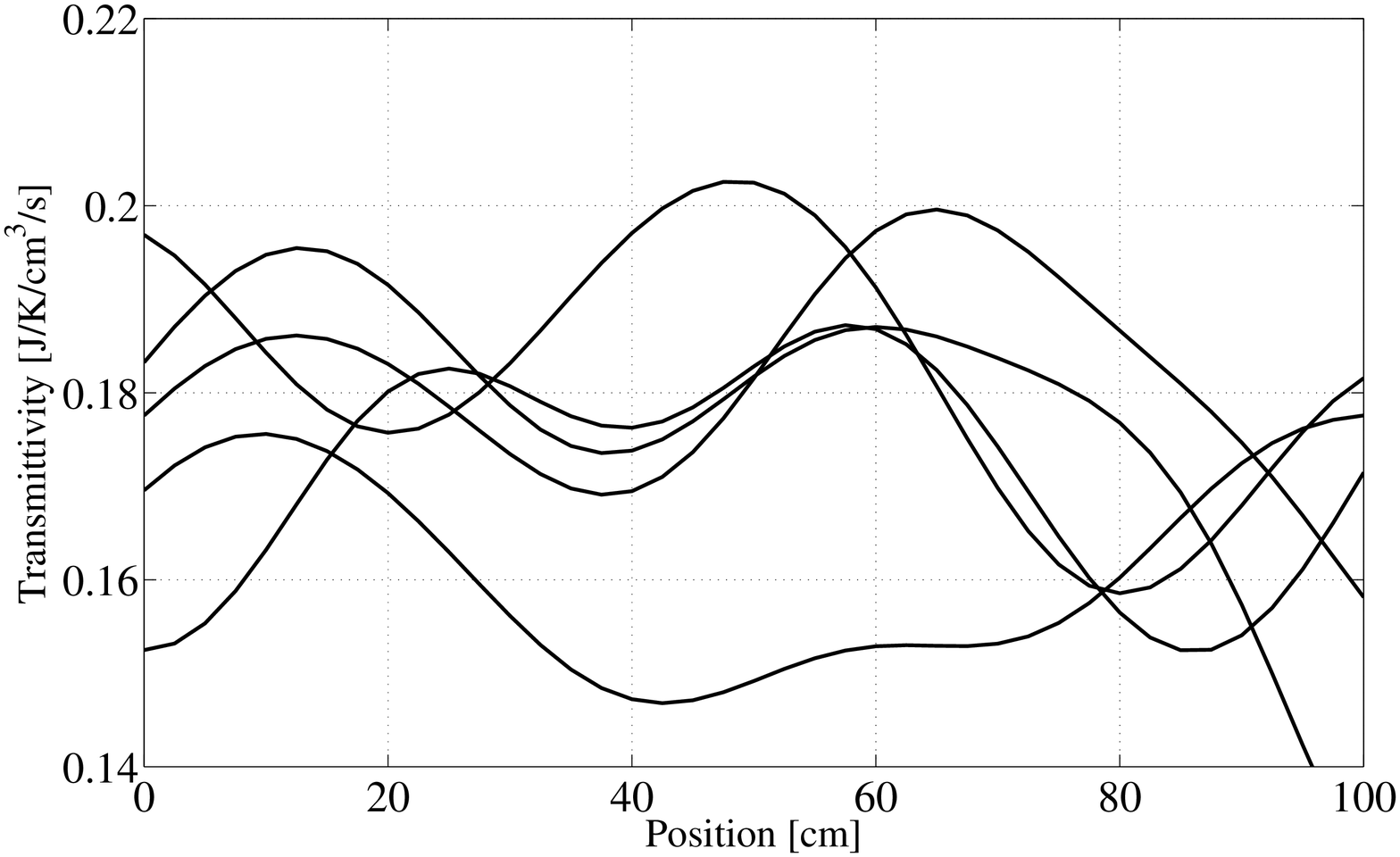}}
    \hfill
    \subfigure[Eigenvalues.]{\includegraphics[width=0.8\textwidth]{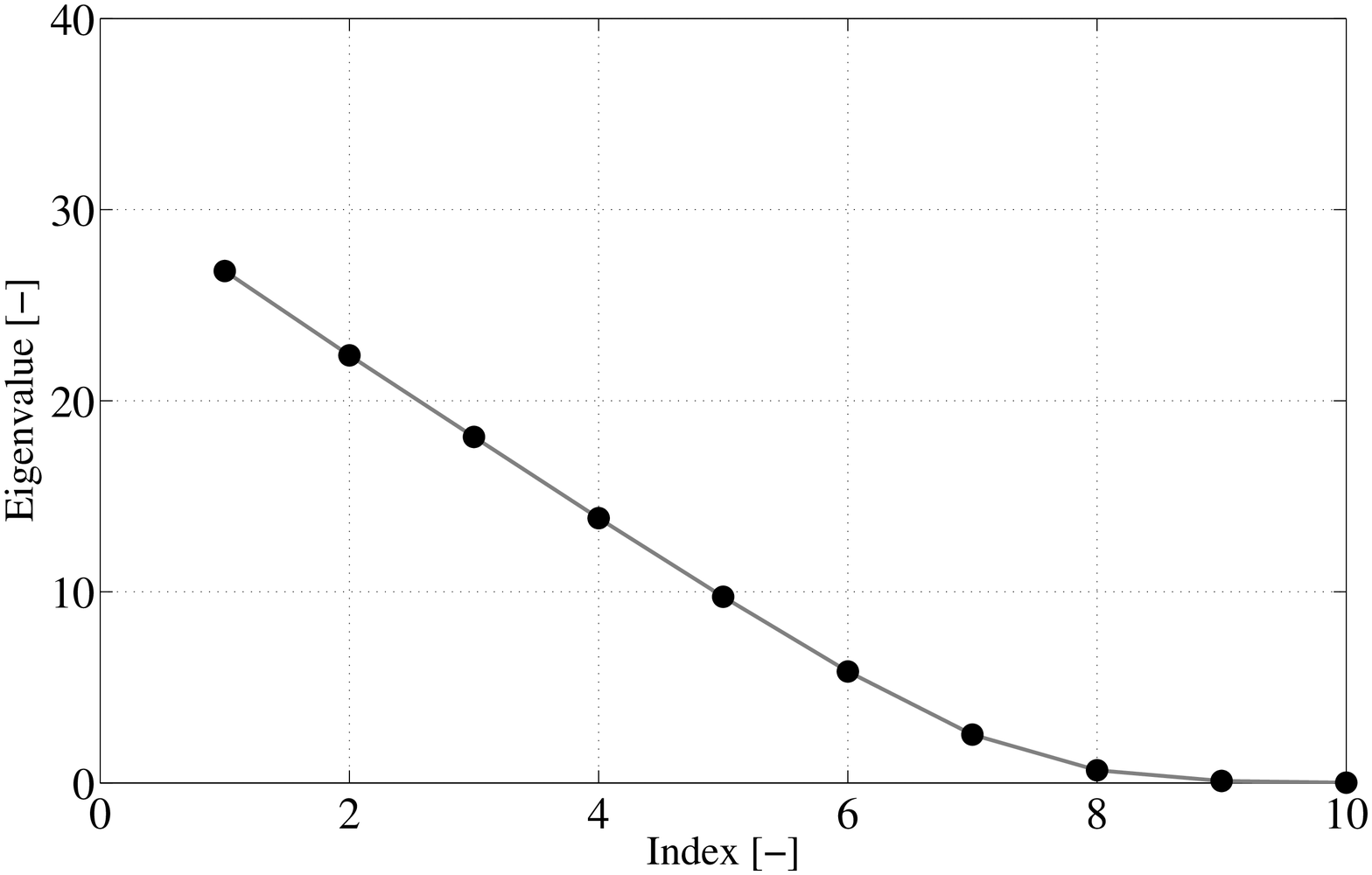}}
    \caption{Thermal transmittivity random field: (a)~five samples and~(b)~ten largest magnitude eigenvalues of the covariance integral operator.}\label{fig:figure1}
  \end{center}
\end{figure}

In addition, we used a thermal transmittivity random field with position-independent mean~$\overline{h}=0.17\,[\text{J/K/cm}^{3}\text{/s}]$, spatial correlation length~$a=15\,[\text{cm}]$, and coefficient of variation~$\delta=10\,\%$.
We retained~$m=10$ terms in expansion~(\ref{eq:hN}).
Figure~\ref{fig:figure1}(a) shows a few sample paths of the random field~$\{h(x,\cdot),\;0\leq x\leq L\}$ thus obtained.
Figure~\ref{fig:figure1}(b) shows the 10 largest magnitude eigenvalues of the covariance integral operator.  

\subsection{Monte Carlo sampling implementation}
\begin{figure}[htp]
  \begin{center}
    \subfigure[Temperature.]{\includegraphics[width=0.8\textwidth]{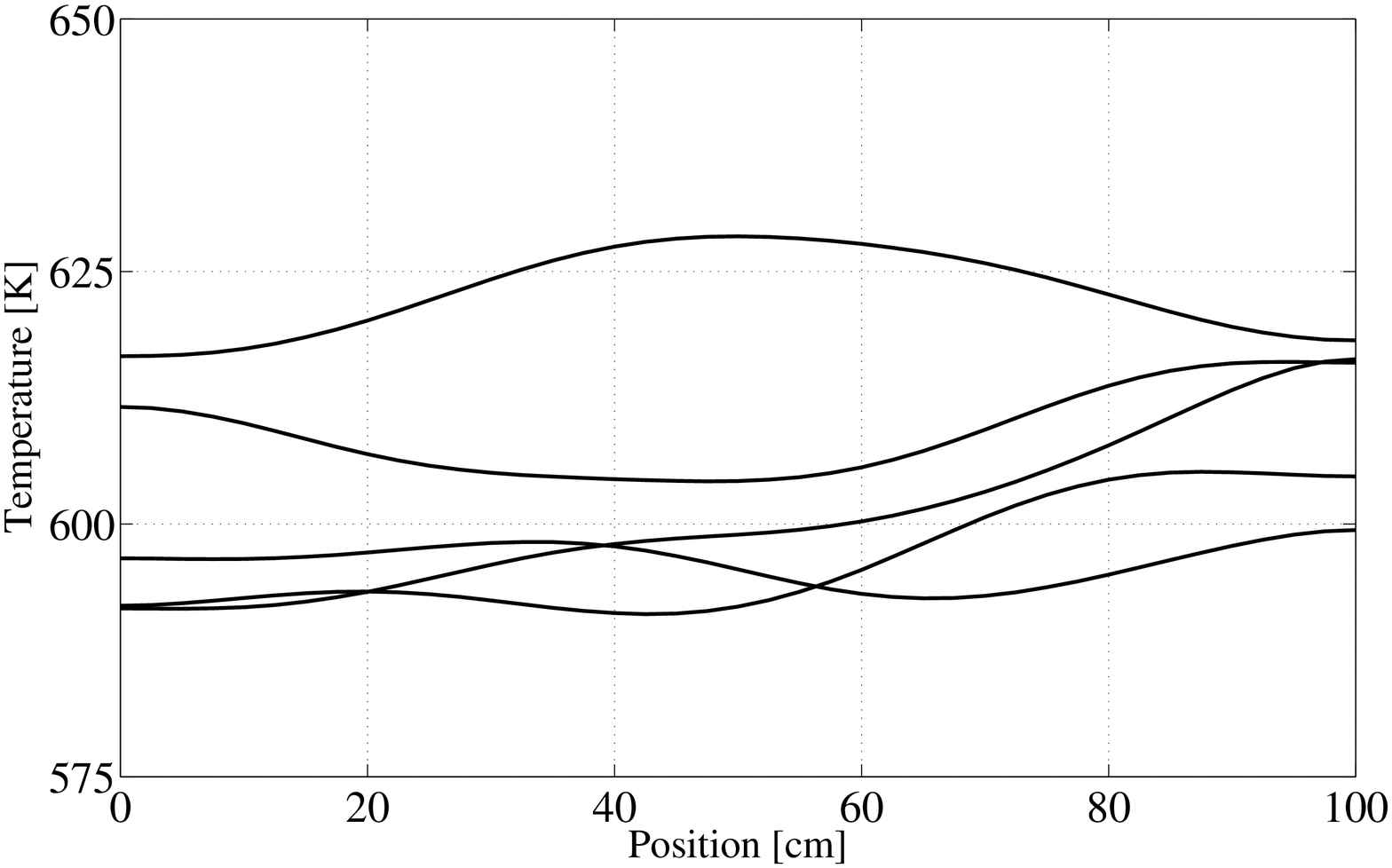}}
    \hfill
    \subfigure[Neutron flux.]{\includegraphics[width=0.8\textwidth]{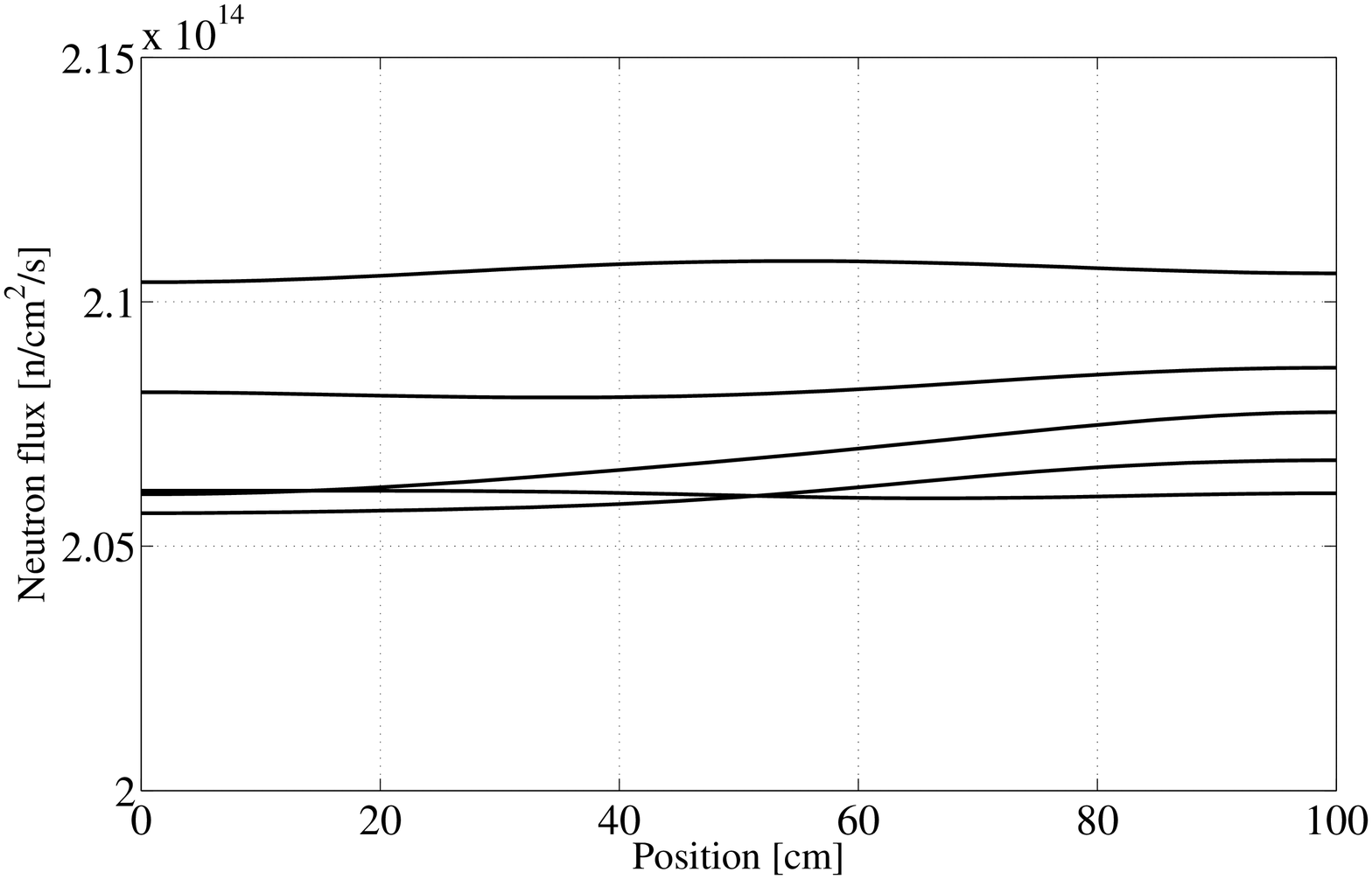}}
    \caption{Monte Carlo simulation: five samples of the solution.}\label{fig:figure2}
  \end{center}
\end{figure}
First, we carried out a Monte Carlo simulation.
We generated $MC=100,000$ sample paths of the thermal transmittivity random field.
Then, for each of these sample paths, we constructed the associated deterministic multiphysics model, each of which we solved using the FE method for the spatial discretization and Gauss-Seidel iteration as the iterative method.
We systematically obtained converged results for~$r-1=40$ finite elements and~$20$ iterations.  

Figure~\ref{fig:figure2} shows a few samples of the random temperature and neutron flux thus obtained.
We can observe that the samples of the random temperature (Fig.~\ref{fig:figure2}(b)) are smoother than the samples of the thermal transmittivity random field~(Fig.~\ref{fig:figure1}(a)); i.e., the former exhibit less rapid oscillations with respect to the position in the reactor than the latter.
We had demonstrated in~\citep{arnst2011a} that this behavior can be attributed to the large magnitude of the diffusion term of the heat subproblem which reduces the nonuniformity of the samples of the random temperature.

\subsection{PC-based implementation involving dimension reduction and measure transformation}
Next, we implemented the proposed PC-based iterative method involving dimension reduction and measure transformation.
This implementation exactly corresponded to Algorithm~\ref{algo:algo6}.
We obtained the results to follow by systematically setting the total degree of the PC expansion of the random temperature to $p=4$ and, with reference to~(\ref{eq:criterion1}) and (\ref{eq:criterion2}), using a range of values for the error tolerance levels $\epsilon_{1}$ and~$\epsilon_{2}$ adopted to determine the reduced dimension and the total degree of the PC expansion of the random neutron flux at each iteration.
We discuss convergence as a function of these error tolerance levels later.
Now, we present detailed results obtained for~$\epsilon_{1}=0.05$ and~$\epsilon_{2}=0.000001$.

\begin{figure}[htp]
  \begin{center}
    \includegraphics[width=0.8\textwidth]{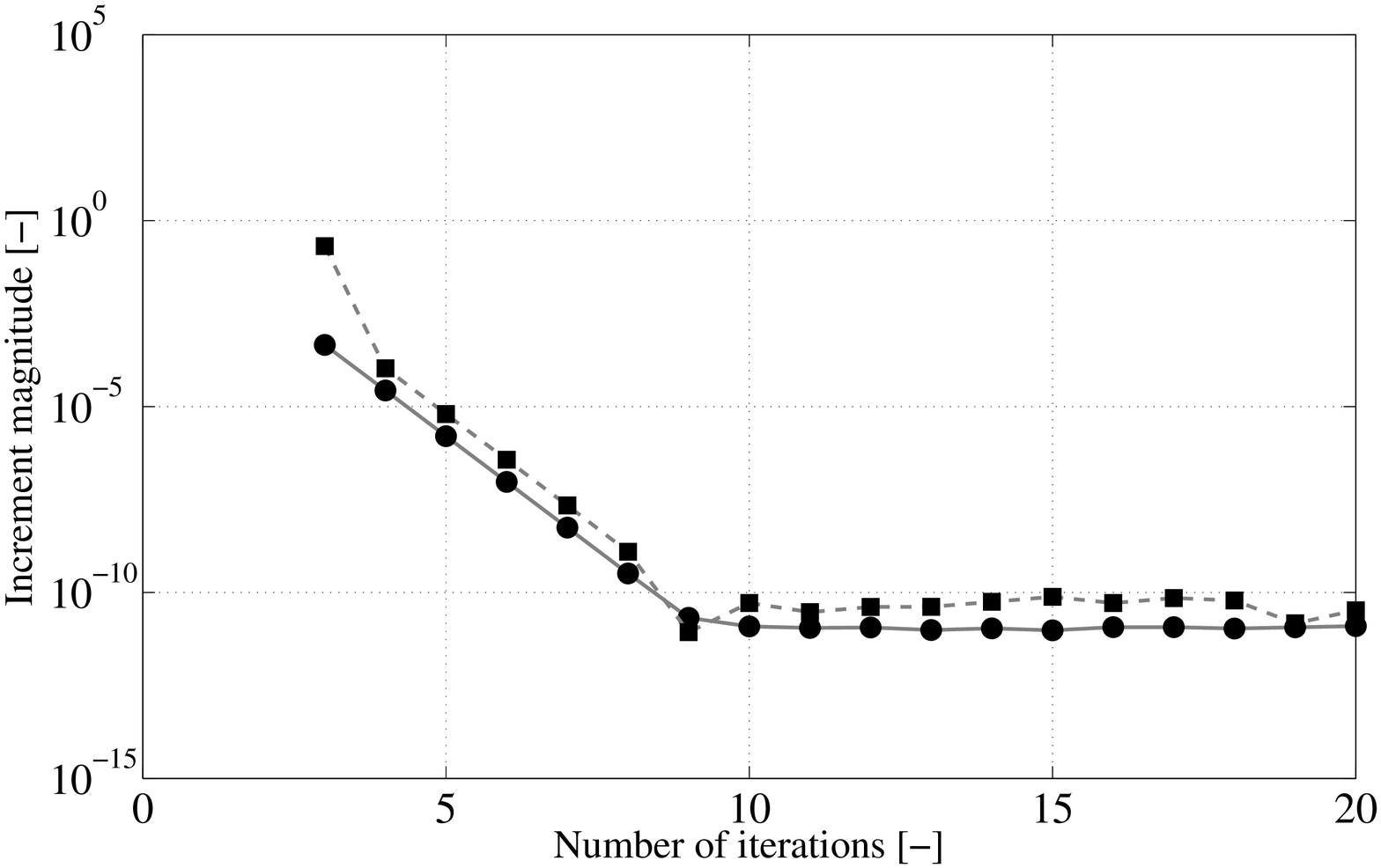}
    \centerline{$\small{\ell\mapsto\sqrt{\frac{1}{MC}\sum_{k=1}^{MC}\|\widehat{\boldsymbol{T}}{}^{\ell,p}(\boldsymbol{\xi}_{k})-\widehat{\boldsymbol{T}}{}^{\ell-1,p}(\boldsymbol{\xi}_{k})\|_{\boldsymbol{W}}^{2}}\Big/\sqrt{\frac{1}{MC}\sum_{k=1}^{MC}\|{\boldsymbol{T}}{}^{\infty}(\boldsymbol{\xi}_{k})\|_{\boldsymbol{W}}^{2}}}\normalsize$ (squares).} 
    \centerline{$\small{\ell\mapsto\sqrt{\frac{1}{MC}\sum_{k=1}^{MC}\|\widehat{\boldsymbol{\Phi}}{}^{\ell,q}(\boldsymbol{\eta}^{\ell,p}(\boldsymbol{\xi}_{k}))-\widehat{\boldsymbol{\Phi}}{}^{\ell-1,q}(\boldsymbol{\eta}^{\ell-1,p}(\boldsymbol{\xi}_{k}))\|_{\boldsymbol{W}}^{2}}\Big/\sqrt{\frac{1}{MC}\sum_{k=1}^{MC}\|{\boldsymbol{\Phi}}{}^{\infty}(\boldsymbol{\xi}_{k})\|_{\boldsymbol{W}}^{2}}}\normalsize$ (circles).}
    \caption{PC-based simulation: convergence with respect to the number of iterations.}\label{fig:figure3}
  \end{center}
\end{figure}

Figure~\ref{fig:figure3} shows the convergence of the iterative method as a function of the number of iterations; note that the superscript $\infty$ is used in the figure captions to indicate convergence with respect to the number of iterations.
The iterative method converged at a linear rate up to approximately iteration~$\ell=10$, after which linear-solver tolerances became dominant and prevented further convergence. 
All results to follow were obtained at iteration~$\ell=20$ and can thus be considered to have converged with respect to the number of iterations.

\begin{figure}[htp]
  \begin{center}
    \subfigure[Mean.]{\includegraphics[width=0.49\textwidth]{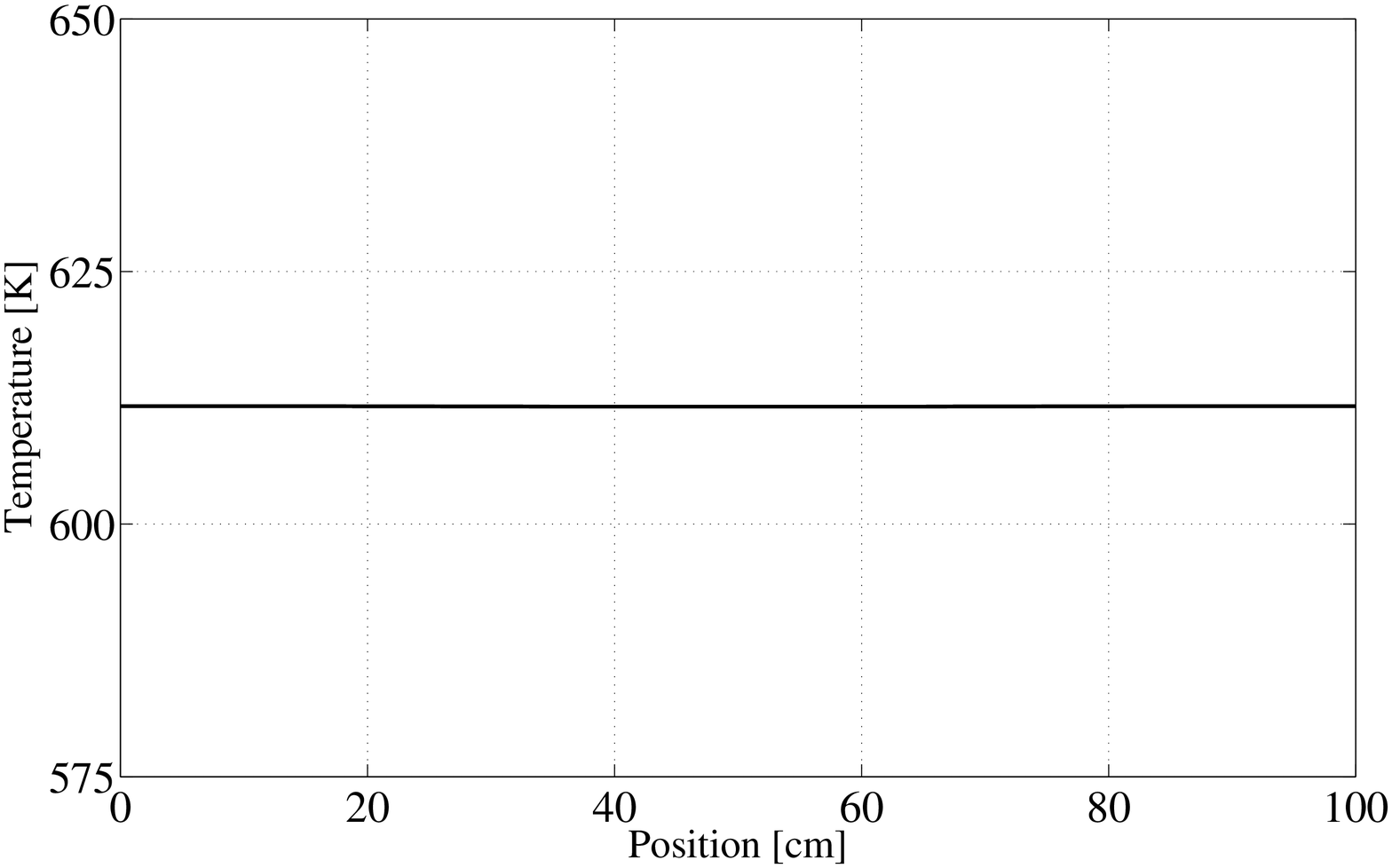}}
    \hfill
    \subfigure[Eigenvalues.]{\includegraphics[width=0.49\textwidth]{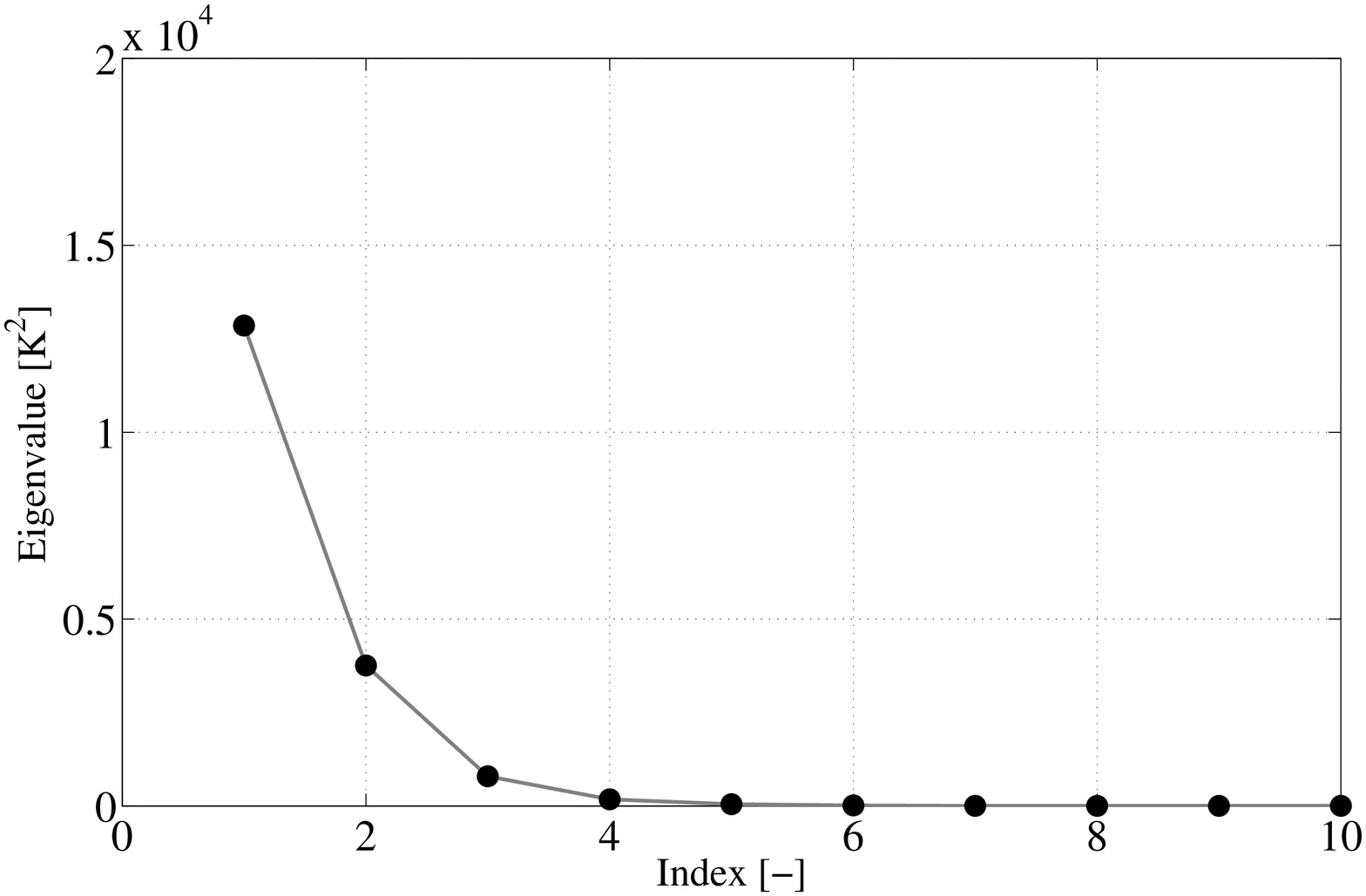}}
    \subfigure[First eigenmode.]{\includegraphics[width=0.49\textwidth]{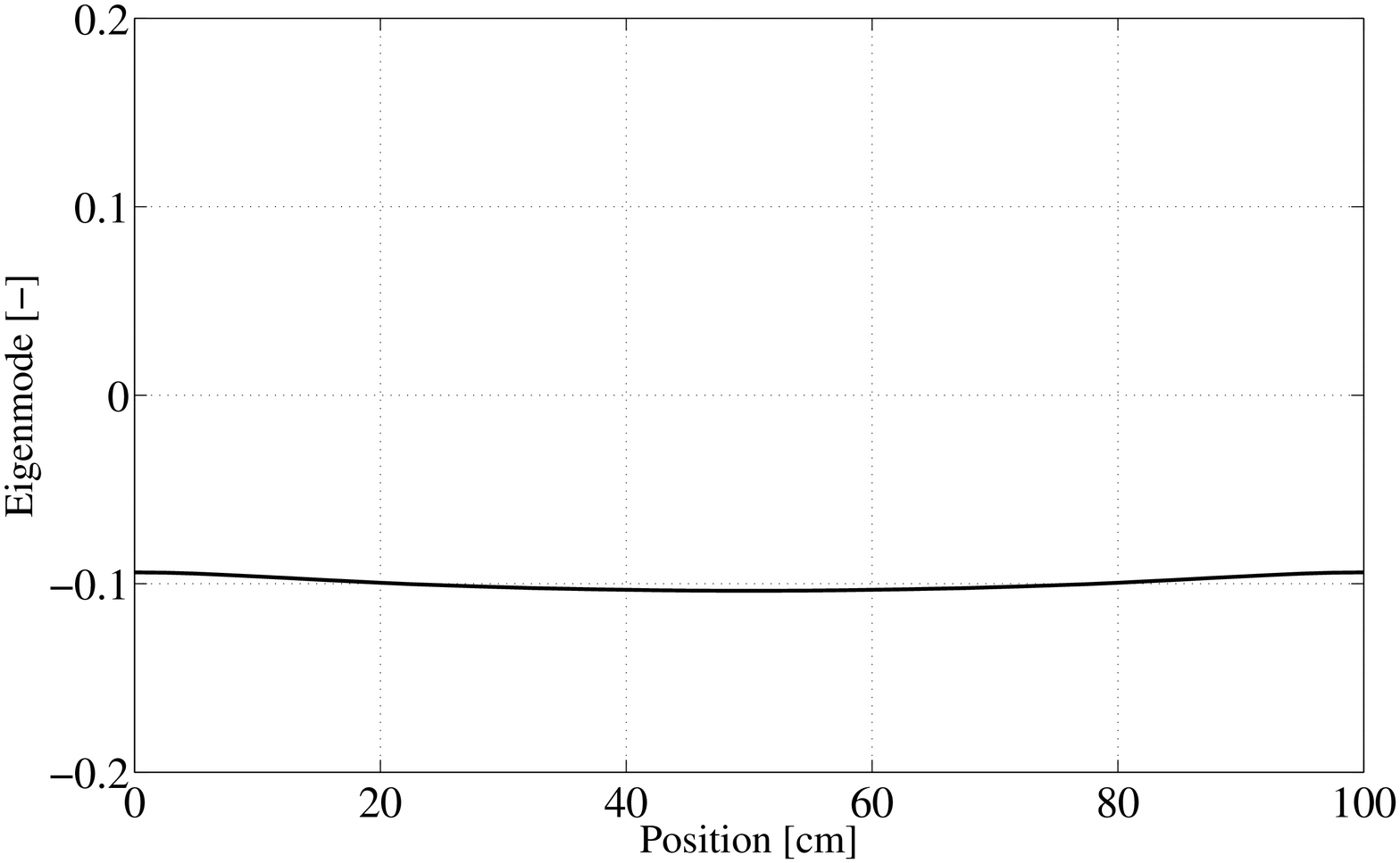}}
    \hfill
    \subfigure[Second eigenmode.]{\includegraphics[width=0.49\textwidth]{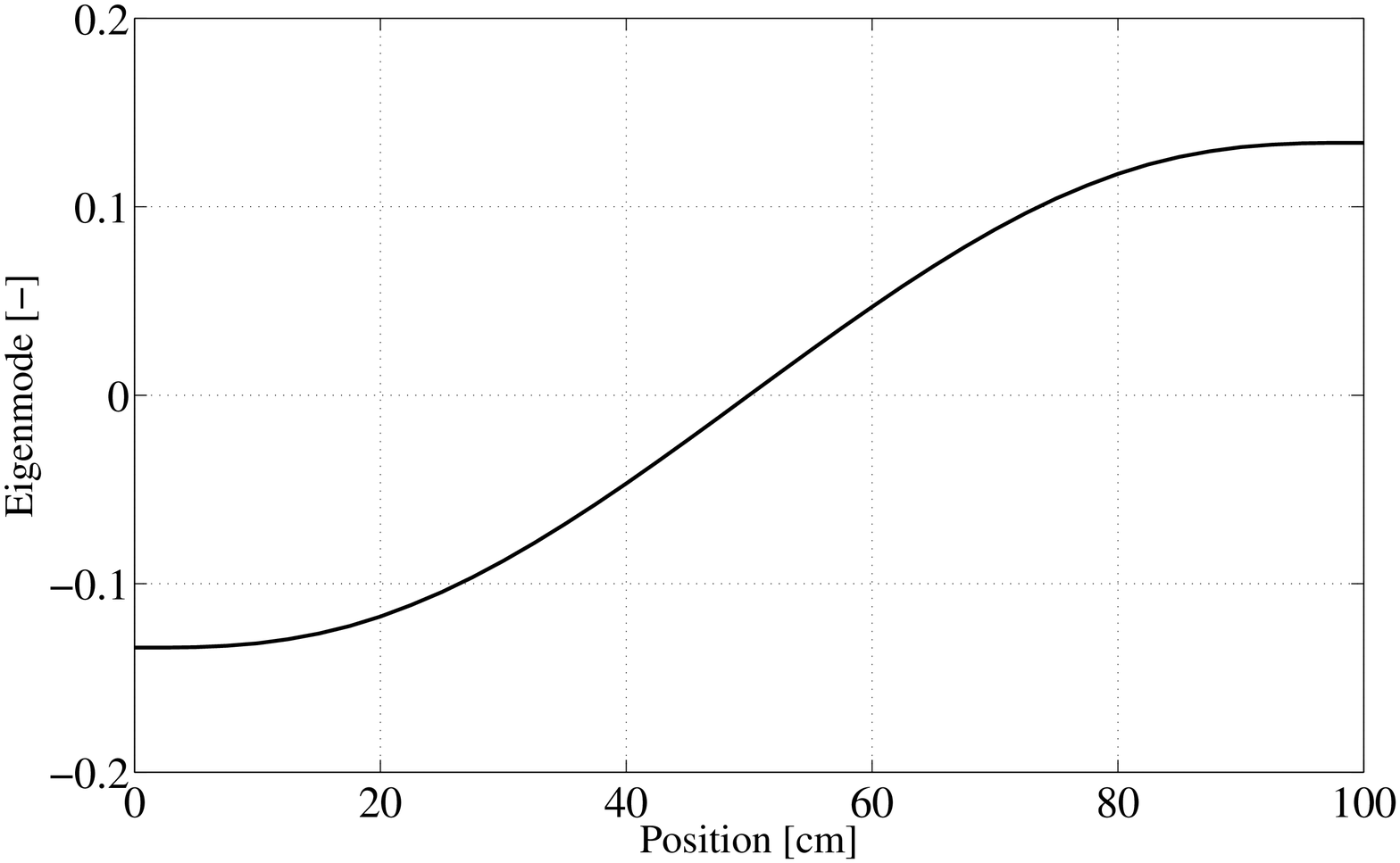}}
    \subfigure[First reduced random variable.]{\includegraphics[width=0.49\textwidth]{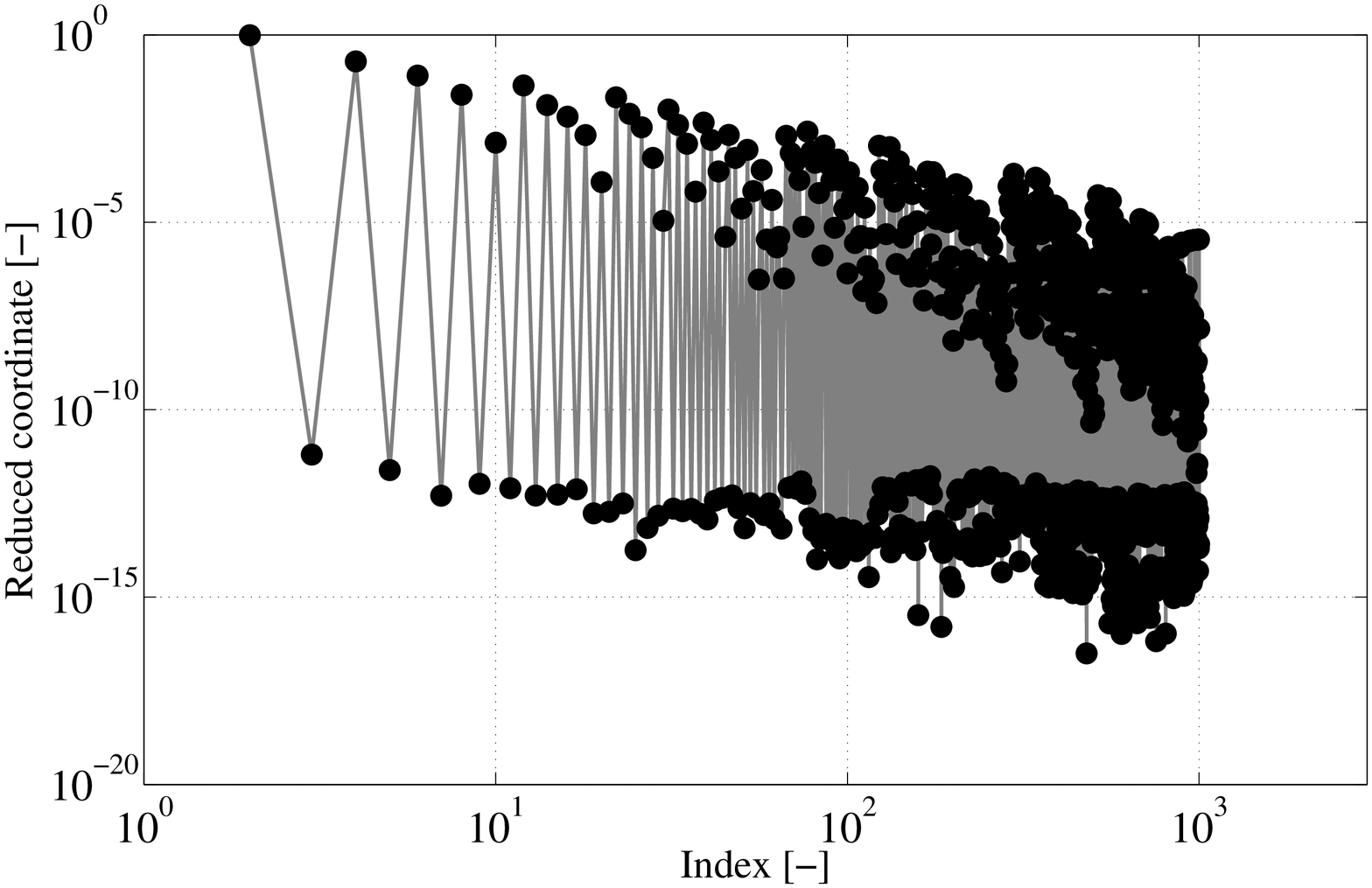}}
    \hfill
    \subfigure[Second reduced random variable.]{\includegraphics[width=0.49\textwidth]{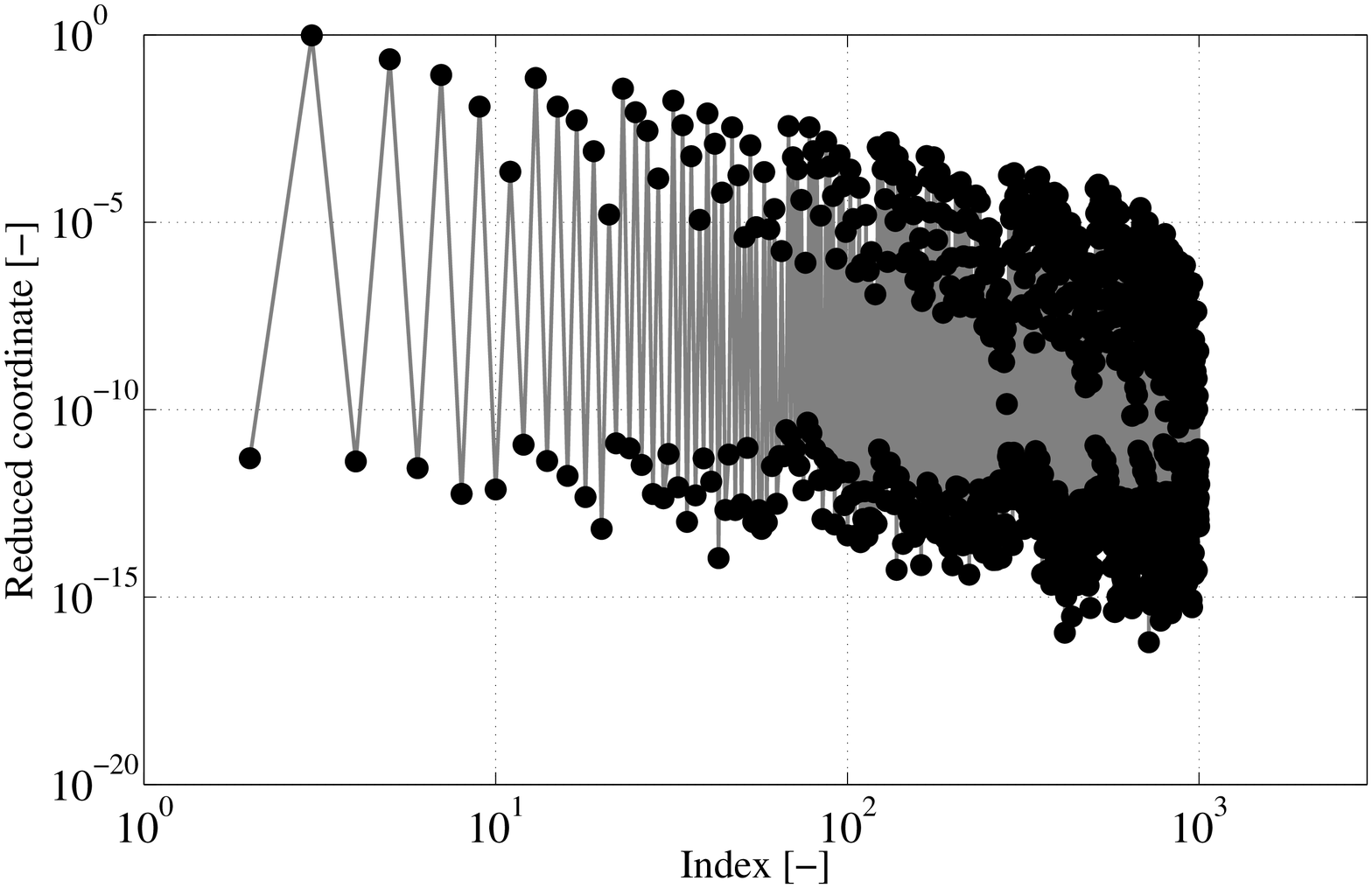}}
    \caption{PC-based simulation: mean, eigenvalues, first and second eigenmode, and PC coordinates of the first and second reduced random variables of the KL~decomposition of the random temperature.}\label{fig:figure4}
  \end{center}
\end{figure}

\begin{figure}[htp]
  \begin{center}
    \subfigure[Joint probability density function.]{\includegraphics[width=0.55\textwidth]{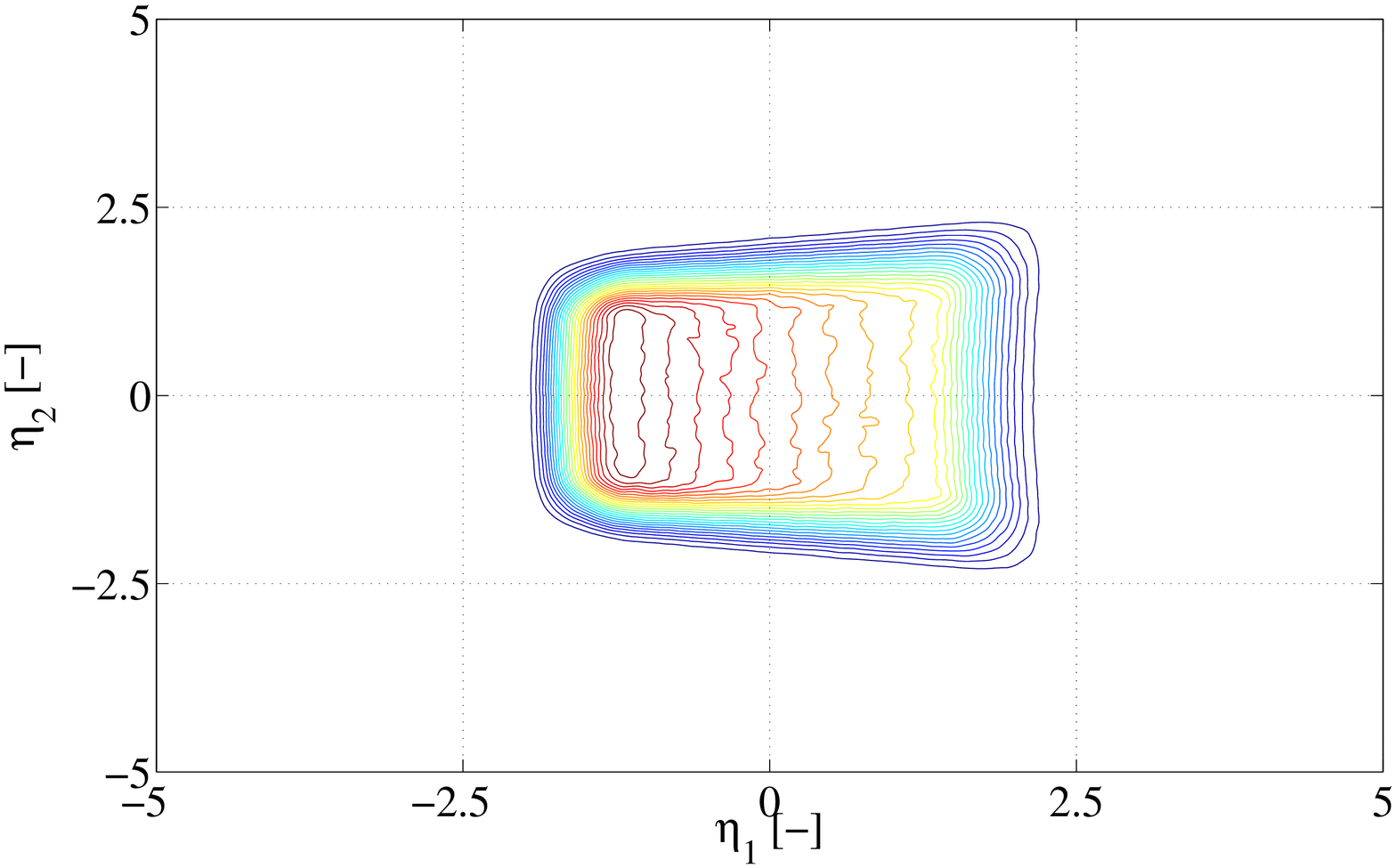}}
    \vfill
    \subfigure[First marginal probability density function.]{\includegraphics[width=0.55\textwidth]{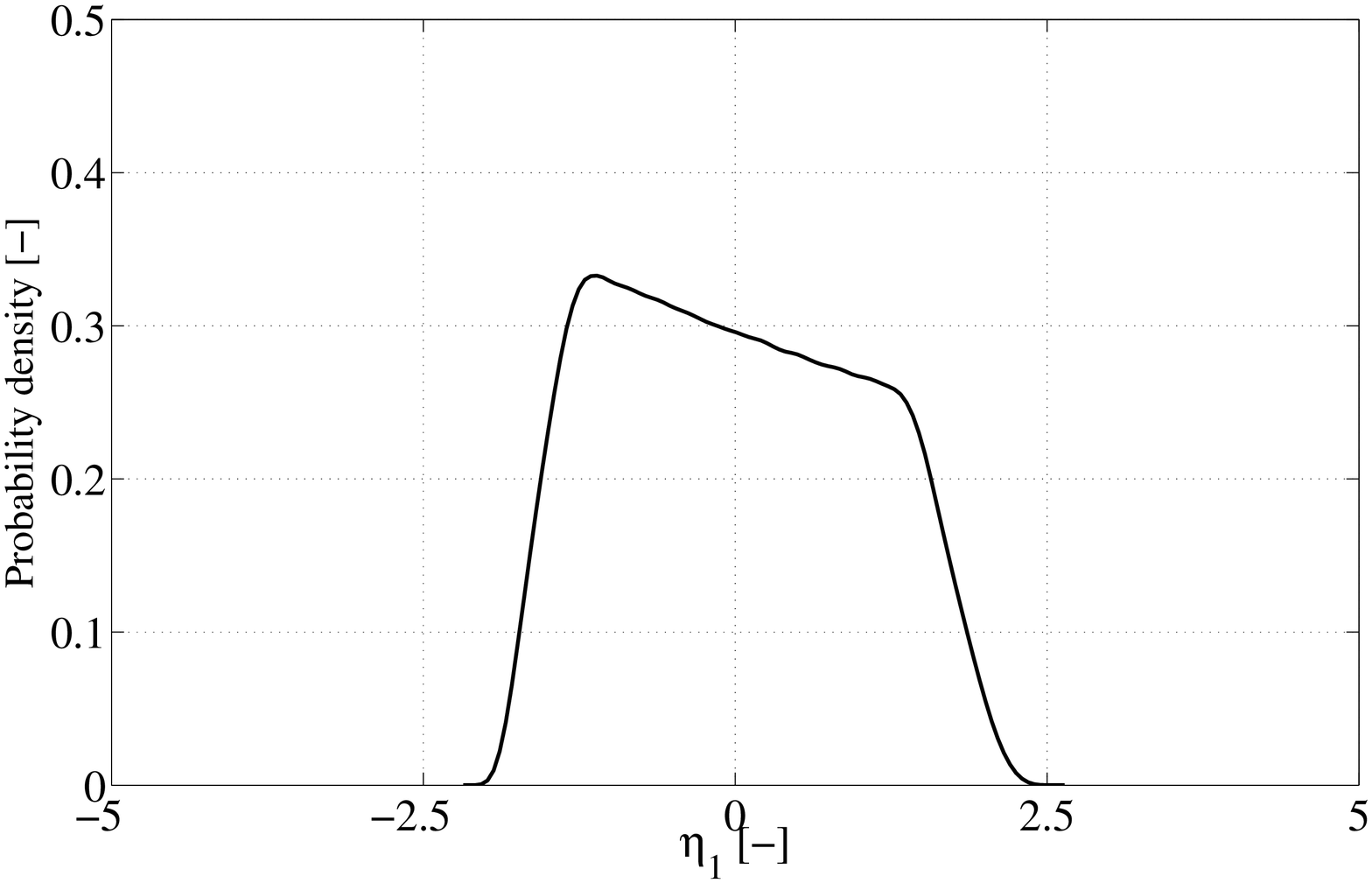}}
     \vfill
    \subfigure[Second marginal probability density function.]{\includegraphics[width=0.55\textwidth]{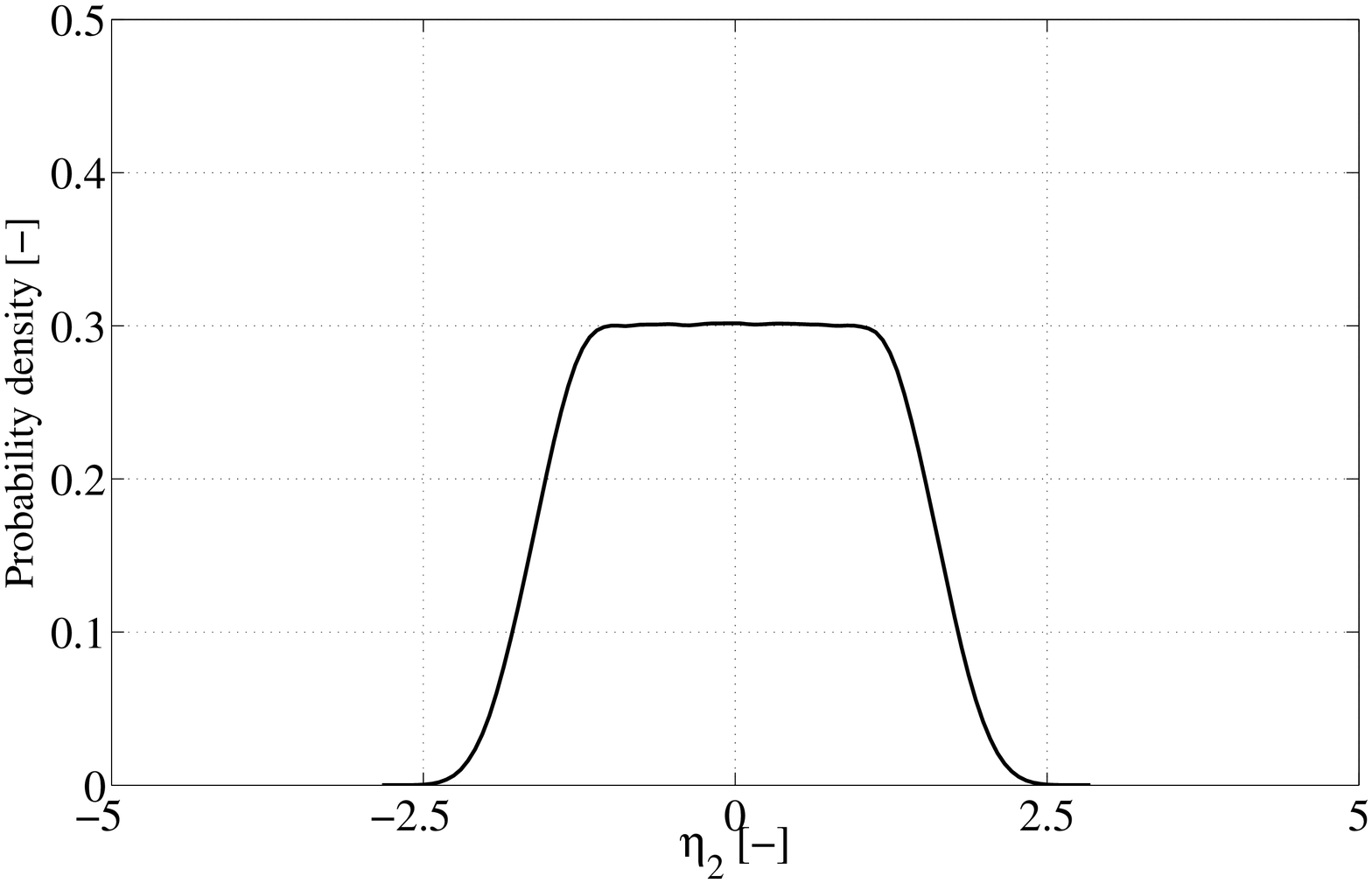}}
    \caption{PC-based simulation: joint and marginal probability density functions of the first and second reduced random variables of the KL~decomposition of the random temperature.}\label{fig:figure5}
  \end{center}
\end{figure}

Figure~\ref{fig:figure4} shows a few components of the KL decomposition of the random temperature, namely, the mean temperature, the 10 largest magnitude eigenvalues, and the eigenmodes and reduced random variables associated with the two largest magnitude eigenvalues.
We can observe that the eigenvalues of the KL decomposition of the random temperature (Fig.~\ref{fig:figure4}(b)) decay at a faster rate than those of the KL decomposition of the thermal transmittivity random field (Fig.~\ref{fig:figure1}(b)), which is consistent with our earlier observation that the samples of the former are smoother than those of the latter.

Figure~\ref{fig:figure5} shows the joint and marginal probability density functions of the reduced random variables.
Clearly, the joint probability density function exhibits statistical dependence and the marginal probability density functions are not ``labeled."

\begin{figure}[htp]
  \begin{center}
     \includegraphics[width=.8\textwidth]{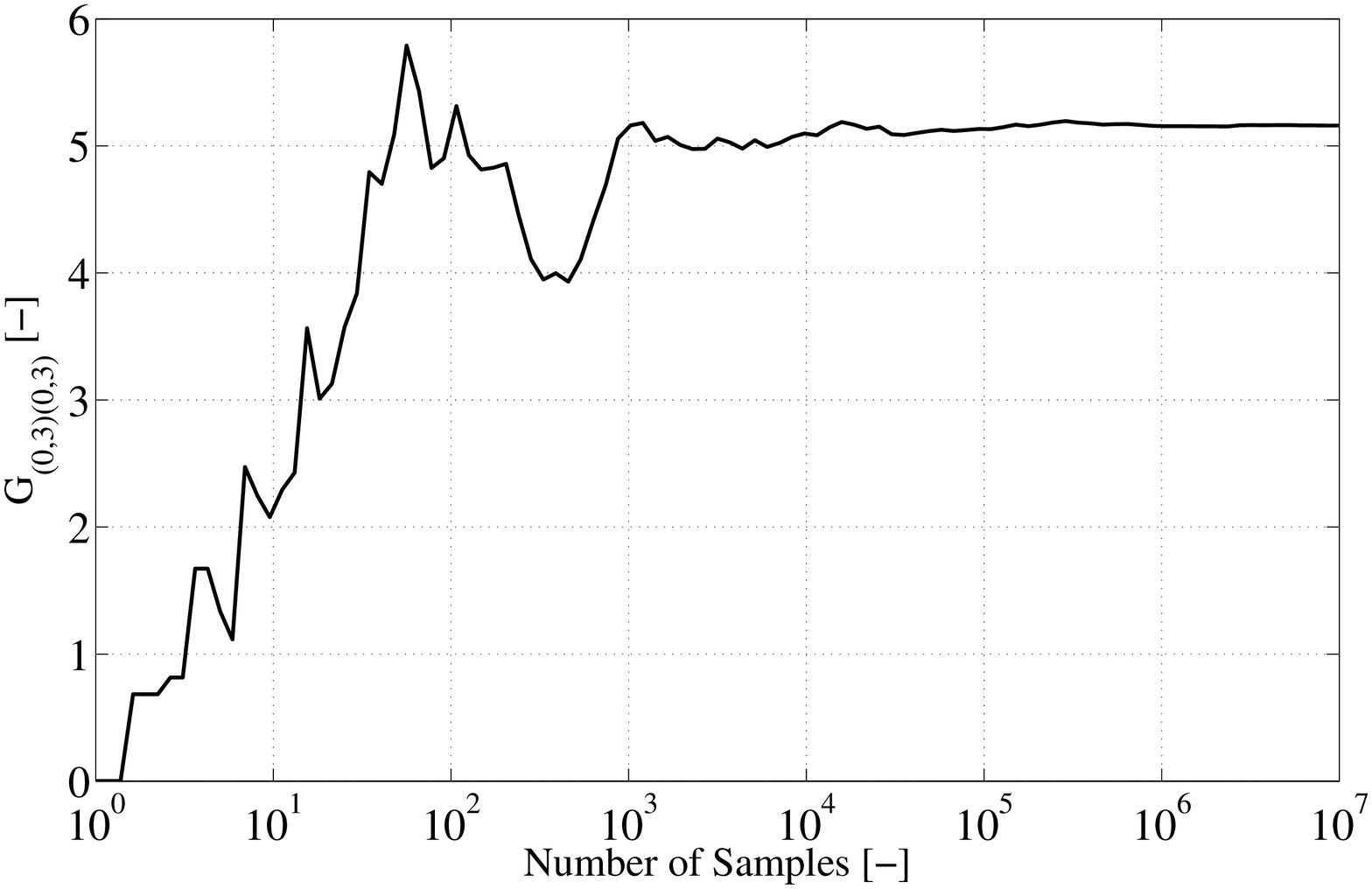}
    \caption{PC-based simulation: convergence of the Monte Carlo estimate of an entry of the Gram matrix with respect to the number of samples.}\label{fig:figure6}
  \end{center}
\end{figure}

\begin{figure}[htp]
  \begin{center}
    \subfigure[$\Gamma_{00}$.]{\includegraphics[width=0.22\textwidth]{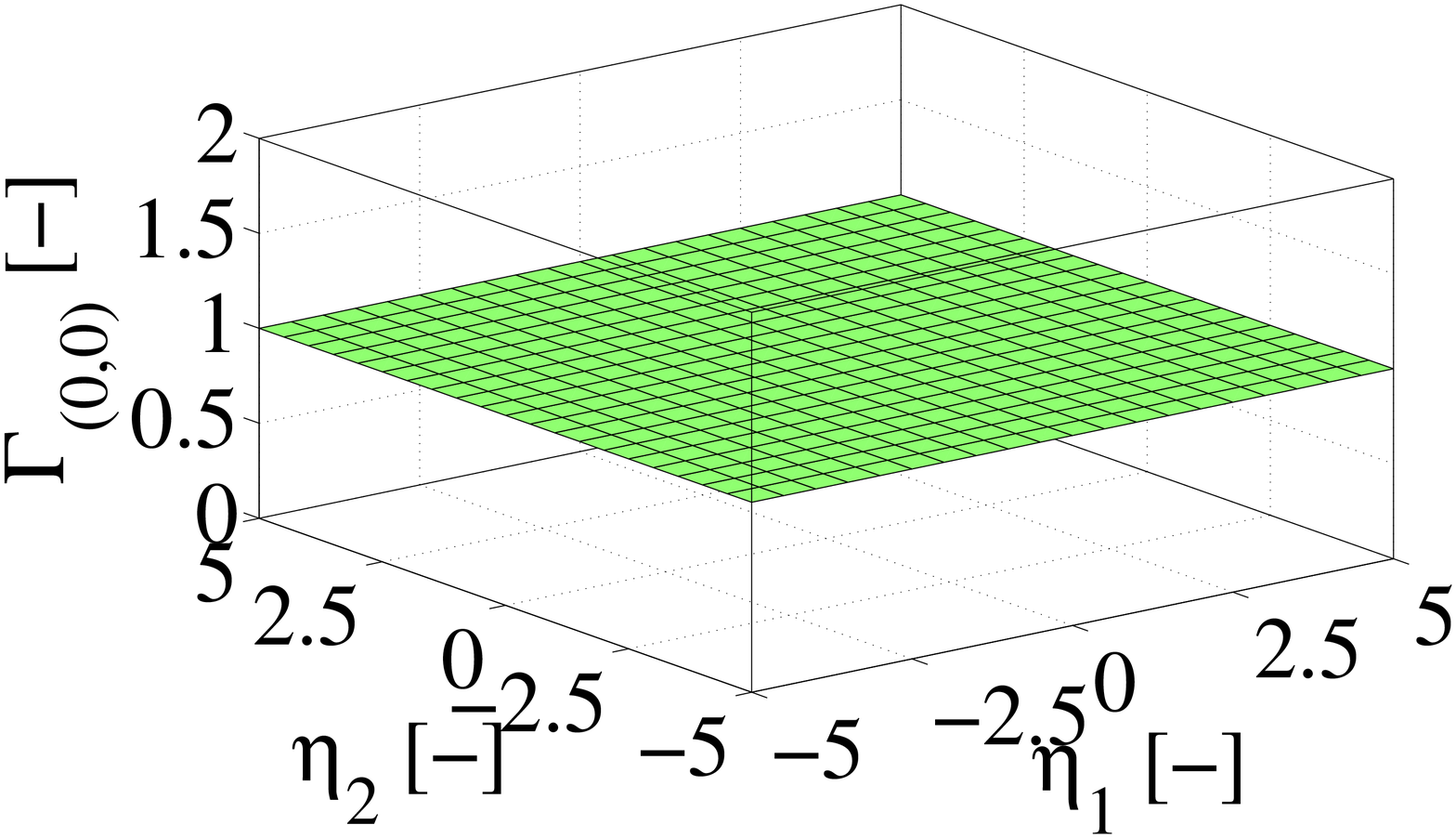}}
    \vfill
    \subfigure[$\Gamma_{10}$.]{\includegraphics[width=0.22\textwidth]{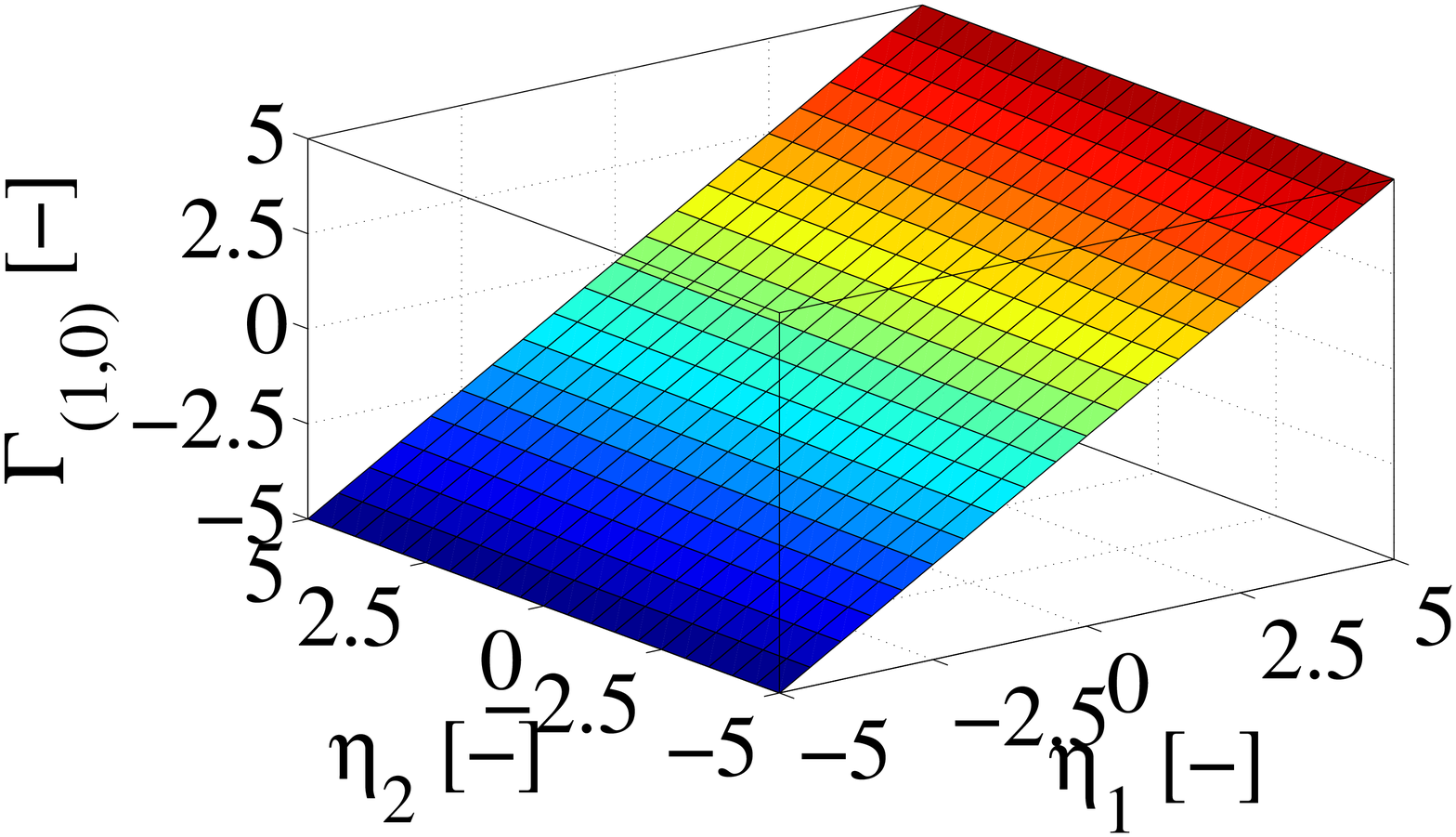}}
    \hspace{1mm}
    \subfigure[$\Gamma_{01}$.]{\includegraphics[width=0.22\textwidth]{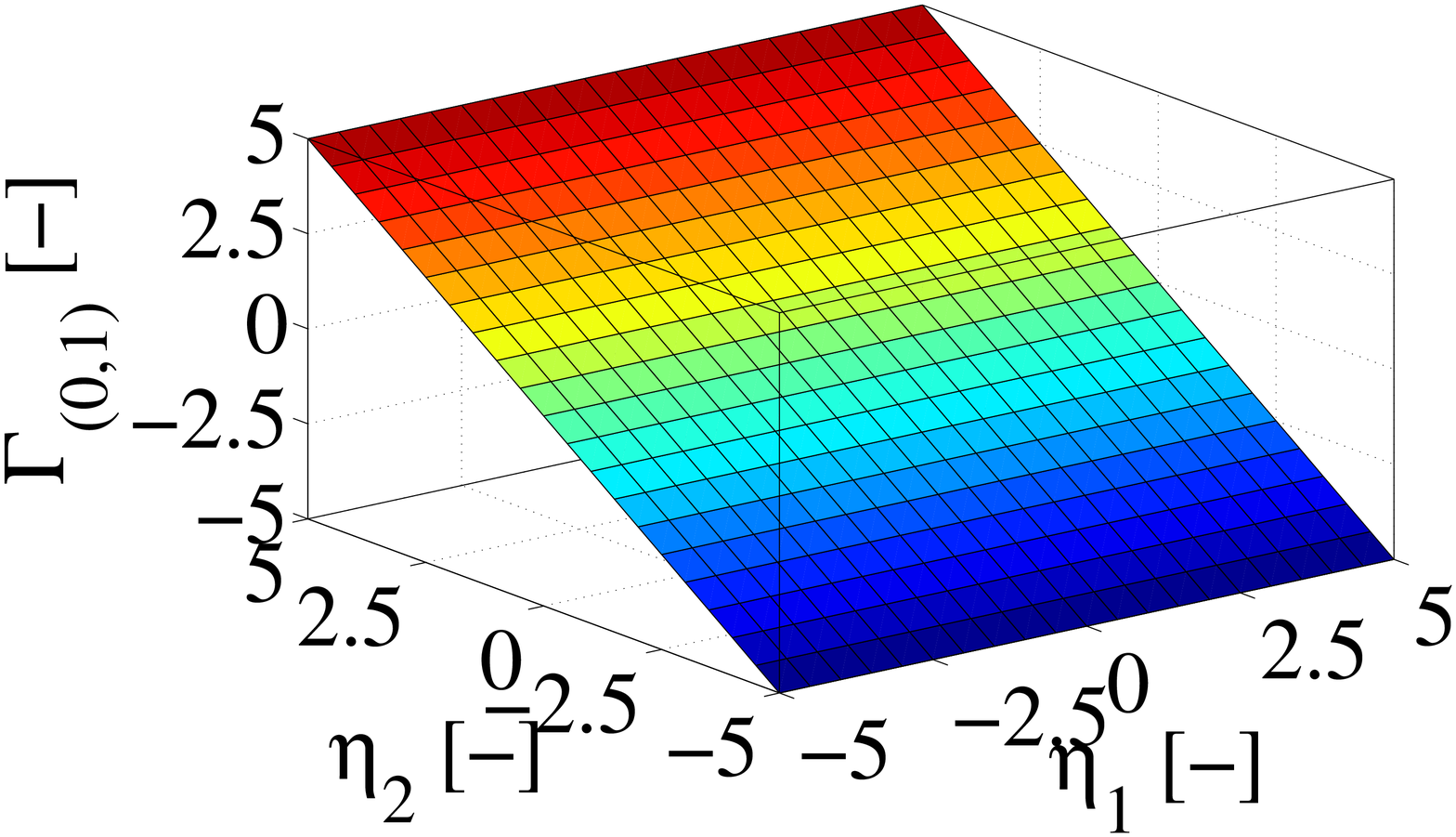}}
    \vfill
    \subfigure[$\Gamma_{20}$.]{\includegraphics[width=0.22\textwidth]{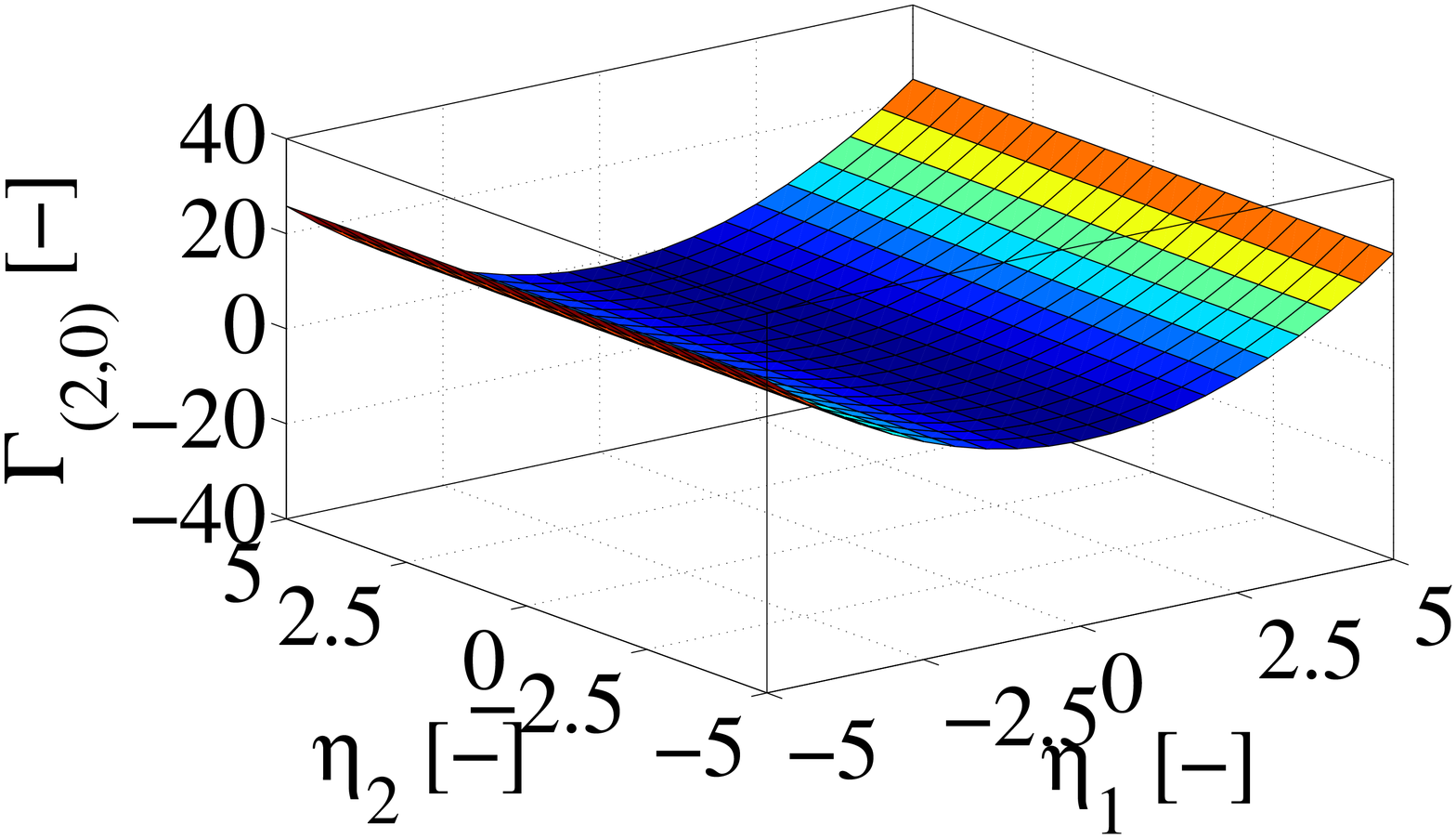}}
    \hspace{1mm}
    \subfigure[$\Gamma_{11}$.]{\includegraphics[width=0.22\textwidth]{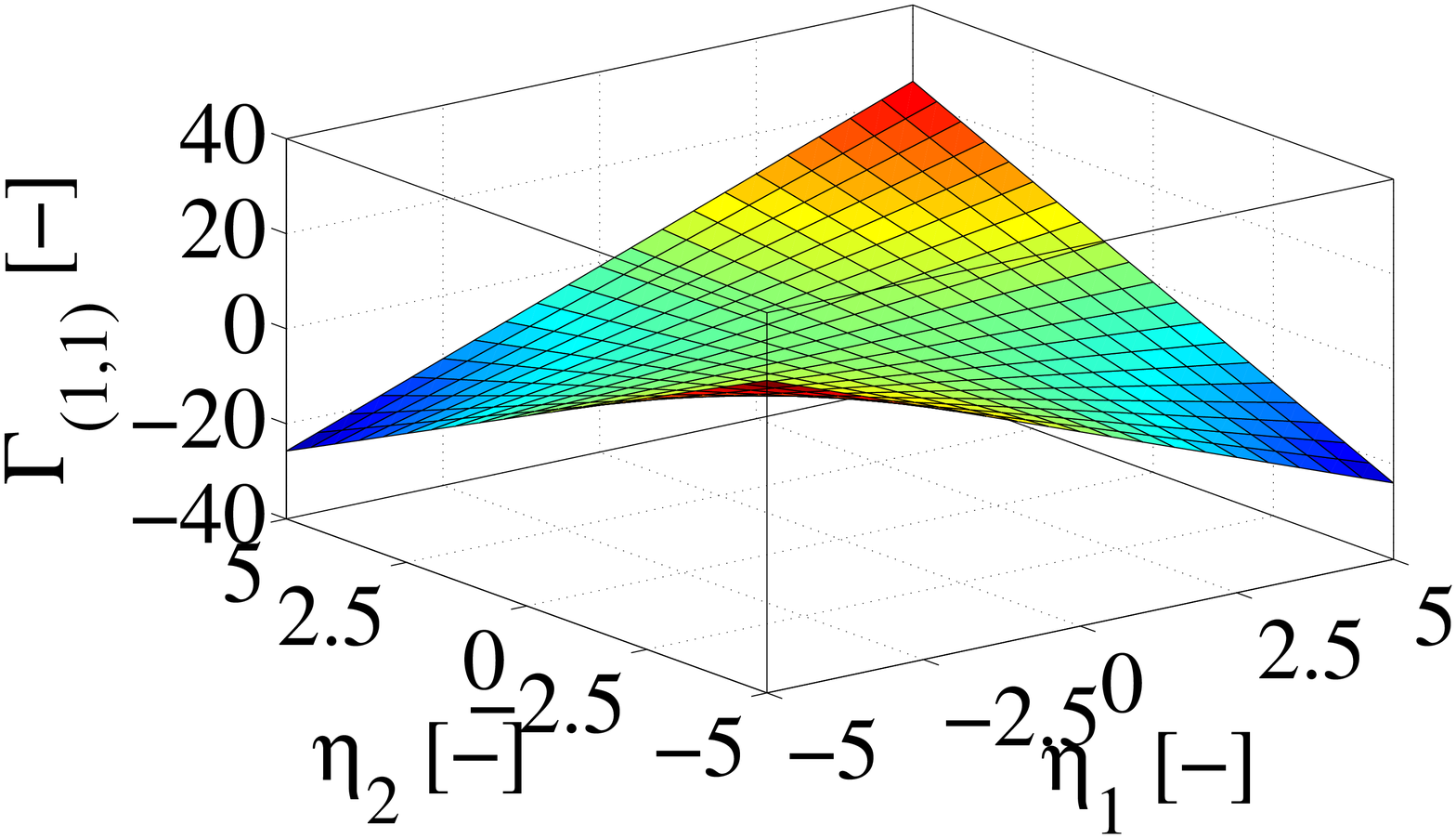}}
    \hspace{1mm}
    \subfigure[$\Gamma_{02}$.]{\includegraphics[width=0.22\textwidth]{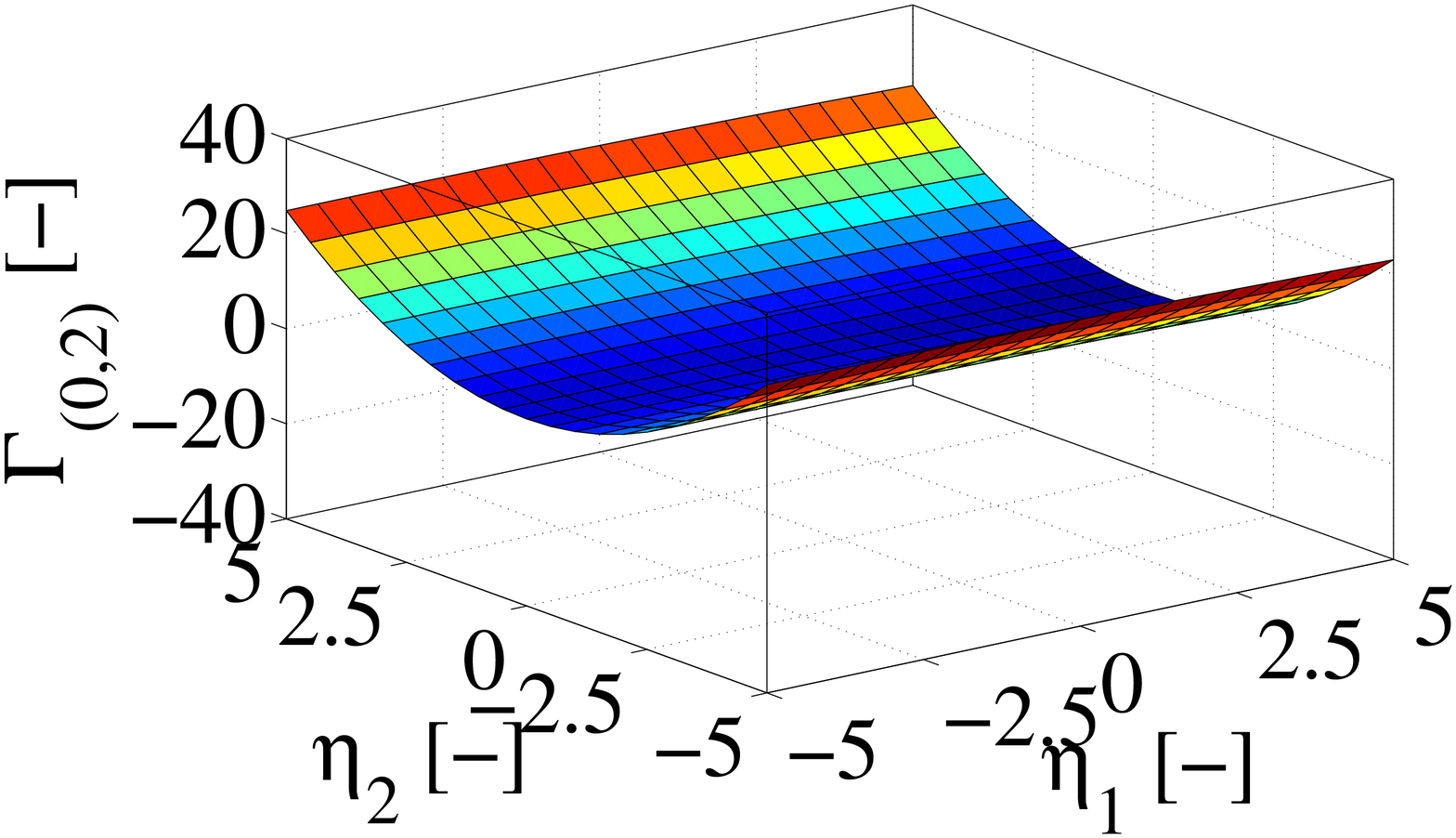}}
    \vfill
    \subfigure[$\Gamma_{30}$.]{\includegraphics[width=0.22\textwidth]{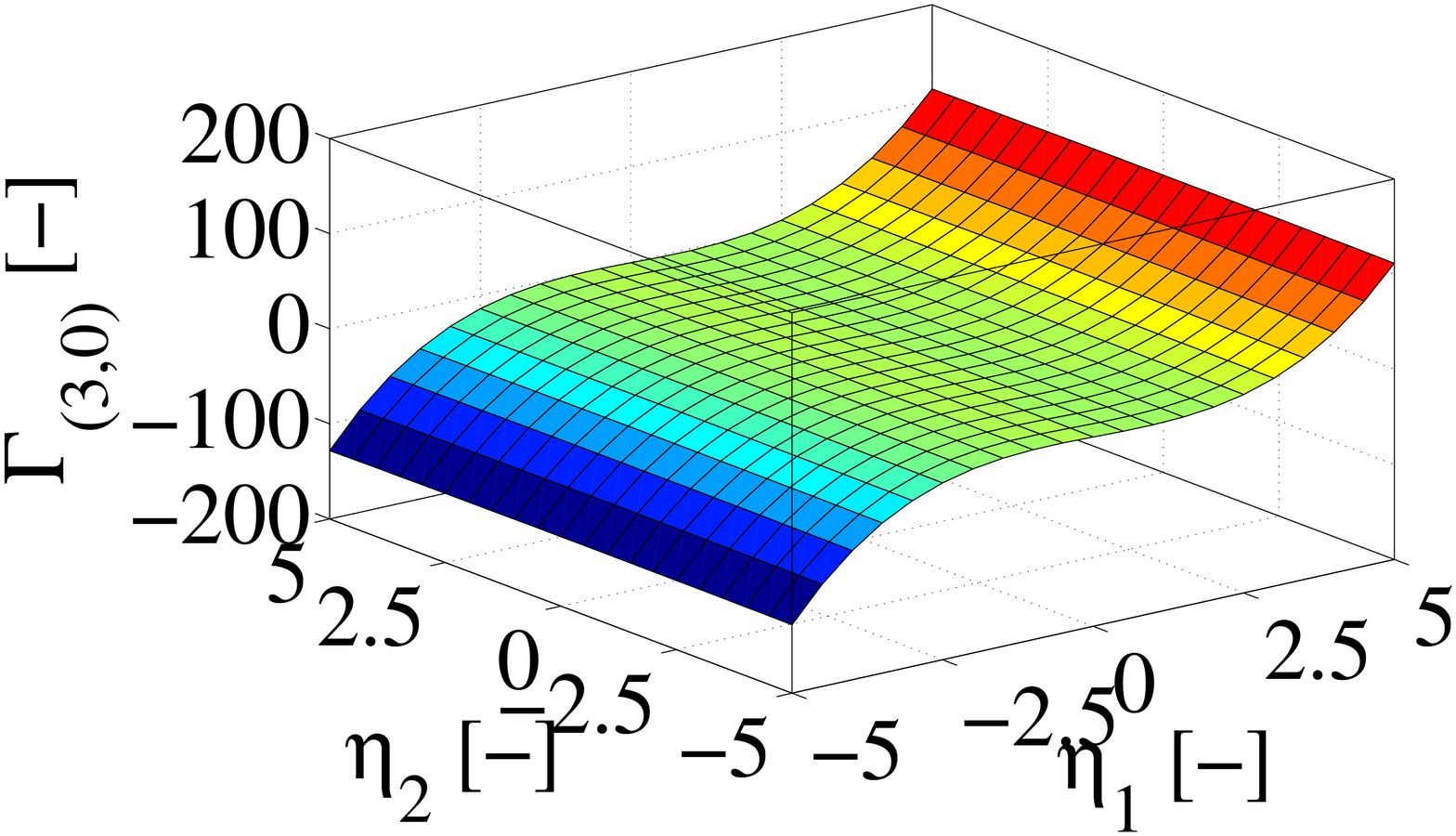}}
    \hspace{1mm}
    \subfigure[$\Gamma_{21}$.]{\includegraphics[width=0.22\textwidth]{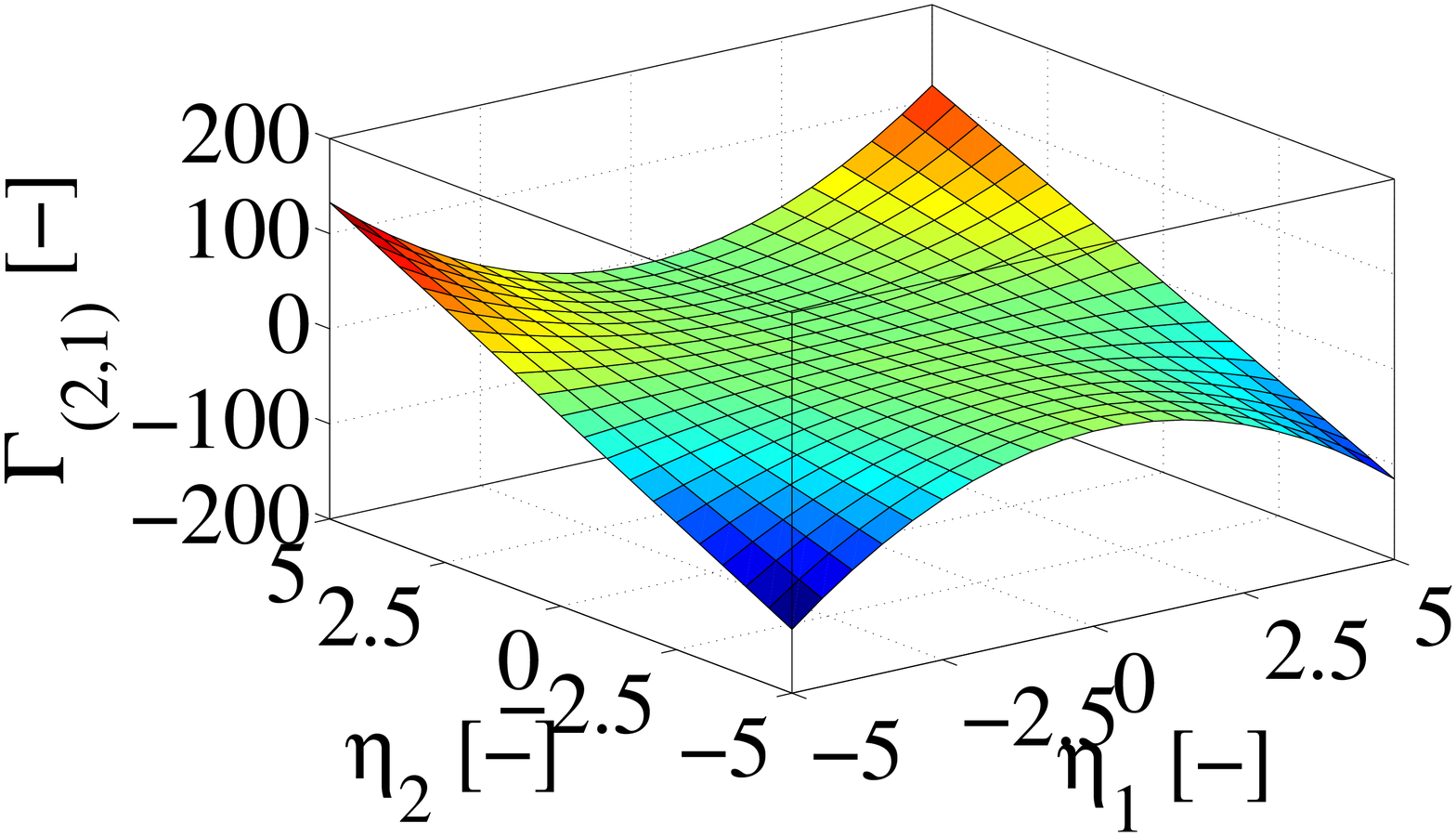}}
    \hspace{1mm}
    \subfigure[$\Gamma_{12}$.]{\includegraphics[width=0.22\textwidth]{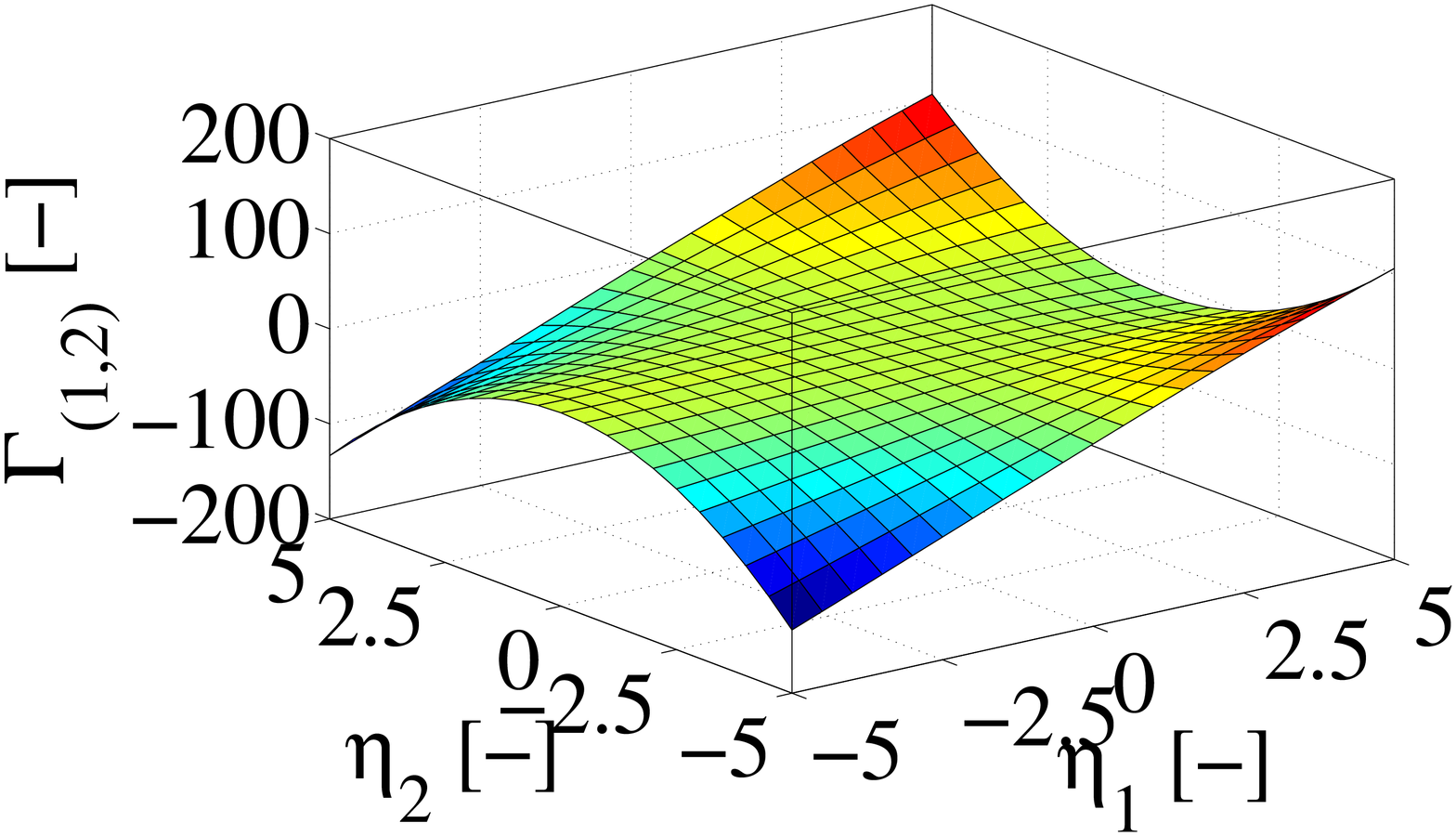}}            
    \hspace{1mm}
    \subfigure[$\Gamma_{03}$.]{\includegraphics[width=0.22\textwidth]{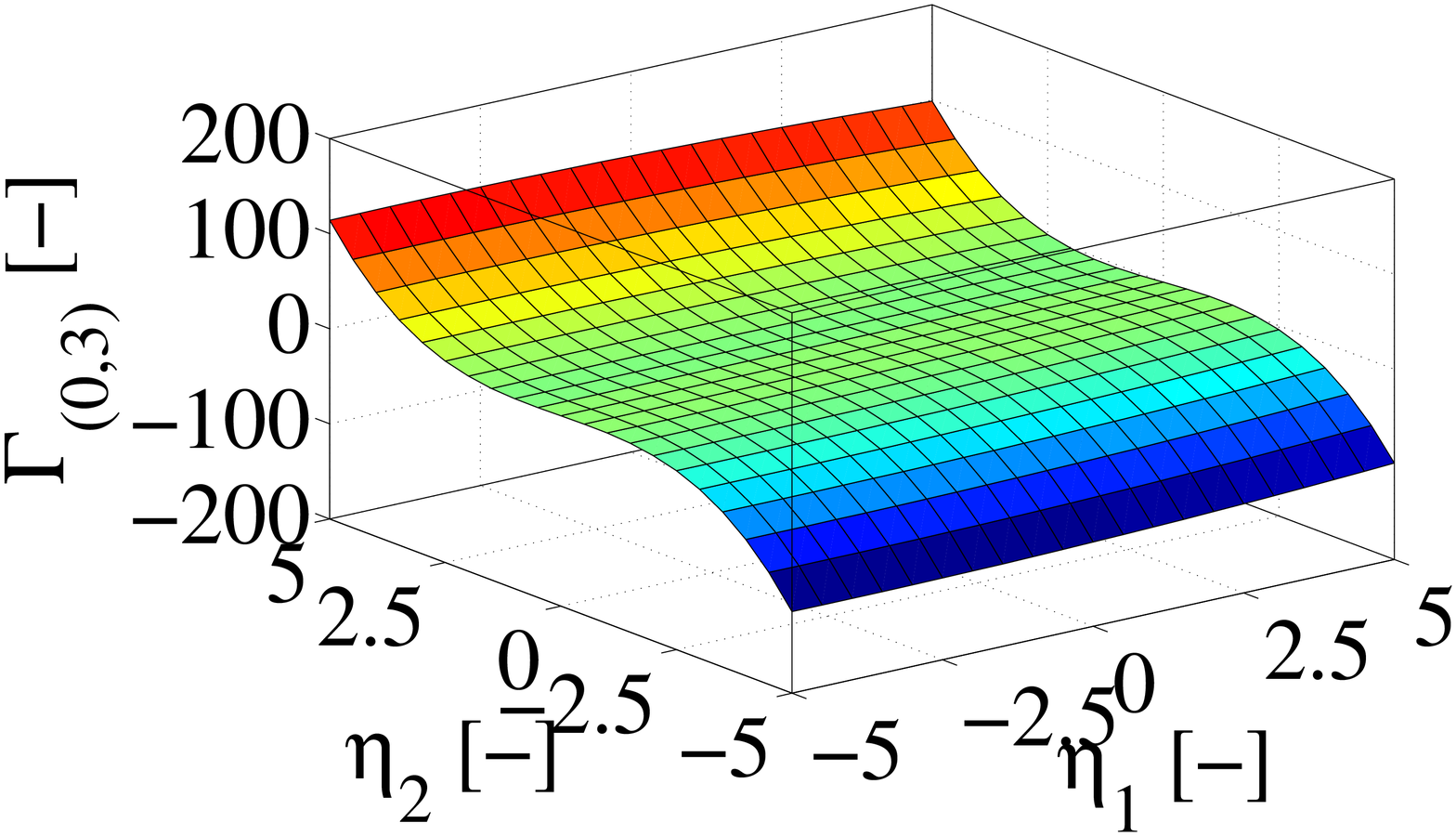}}    
    \caption{PC-based simulation: computed orthonormal polynomials up to a total degree of $q=3$.}\label{fig:figure7}
  \end{center}
\end{figure}

At iteration~$\ell=20$, a KL decomposition retaining only $d=2$ terms was found to be sufficiently accurate to satisfy~(\ref{eq:criterion1}) for~$\epsilon_{1}=0.05$; thus, at iteration~$\ell=20$, the measure transformation necessitated the construction of orthonormal polynomials and quadrature rules with respect to the probability distribution of the first and second reduced random variables.

Figures~\ref{fig:figure6} and~\ref{fig:figure7} illustrate the proposed computational construction of orthonormal polynomials.
With reference to Sec.~\ref{sec:illusimp}, Fig.~\ref{fig:figure6} demonstrates the convergence of the Monte Carlo estimate of an entry of the requisite Gram matrix with respect to the number of samples.  
Figure~\ref{fig:figure7} shows the obtained orthonormal polynomials up to a total degree of~$q=3$.

\begin{figure}[htp]
  \begin{center}
     \includegraphics[width=.8\textwidth]{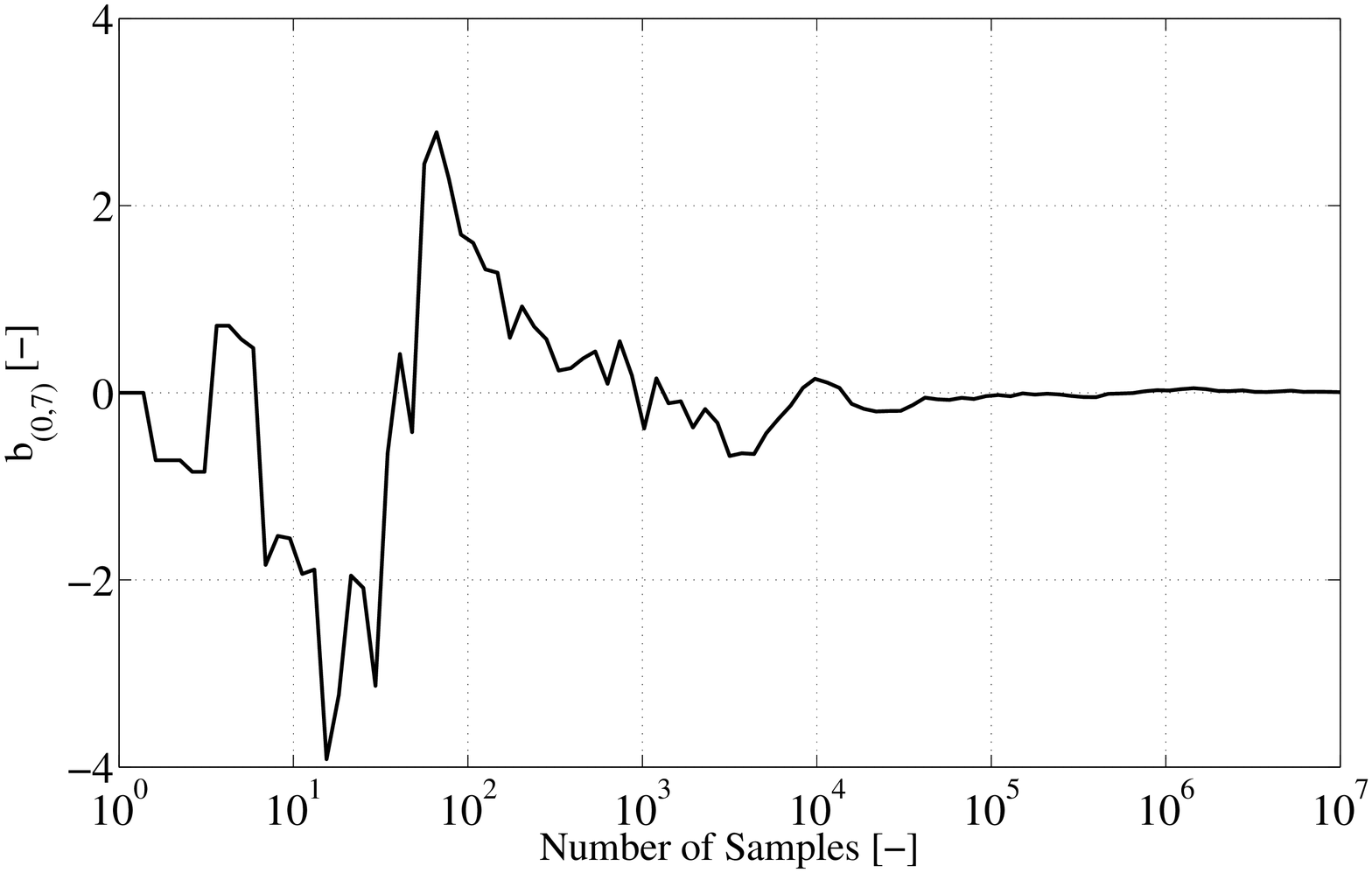}
    \caption{PC-based simulation: convergence of the Monte Carlo estimate of a component of the moment vector with respect to the number of samples.}\label{fig:figure8}
  \end{center}
\end{figure}

\begin{figure}[htp]
  \begin{center}
    \subfigure[$\tilde{\boldsymbol{\eta}}{}_{k}^{\ell}=\boldsymbol{\eta}^{\ell,p}(\tilde{\boldsymbol{\xi}}_{k})$ for $\tilde{\lambda}=3$.]{\includegraphics[width=0.31\textwidth]{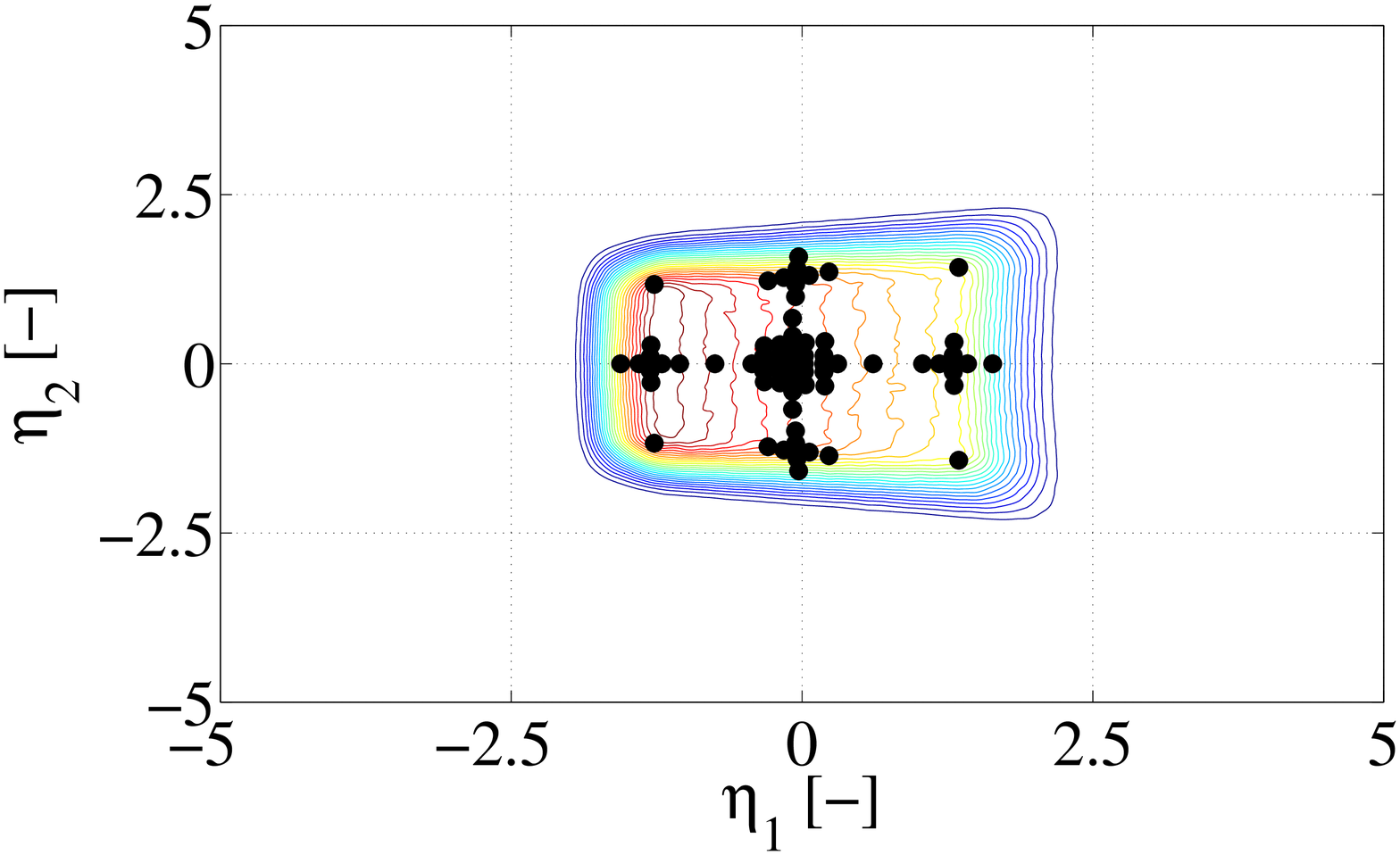}}
    \hfill
    \subfigure[$\boldsymbol{\eta}_{k}^{\ell}$ for $\lambda=4$ and $\tilde{\lambda}=3$.]{\includegraphics[width=0.31\textwidth]{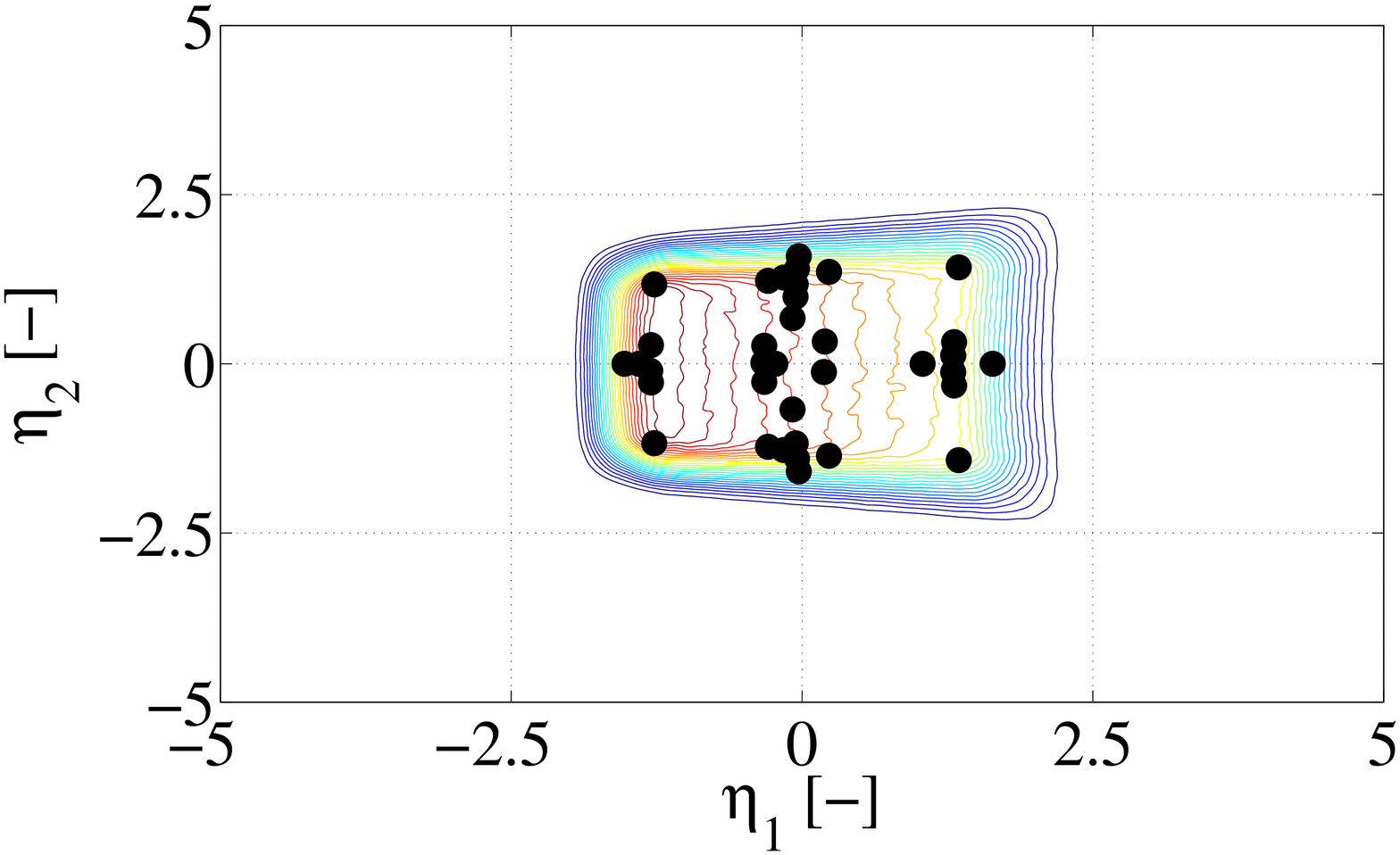}}
    \hfill
    \subfigure[$w^{\ell}$ for $\lambda=4$ and $\tilde{\lambda}=3$.]{\includegraphics[width=0.31\textwidth]{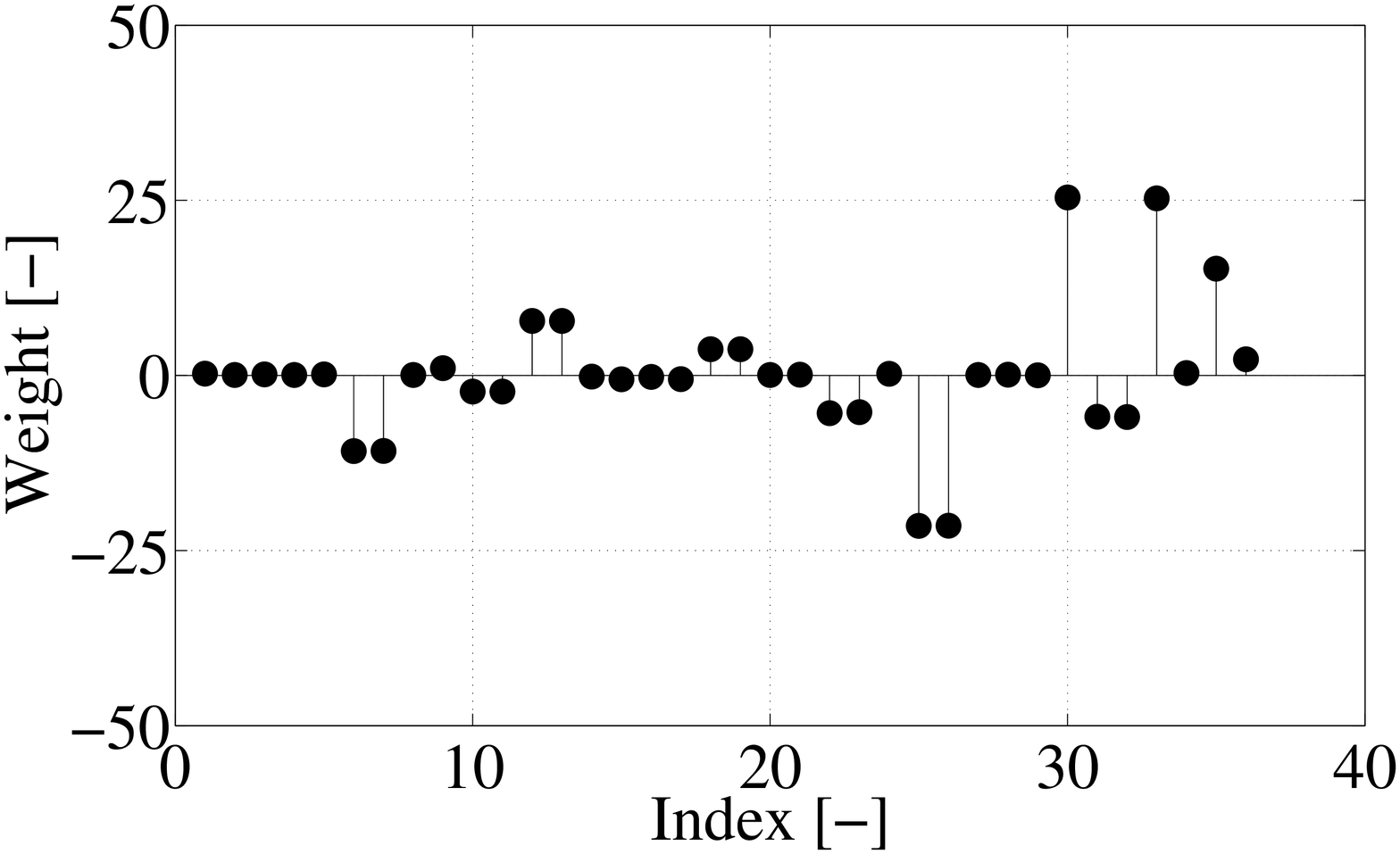}}
    \subfigure[$\tilde{\boldsymbol{\eta}}{}_{k}^{\ell}=\boldsymbol{\eta}^{\ell,p}(\tilde{\boldsymbol{\xi}}_{k})$ for $\tilde{\lambda}=4$.]{\includegraphics[width=0.31\textwidth]{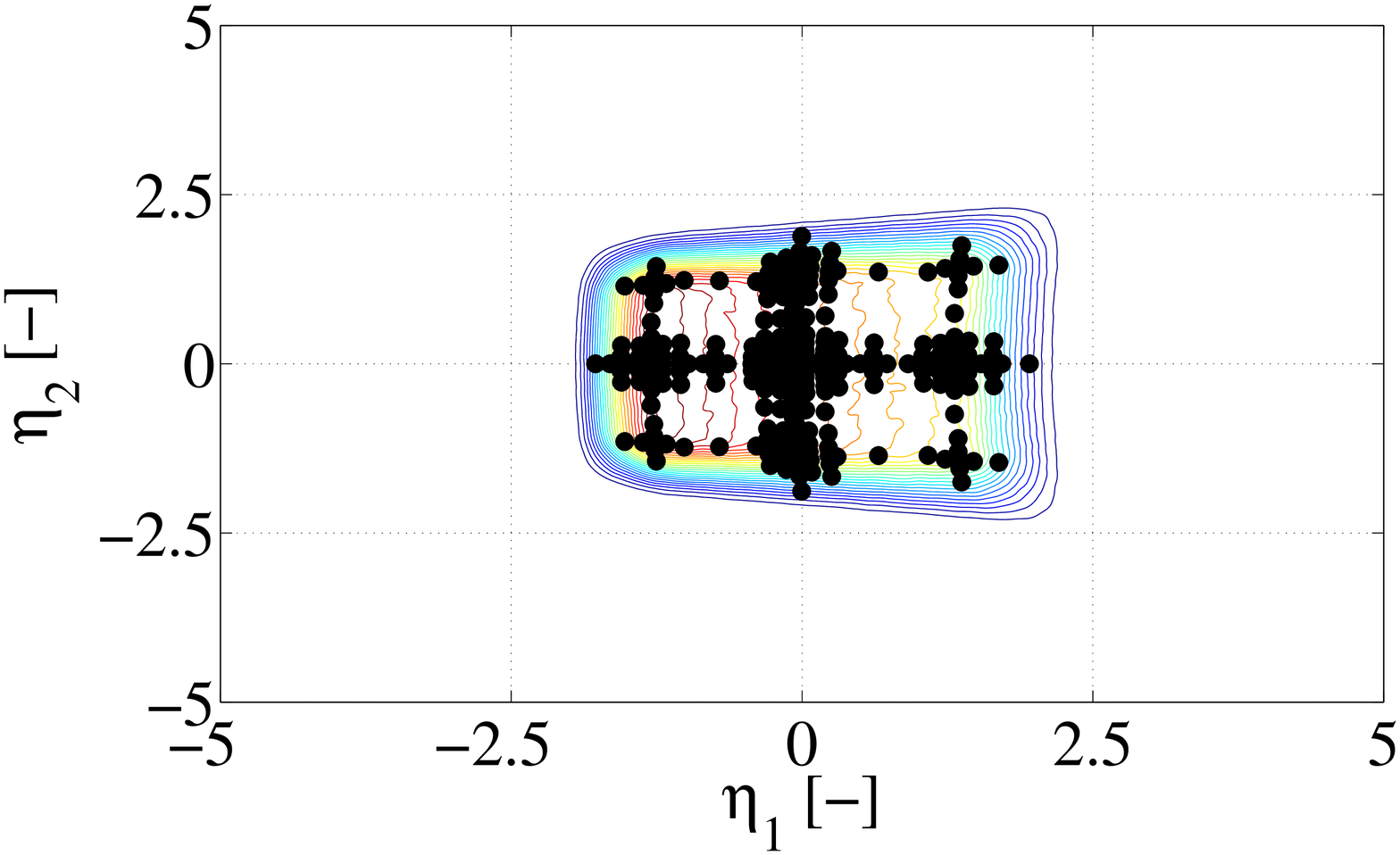}}
    \hfill
    \subfigure[$\boldsymbol{\eta}_{k}^{\ell}$ for $\lambda=4$ and $\tilde{\lambda}=4$.]{\includegraphics[width=0.31\textwidth]{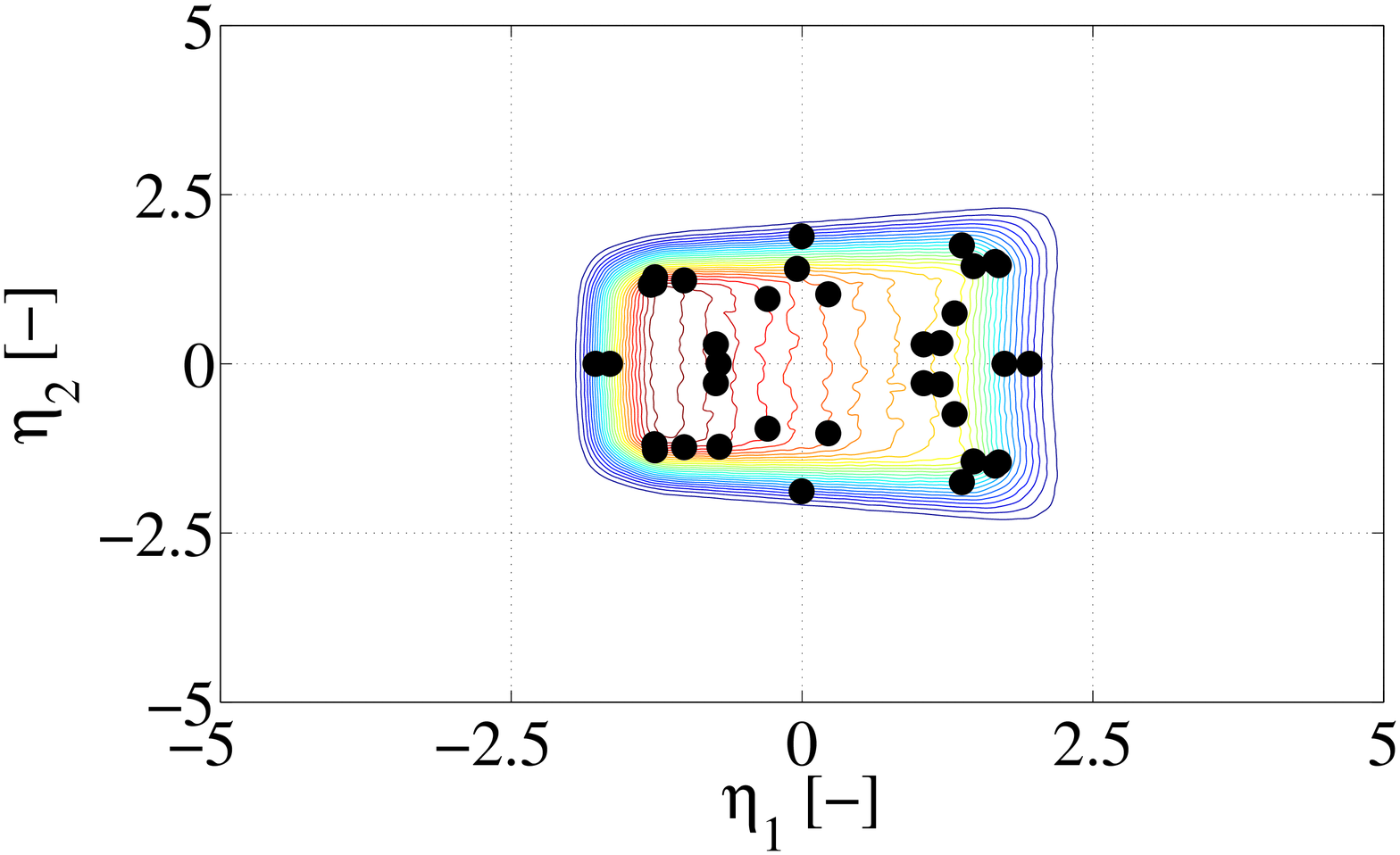}}
    \hfill
    \subfigure[$w^{\ell}$ for $\lambda=4$ and $\tilde{\lambda}=4$.]{\includegraphics[width=0.31\textwidth]{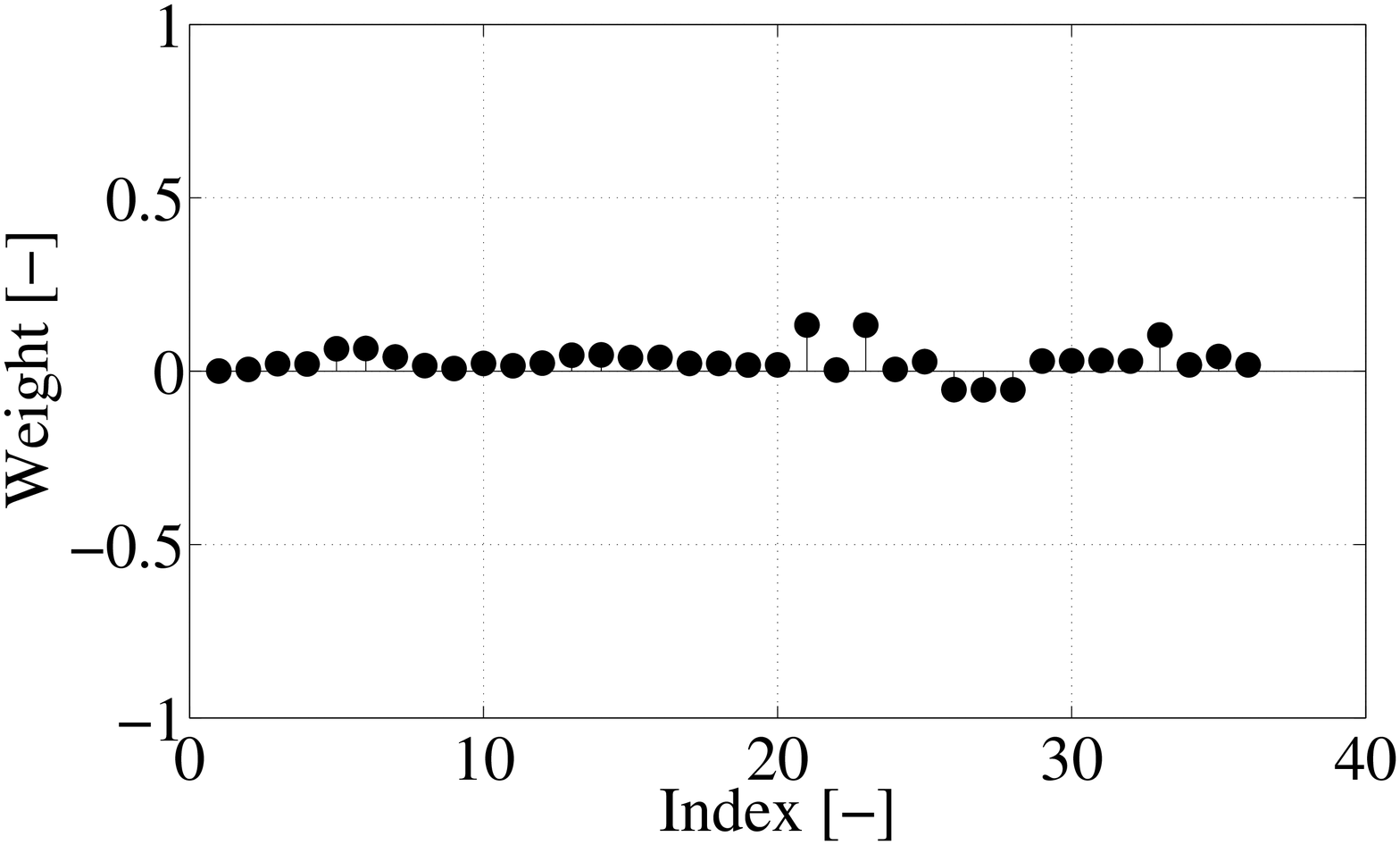}}
    \subfigure[$\tilde{\boldsymbol{\eta}}{}_{k}^{\ell}=\boldsymbol{\eta}^{\ell,p}(\tilde{\boldsymbol{\xi}}_{k})$ for $\tilde{\lambda}=5$.]{\includegraphics[width=0.31\textwidth]{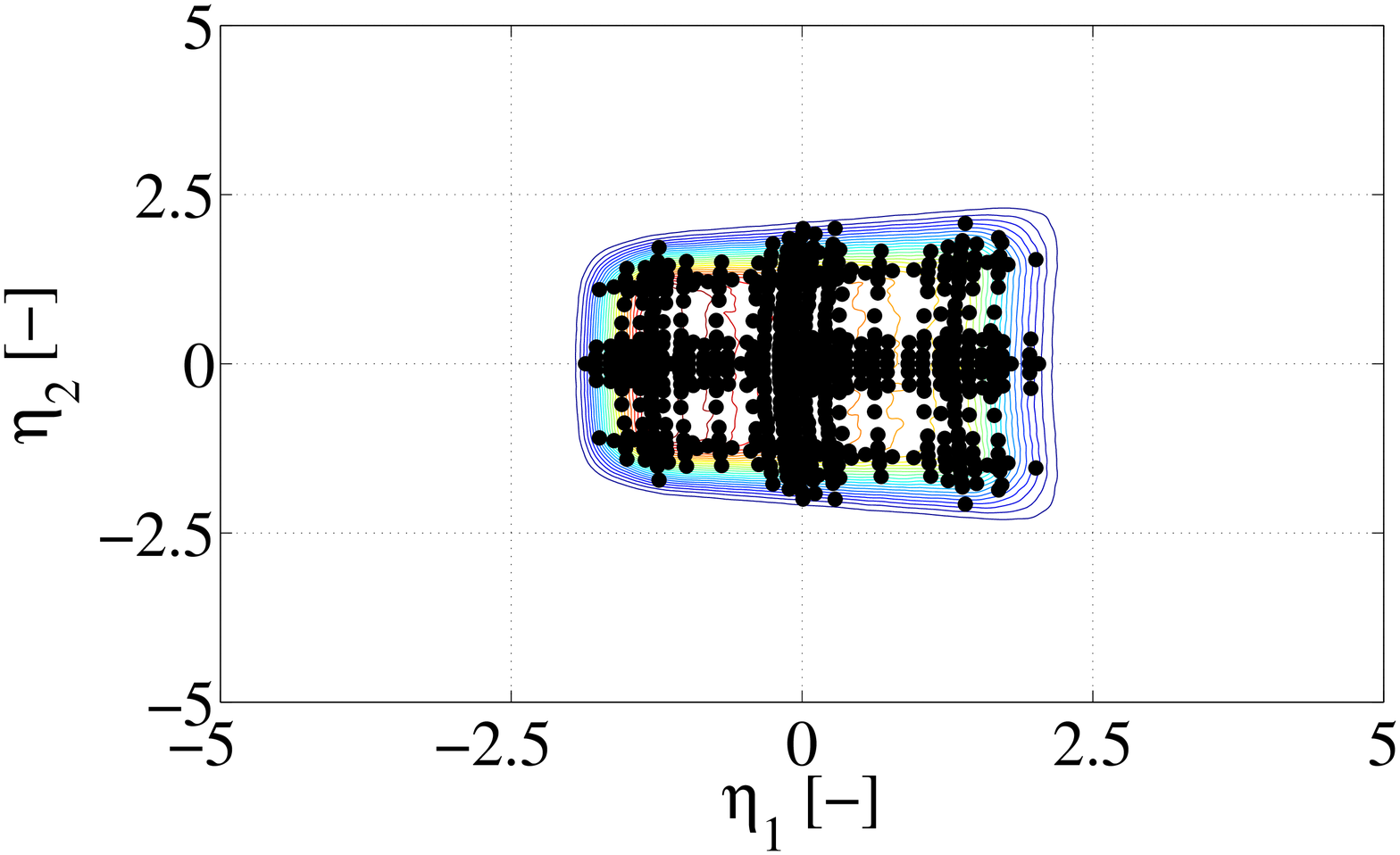}}
    \hfill
    \subfigure[$\boldsymbol{\eta}_{k}^{\ell}$ for $\lambda=4$ and $\tilde{\lambda}=5$.]{\includegraphics[width=0.31\textwidth]{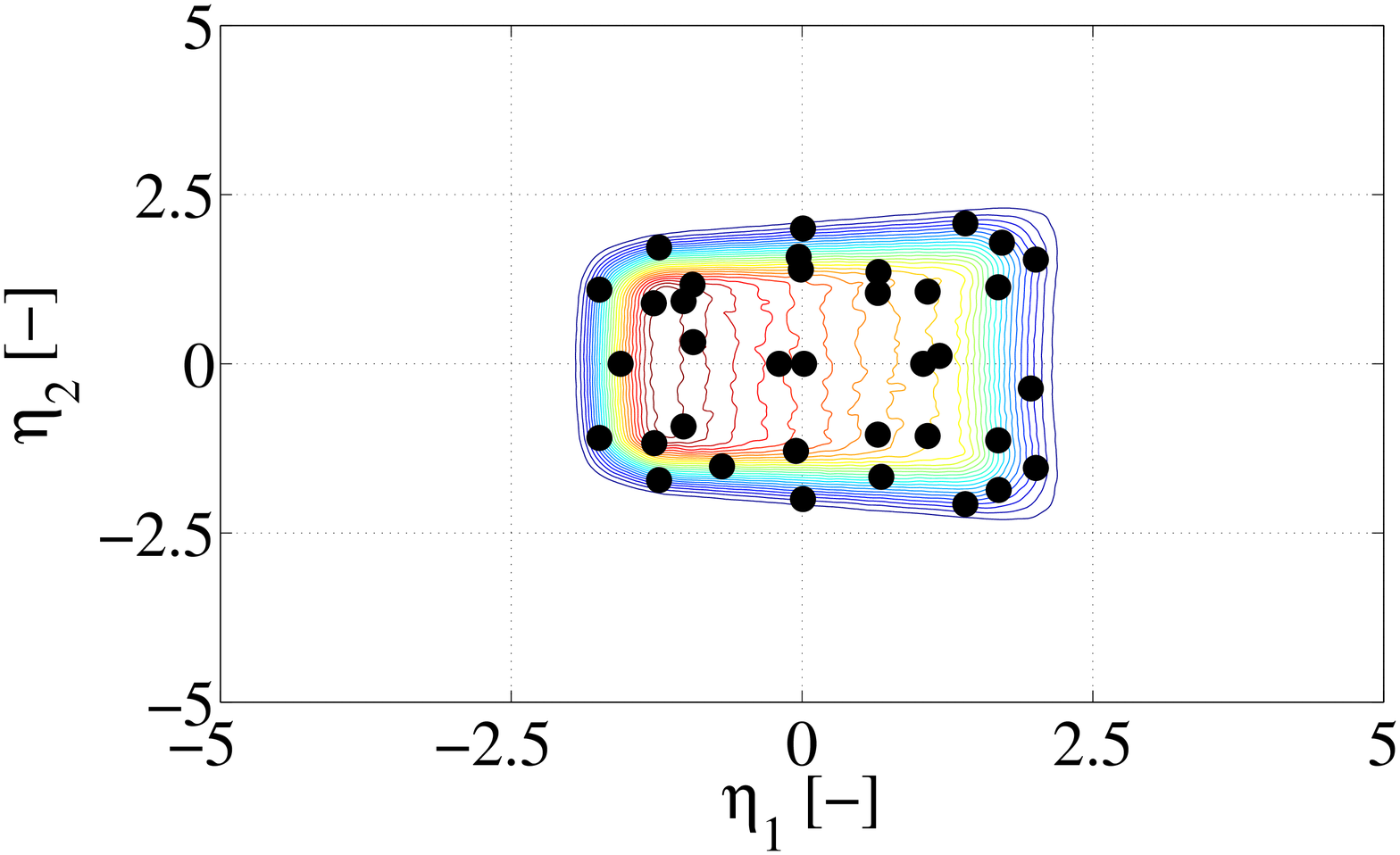}}
    \hfill
    \subfigure[$w^{\ell}$ for $\lambda=4$ and $\tilde{\lambda}=5$.]{\includegraphics[width=0.31\textwidth]{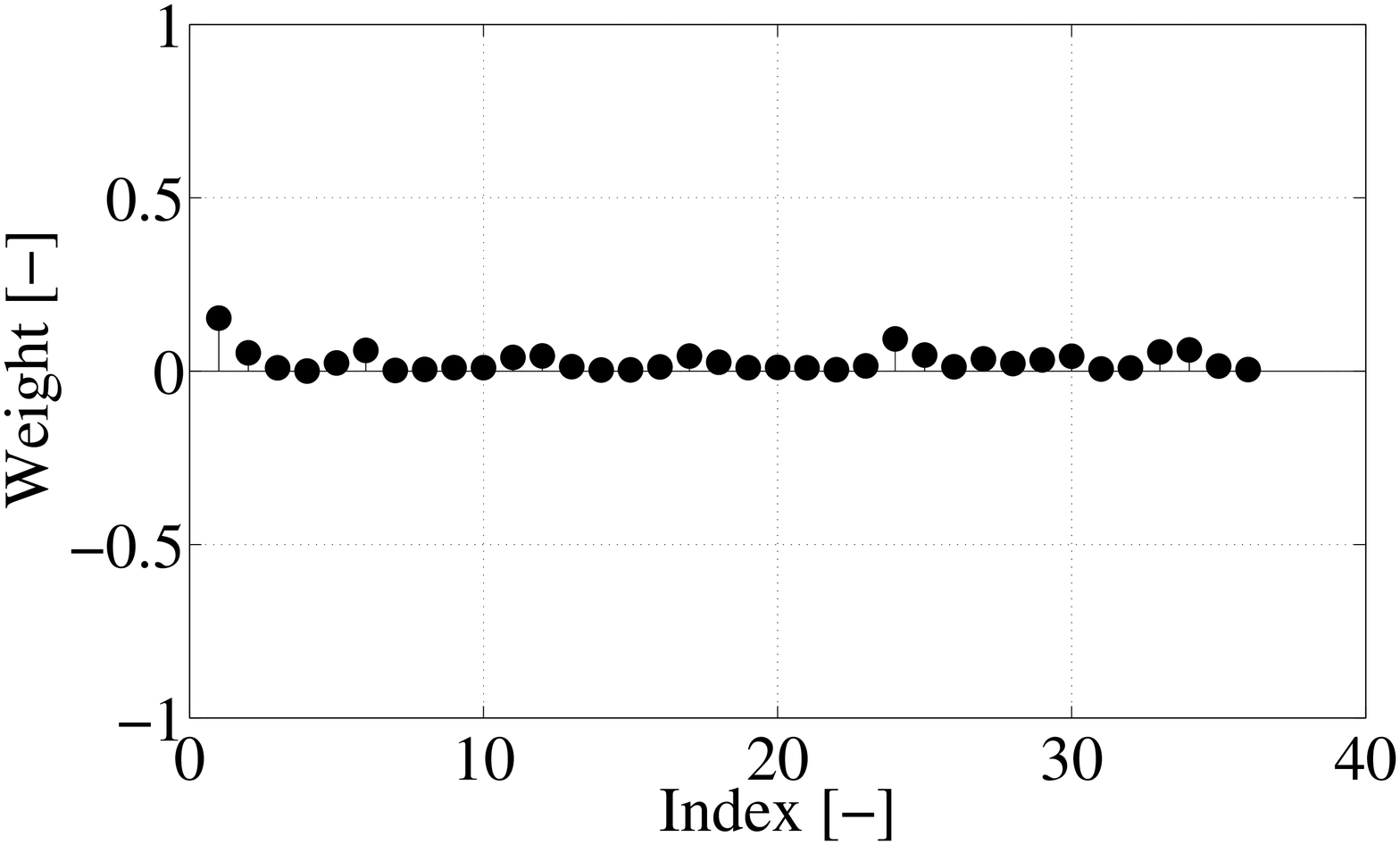}}
    \caption{PC-based simulation: nodes of the quadrature rules~$\{(\tilde{\boldsymbol{\eta}}{}^{\ell}_{k}=\boldsymbol{\eta}^{\ell,p}(\tilde{\boldsymbol{\xi}}_{k}),\tilde{w}_{k}),1\leq k\leq\tilde{\nu}\}$ and nodes and weights of the embedded quadrature rules~$\{(\boldsymbol{\eta}_{k}^{\ell},w_{k}),1\leq k\leq\nu^{\ell}\}$ of level $\lambda=4$ obtained from sparse-grid Gauss-Legendre quadrature rules~$\{(\tilde{\boldsymbol{\xi}}_{k},\tilde{w}_{k}),1\leq k\leq\tilde{\nu}\}$ of levels $\tilde{\lambda}=3$, $4$, and $5$.}\label{fig:figure9}
  \end{center}
\end{figure}

\begin{figure}[htp]
  \begin{center}
     \includegraphics[width=.8\textwidth]{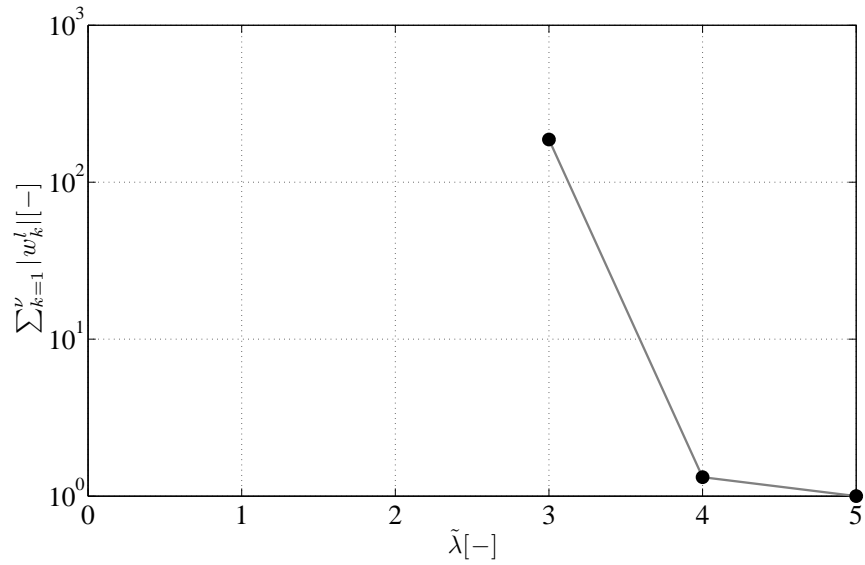}
    \caption{PC-based simulation: convergence of the sum of the absolute values of the weights of the embedded quadrature rule of level $\lambda=4$ with respect to the level~$\tilde{\lambda}$ of the sparse-grid Gauss-Legendre quadrature rule.}\label{fig:figure9j}
  \end{center}
\end{figure}

\begin{figure}[htp]
  \begin{center}
    \subfigure[$\boldsymbol{\eta}_{k}^{\ell}$ for $\lambda=1$.]{\includegraphics[width=0.49\textwidth]{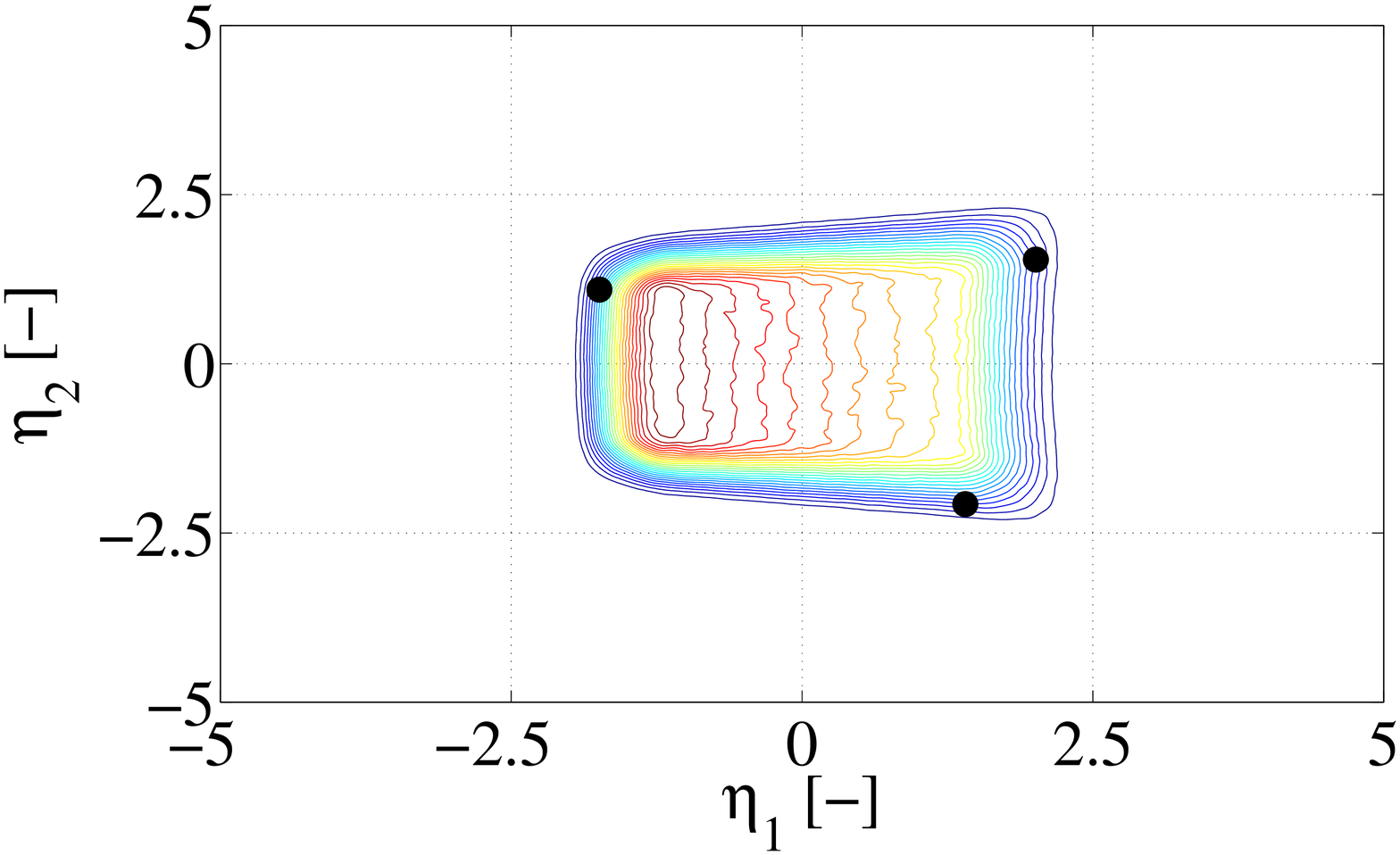}}
    \hfill
    \subfigure[$w_{k}^{\ell}$ for $\lambda=1$.]{\includegraphics[width=0.49\textwidth]{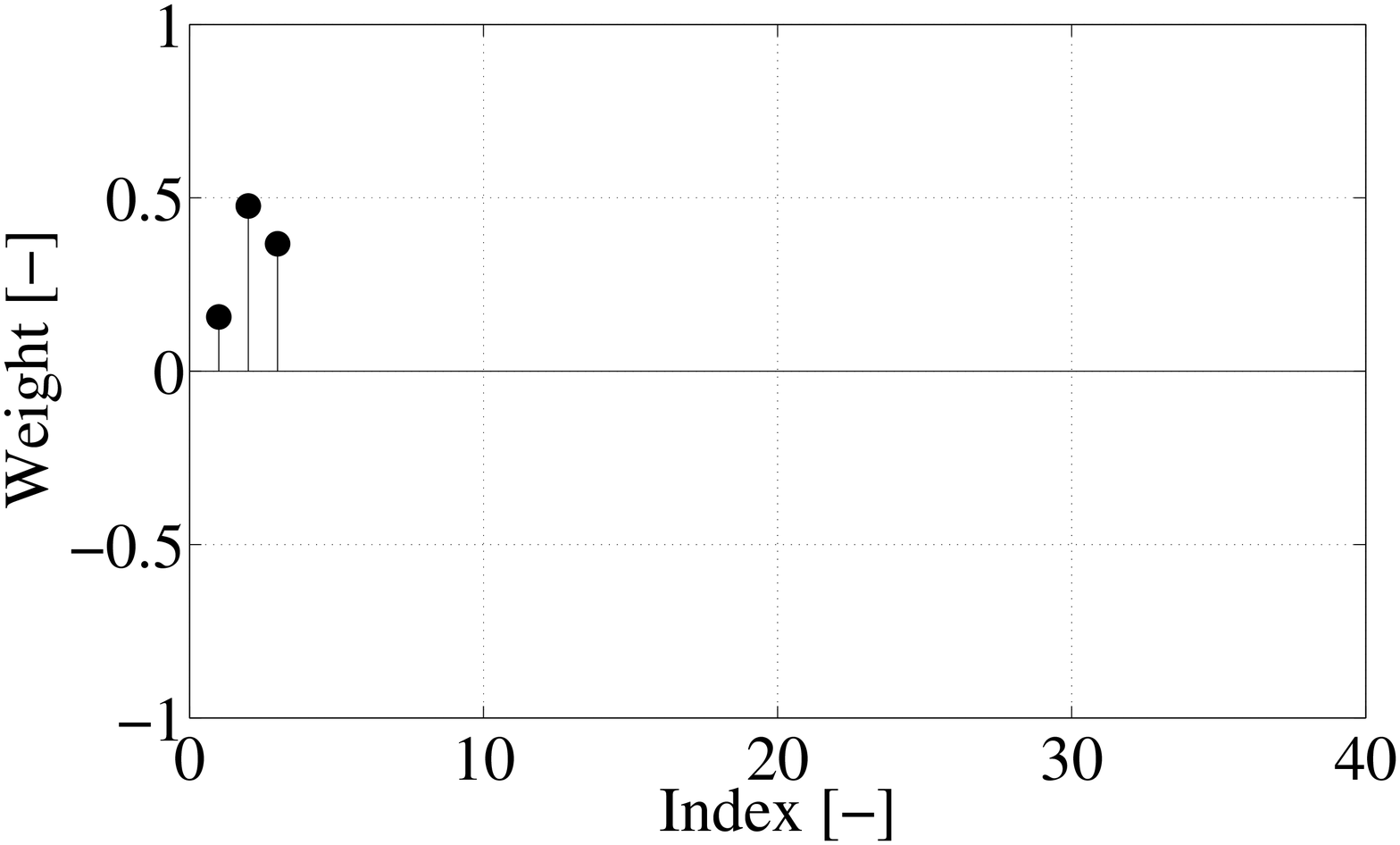}}
    \subfigure[$\boldsymbol{\eta}_{k}^{\ell}$ for $\lambda=2$.]{\includegraphics[width=0.49\textwidth]{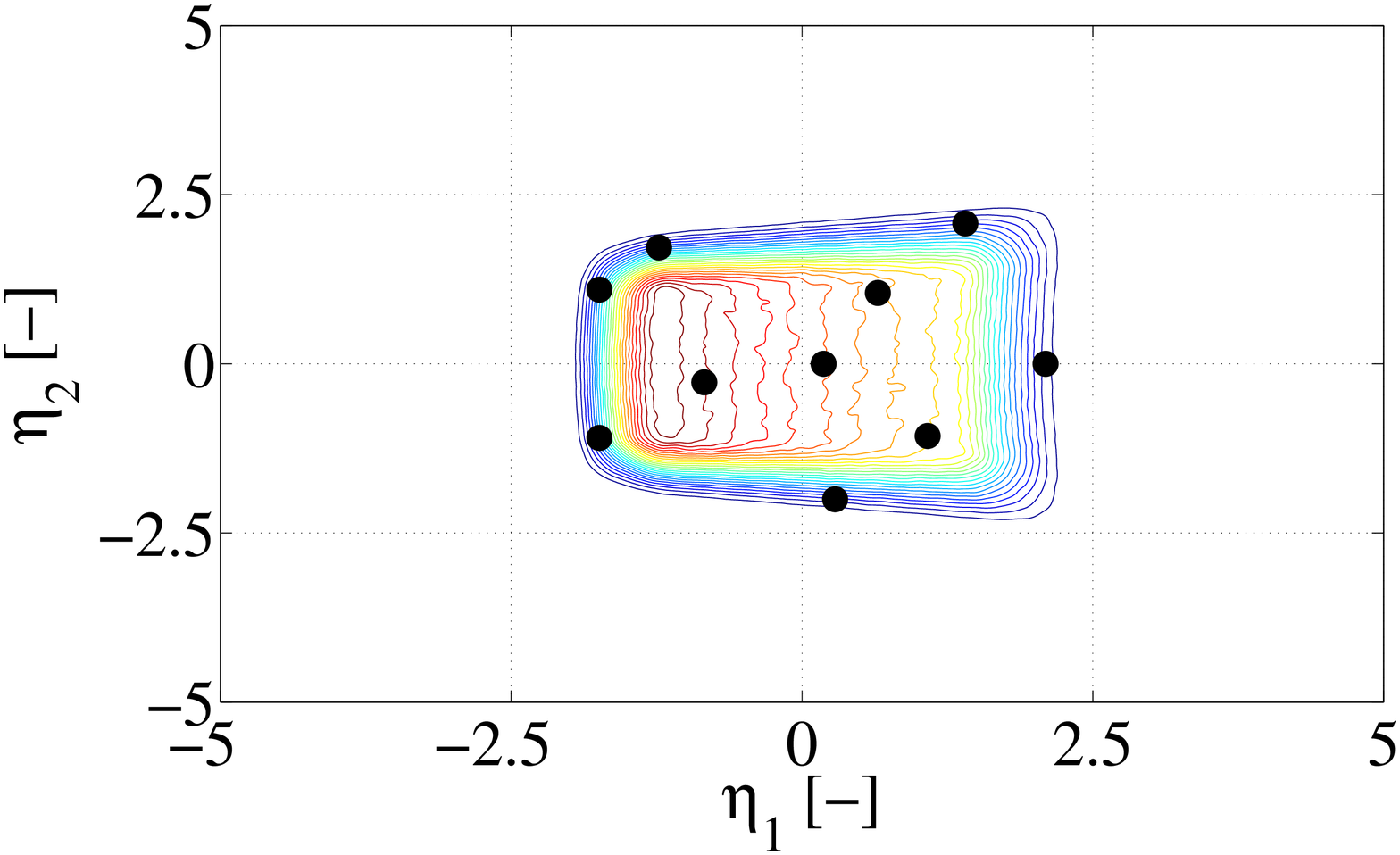}}
    \hfill
    \subfigure[$w_{k}^{\ell}$ for $\lambda=2$.]{\includegraphics[width=0.49\textwidth]{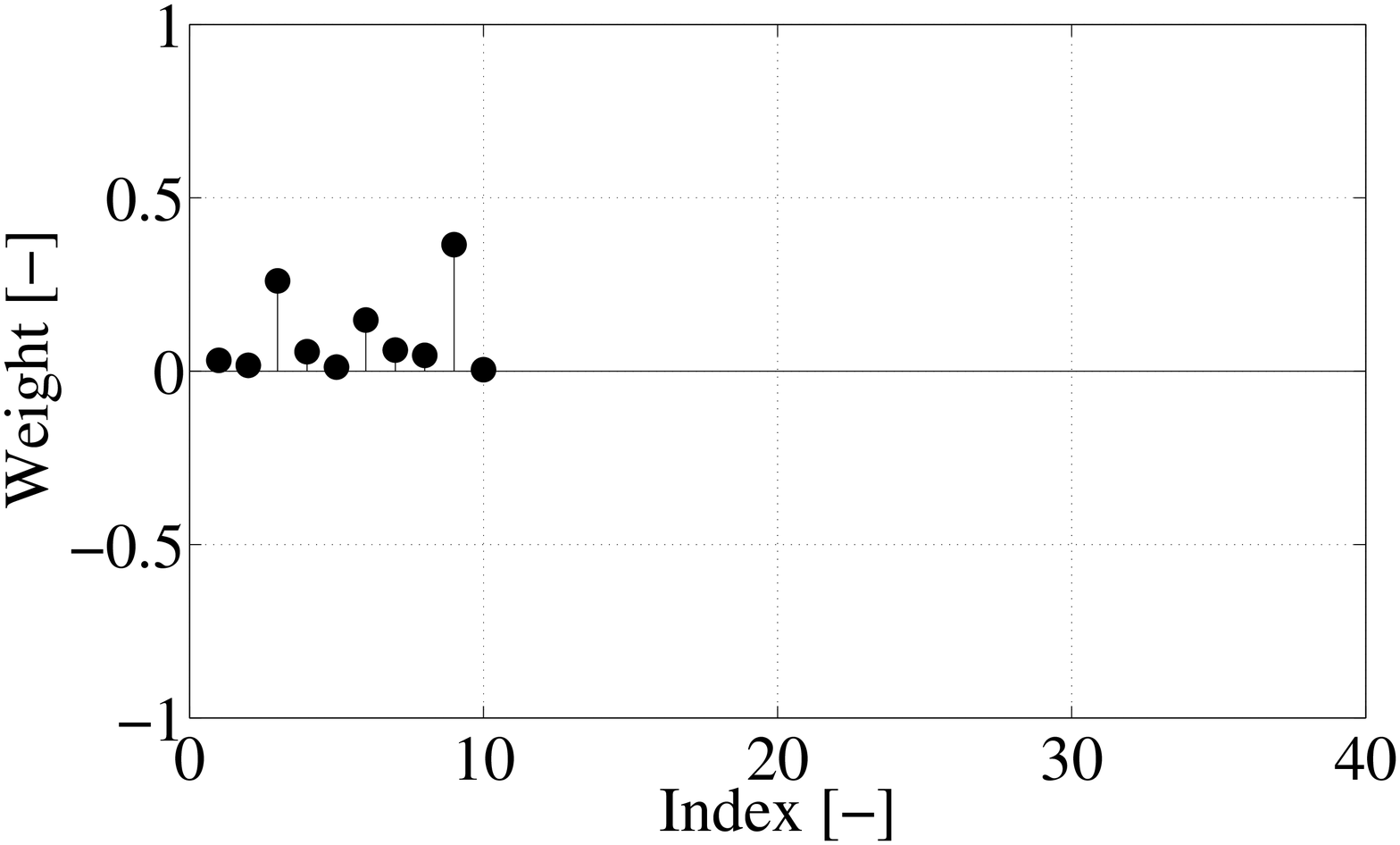}}
    \subfigure[$\boldsymbol{\eta}_{k}^{\ell}$ for $\lambda=3$.]{\includegraphics[width=0.49\textwidth]{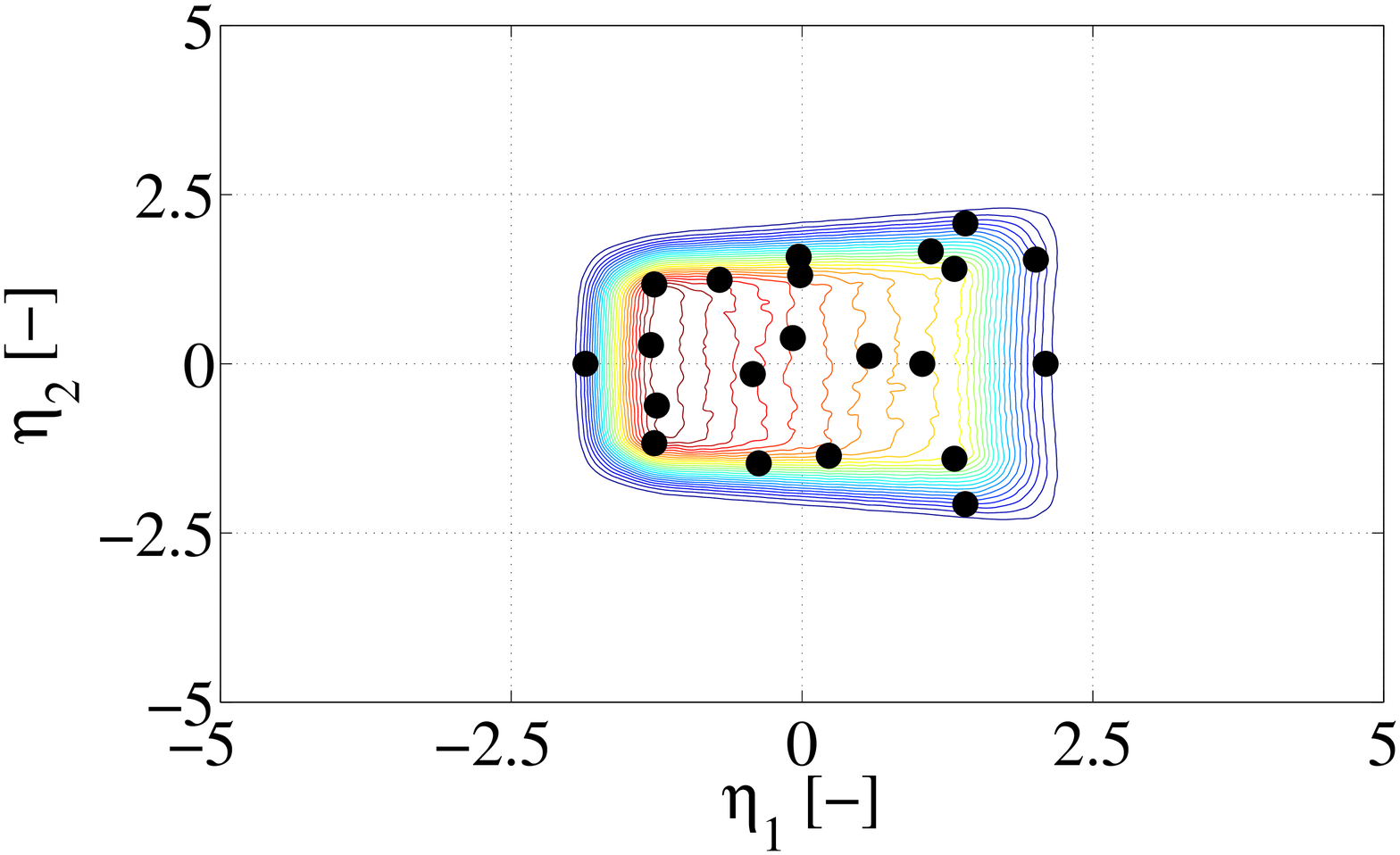}}
    \hfill
    \subfigure[$w_{k}^{\ell}$ for $\lambda=3$.]{\includegraphics[width=0.49\textwidth]{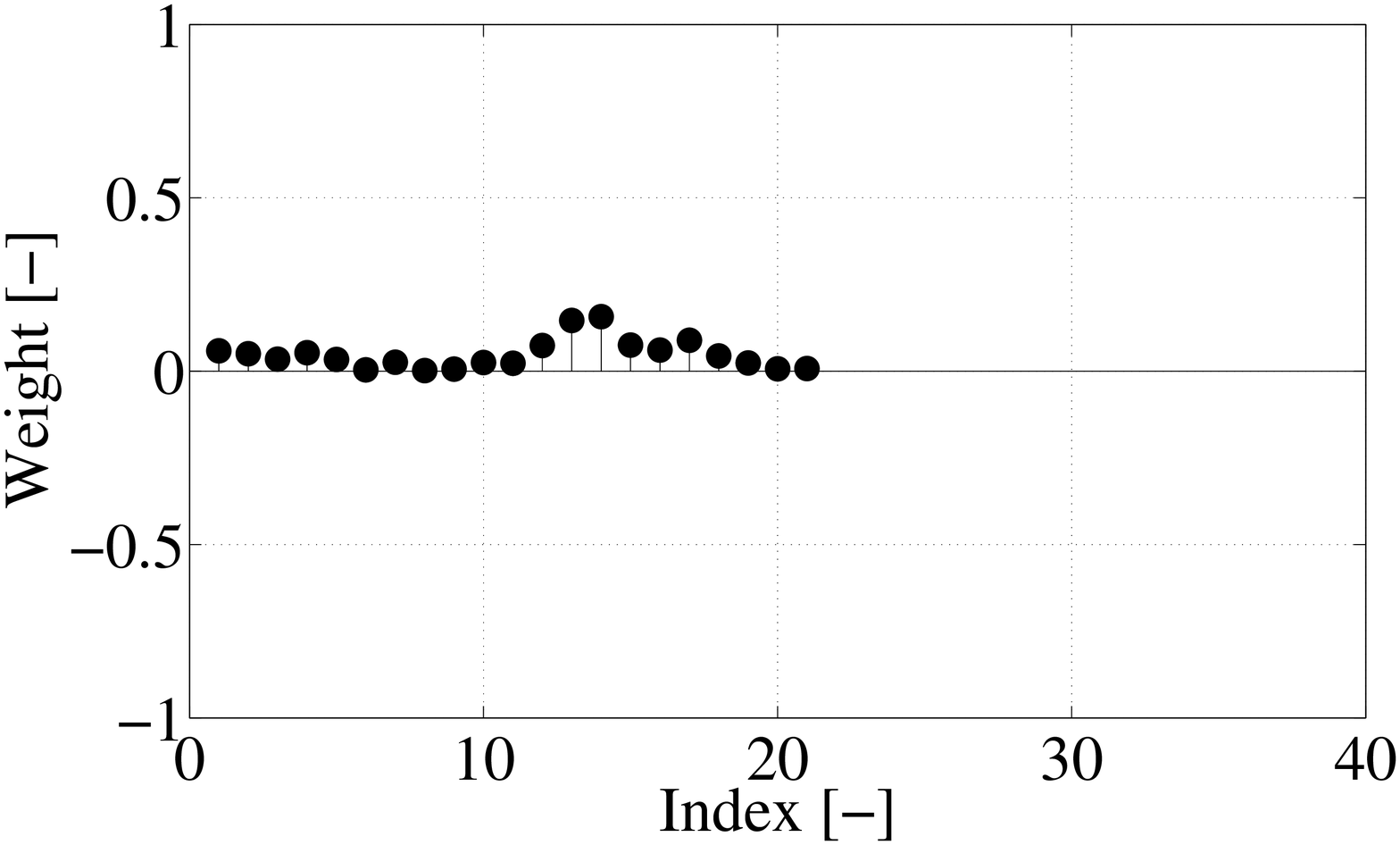}}
    \caption{PC-based simulation: embedded quadrature rules~$\{(\boldsymbol{\eta}_{k}^{\ell},w_{k}),1\leq k\leq\nu^{\ell}\}$ of levels $\lambda=1$, $2$, and $3$ obtained from the sparse-grid Gauss-Legendre quadrature rule of level $\tilde{\lambda}=5$.}\label{fig:figure10}
  \end{center}
\end{figure}

Figures~\ref{fig:figure8}--\ref{fig:figure10} illustrate the proposed computational construction of quadrature rules.
With reference to Sec.~\ref{sec:illusimp}, Fig.~\ref{fig:figure8} demonstrates the convergence of the Monte Carlo estimate of a component of the requisite moment vector with respect to the number of samples.  

We obtained the results shown in Fig.~\ref{fig:figure9} as follows.
First, we constructed three sparse-grid Gauss-Legendre quadrature rules of the form~$\{(\tilde{\boldsymbol{\xi}}_{k},\tilde{w}_{k}),1\leq k\leq\tilde{\nu}\}$ for integration with respect to the probability distribution of the input random variables: the first rule was of the level~$\tilde{\lambda}=3$ and had~$\tilde{\nu}=261$ nodes; the second was of the level $\tilde{\lambda}=4$ and had $\tilde{\nu}=2,441$ nodes; and the third was of the level~$\tilde{\lambda}=5$ and had~$\tilde{\nu}=18,881$ nodes.
Next, we carried out a ``change of variables" to transform each of these sparse-grid Gauss-Legendre quadrature rules into a corresponding quadrature rule of the form~$\{(\tilde{\boldsymbol{\eta}}{}^{\ell}_{k}=\boldsymbol{\eta}^{\ell,p}(\tilde{\boldsymbol{\xi}}_{k}),\tilde{w}_{k}),1\leq k\leq\tilde{\nu}\}$ for integration with respect to the probability distribution of the reduced random variables; Fig.~\ref{fig:figure9}((a), (d), and(g)) shows the nodes of the quadrature rules thus obtained.
Finally, on the basis of each of these quadrature rules, we constructed an embedded quadrature rule~$\{(\boldsymbol{\eta}_{k}^{\ell},w_{k}),1\leq k\leq\nu^{\ell}\}$ of level~$\lambda=4$ using Algorithm~1; Fig.~\ref{fig:figure9}((b), (e), and(h)) and~\ref{fig:figure9}((c), (f), and(i)) shows the nodes and weights of the embedded quadrature rules thus obtained.
Figure~\ref{fig:figure9}((c), (f), and(i)) indicates that as the level~$\tilde{\lambda}$ was increased, Algorithm~1 was provided with a greater choice of candidate nodes from which the nodes of the embedded quadrature rule could be selected and thus provided a better embedded quadrature rule with a smaller value for the sum of the absolute values of the weights.
Figure~\ref{fig:figure9j} demonstrates the convergence of the sum of the absolute values of the weights with respect to the level~$\tilde{\lambda}$, indicating convergence at~$\tilde{\lambda}=5$. 

Figure~\ref{fig:figure10} shows the nodes and weights of the embedded quadrature rules of levels~$\lambda=1$, $2$, and $3$ constructed similarly on the basis of the sparse-grid Gauss-Legendre quadrature rule of level $\tilde{\lambda}=5$.
This figure indicates that as the level $\lambda$ was increased, higher accuracy was required, and thus an embedded quadrature rule with more nodes was systematically obtained.

\begin{figure}[htp]
  \begin{center}
   \subfigure[Temperature.]{\includegraphics[width=.8\textwidth]{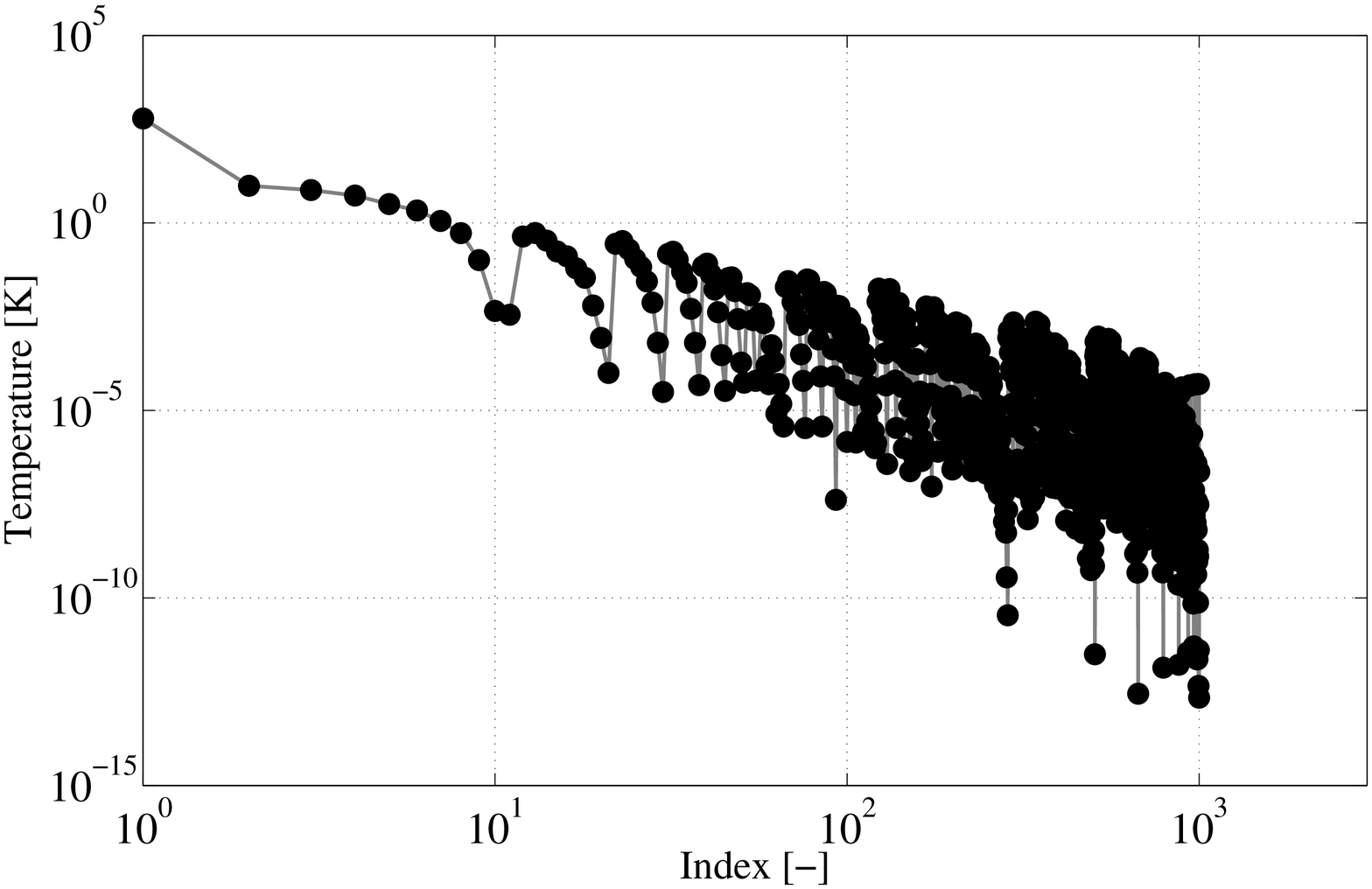}}
    \hfill
    \subfigure[Neutron flux.]{\includegraphics[width=.8\textwidth]{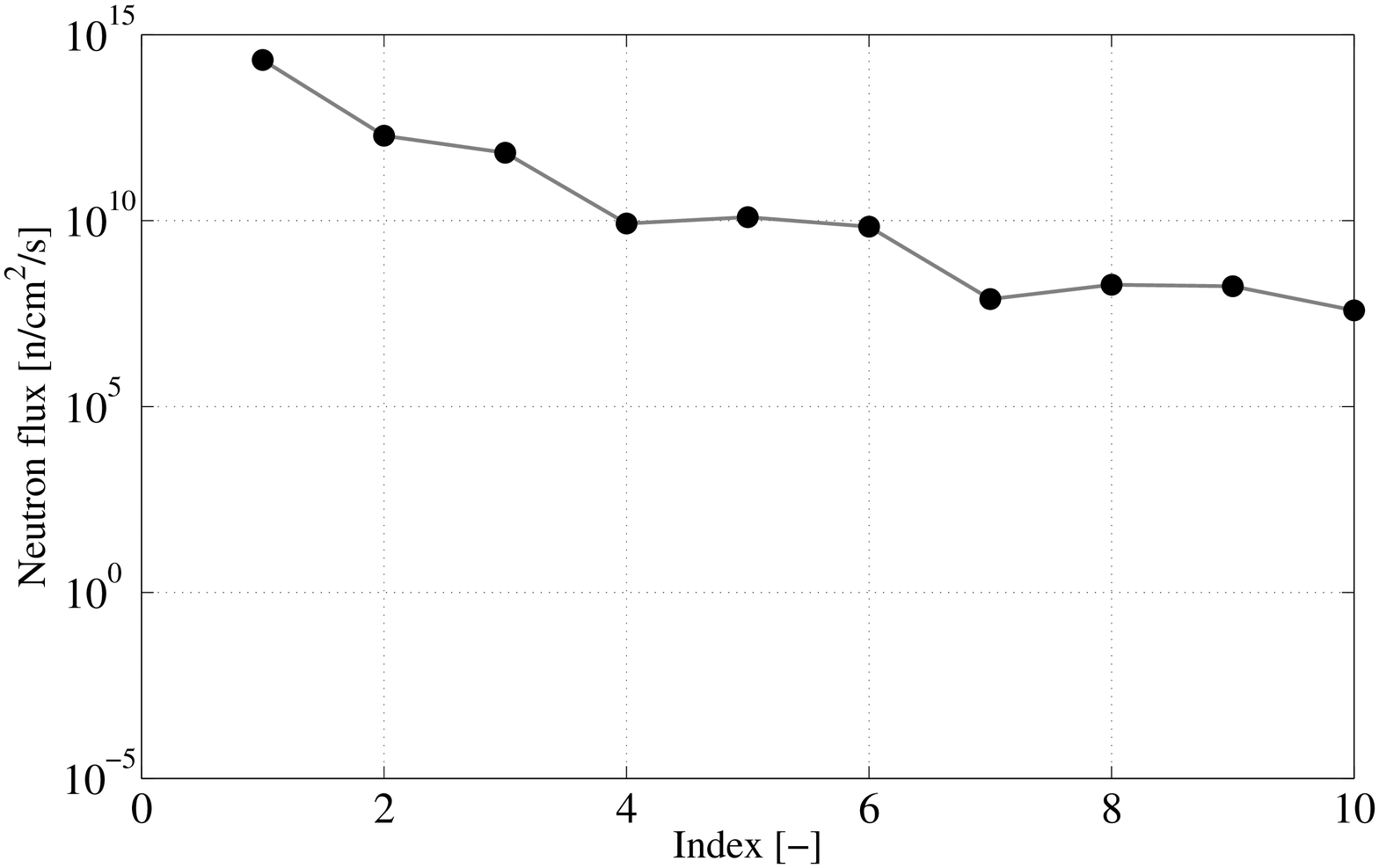}}
    \caption{PC-based simulation: PC coordinates at $x=10\,[\mbox{cm}]$ of the solution.}\label{fig:figure11}
  \end{center}
\end{figure}

At iteration~$\ell=20$, a PC expansion truncated at a total degree of $q=3$ was sufficiently accurate to satisfy~(\ref{eq:criterion2}) for~$\epsilon_{2}= 0.000001$.
The representation of the random temperature by a PC expansion of dimension $m=10$ and total degree $p=4$ requires~$1,001=14!/10!/4!$ terms, whereas the representation of the random neutron flux by a PC expansion of dimension~$d=2$ and total degree~$q=3$ requires only~$10=5!/2!/3!$ terms.
Figure~\ref{fig:figure11} shows a few PC coordinates.

\begin{figure}[htp]
  \begin{center}
    \subfigure[Temperature.]{\includegraphics[width=0.8\textwidth]{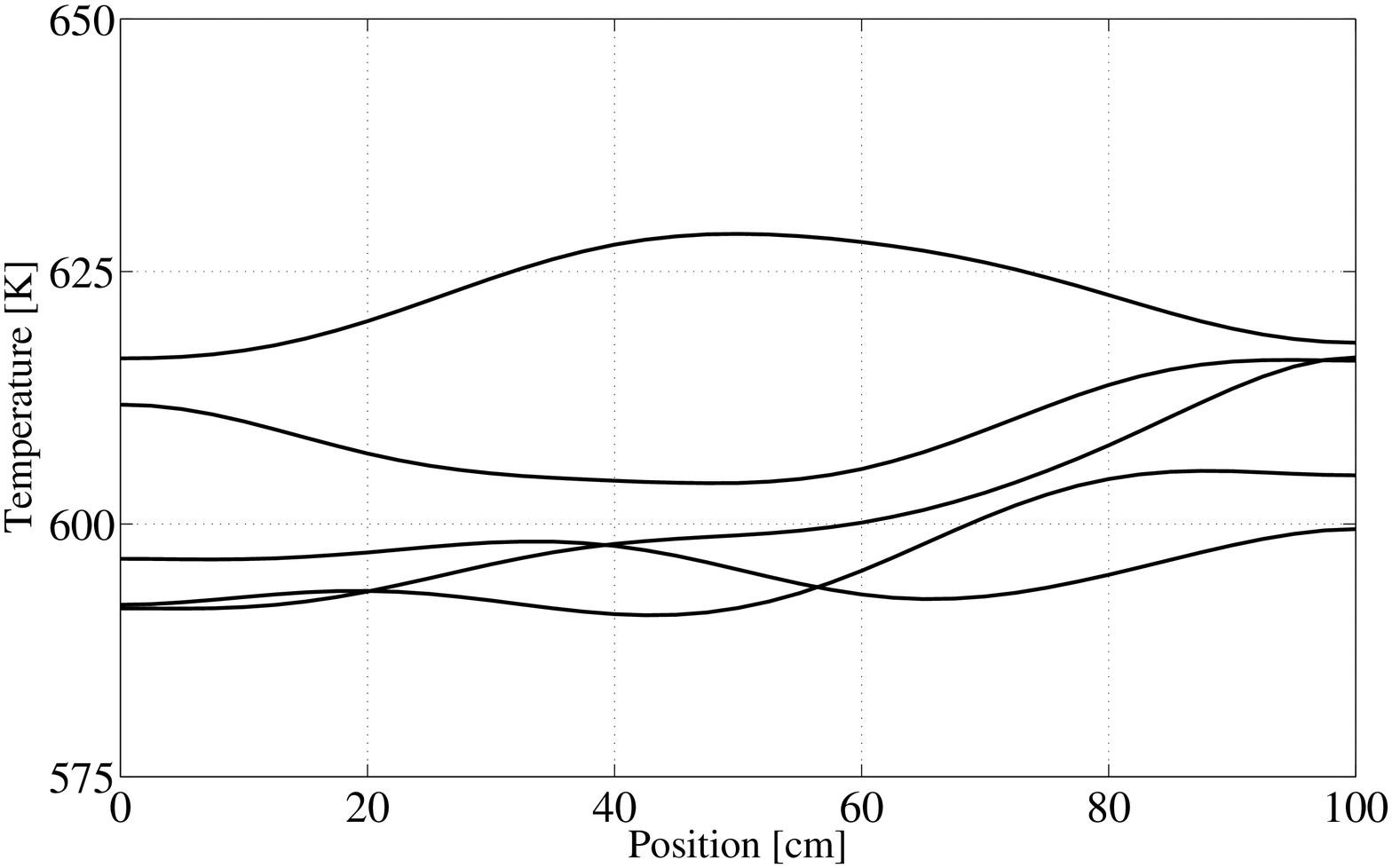}}
    \hfill
    \subfigure[Neutron flux.]{\includegraphics[width=0.8\textwidth]{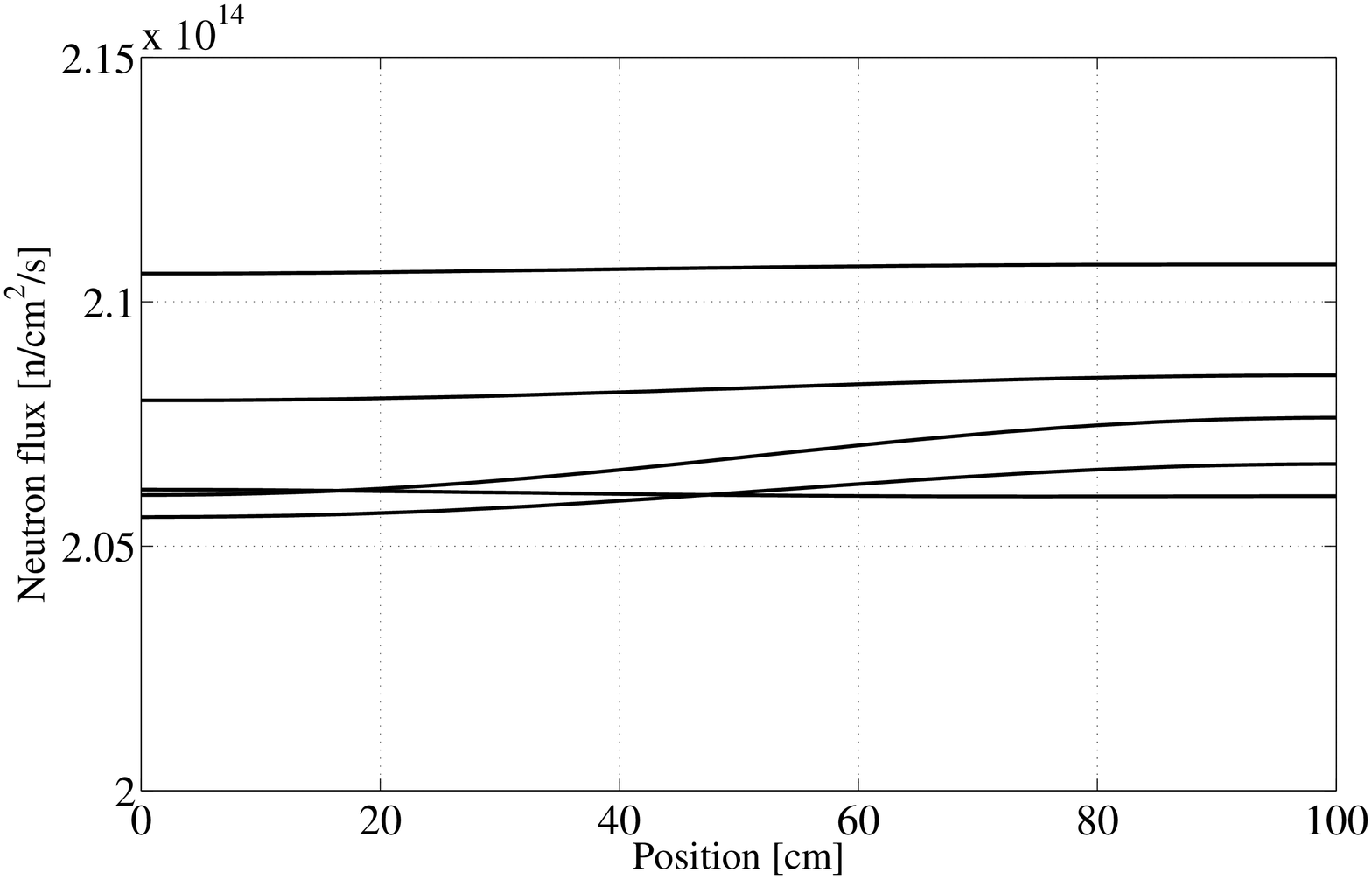}}
    \caption{PC-based simulation: five samples of the solution.}\label{fig:figure12}
  \end{center}
\end{figure}

Figure~\ref{fig:figure12} shows a few samples of the random temperature and neutron flux deduced from the PC expansions obtained as the output of the solution algorithm.
The samples of the input random variables used to synthesize the samples of the random temperature and neutron flux shown in Fig.~\ref{fig:figure12} were identical to those used to generate the samples shown in Fig.~\ref{fig:figure2}.
The similarity of the samples in Figs.~\ref{fig:figure2} and~\ref{fig:figure12} indicates that the PC-based surrogate model not only provides an accurate \textit{global} representation of the multiphysics model but also is capable of accurately reproducing a \textit{sample-wise} response.

\subsection{Convergence analysis}

\begin{figure}[htp]
  \begin{center}
    \subfigure[Reduced dimension $d$.]{\includegraphics[width=0.55\textwidth]{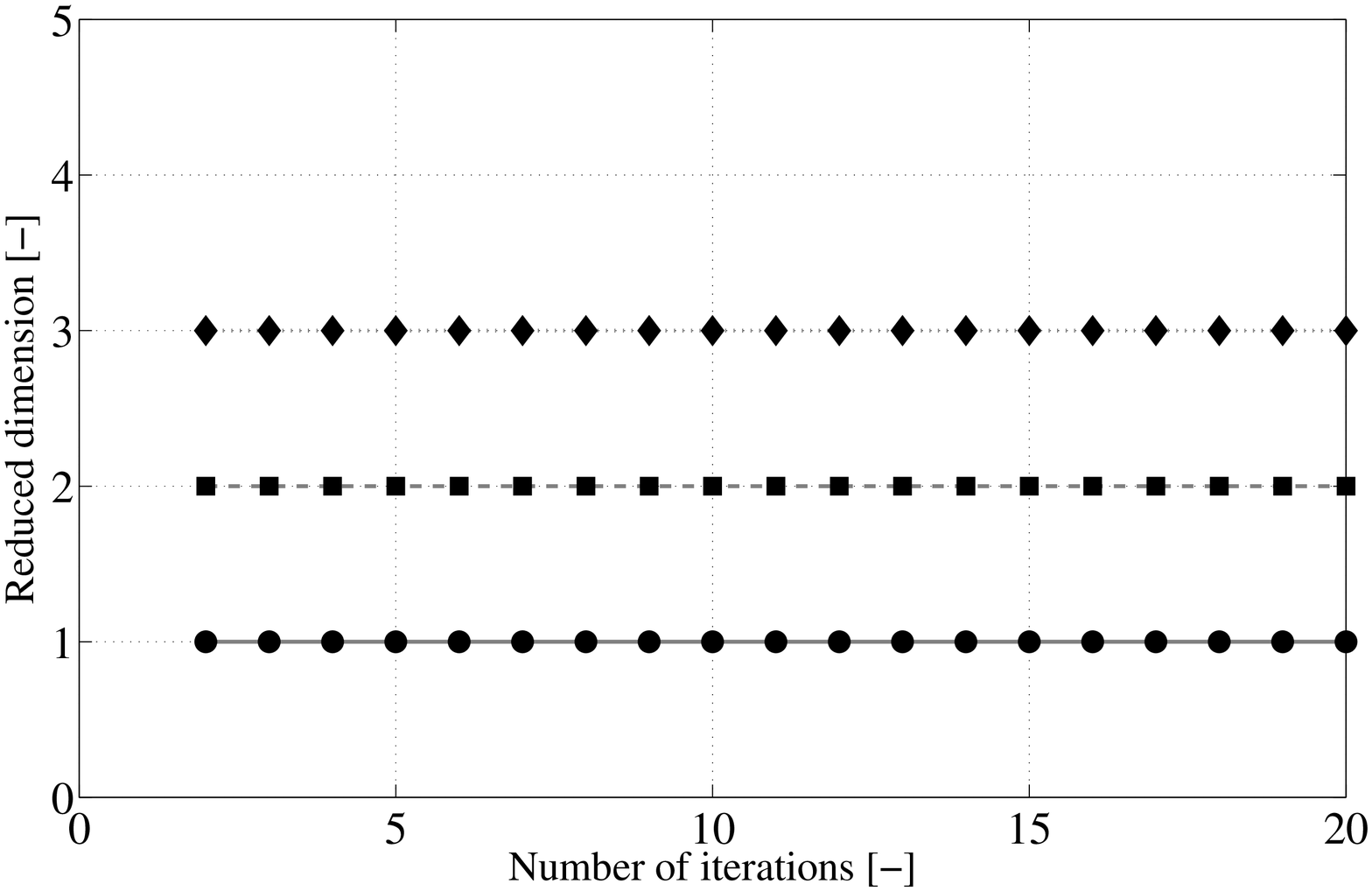}}
    \vfill
    \subfigure[$\ell\mapsto\sqrt{\frac{1}{MC}\sum_{k=1}^{MC}\|\boldsymbol{T}^{\ell}(\boldsymbol{\xi}_{k})-\widehat{\boldsymbol{T}}{}^{\ell,p}(\boldsymbol{\xi}_{k})\|_{\boldsymbol{W}}^{2}}\Big/\sqrt{\frac{1}{MC}\sum_{k=1}^{MC}\|\boldsymbol{T}^{\ell}(\boldsymbol{\xi}_{k})\|_{\boldsymbol{W}}^{2}}$.]{\includegraphics[width=0.55\textwidth]{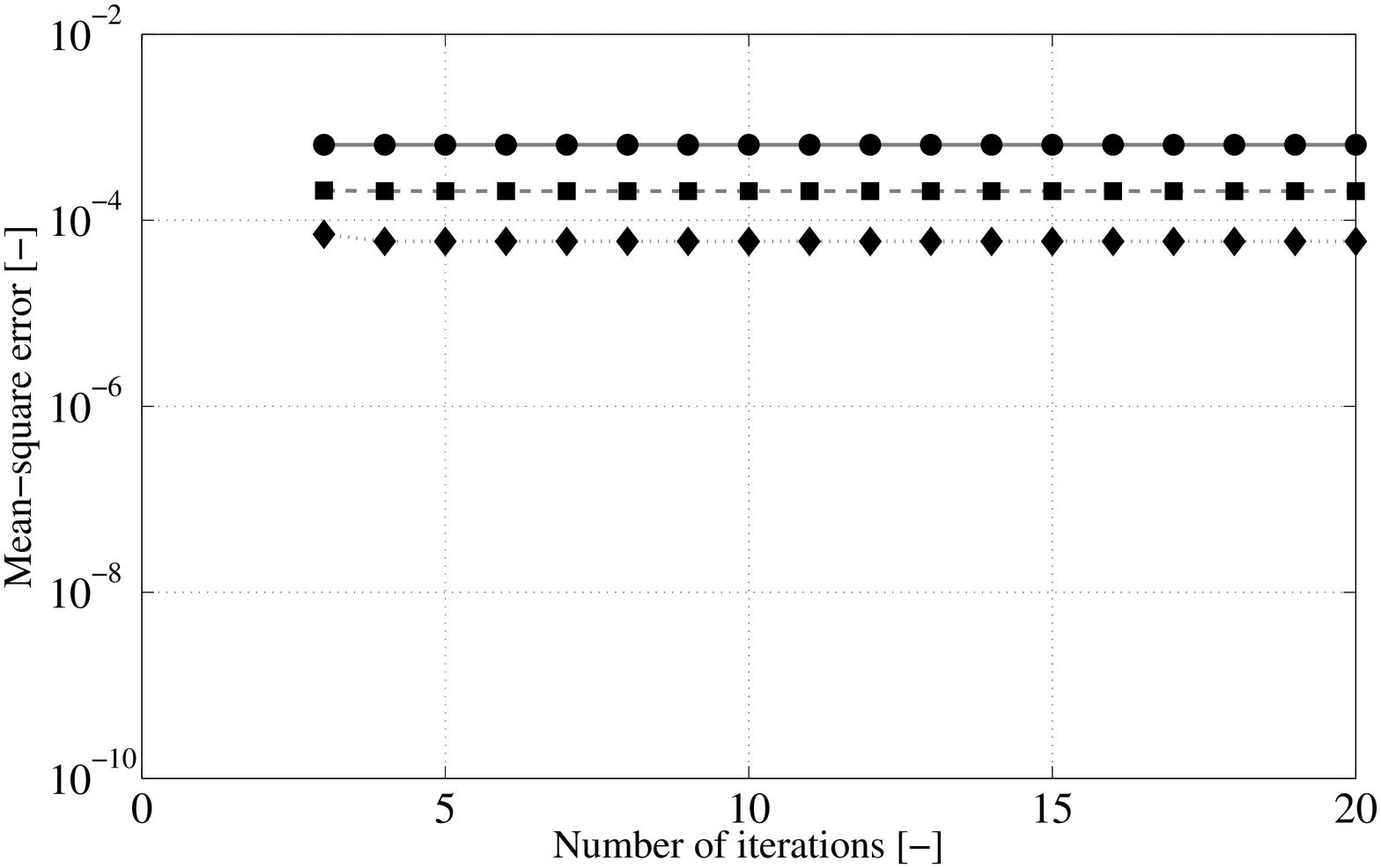}}
    \vfill
    \subfigure[$\ell\mapsto\sqrt{\frac{1}{MC}\sum_{k=1}^{MC}\|\boldsymbol{\Phi}^{\ell}(\boldsymbol{\xi}_{k})-\widehat{\boldsymbol{\Phi}}{}^{\ell,q}(\boldsymbol{\eta}^{\ell,p}(\boldsymbol{\xi}_{k}))\|_{\boldsymbol{W}}^{2}}\Big/\sqrt{\frac{1}{MC}\sum_{k=1}^{MC}\|\boldsymbol{\Phi}^{\ell}(\boldsymbol{\xi}_{k})\|_{\boldsymbol{W}}^{2}}$.]{\includegraphics[width=0.55\textwidth]{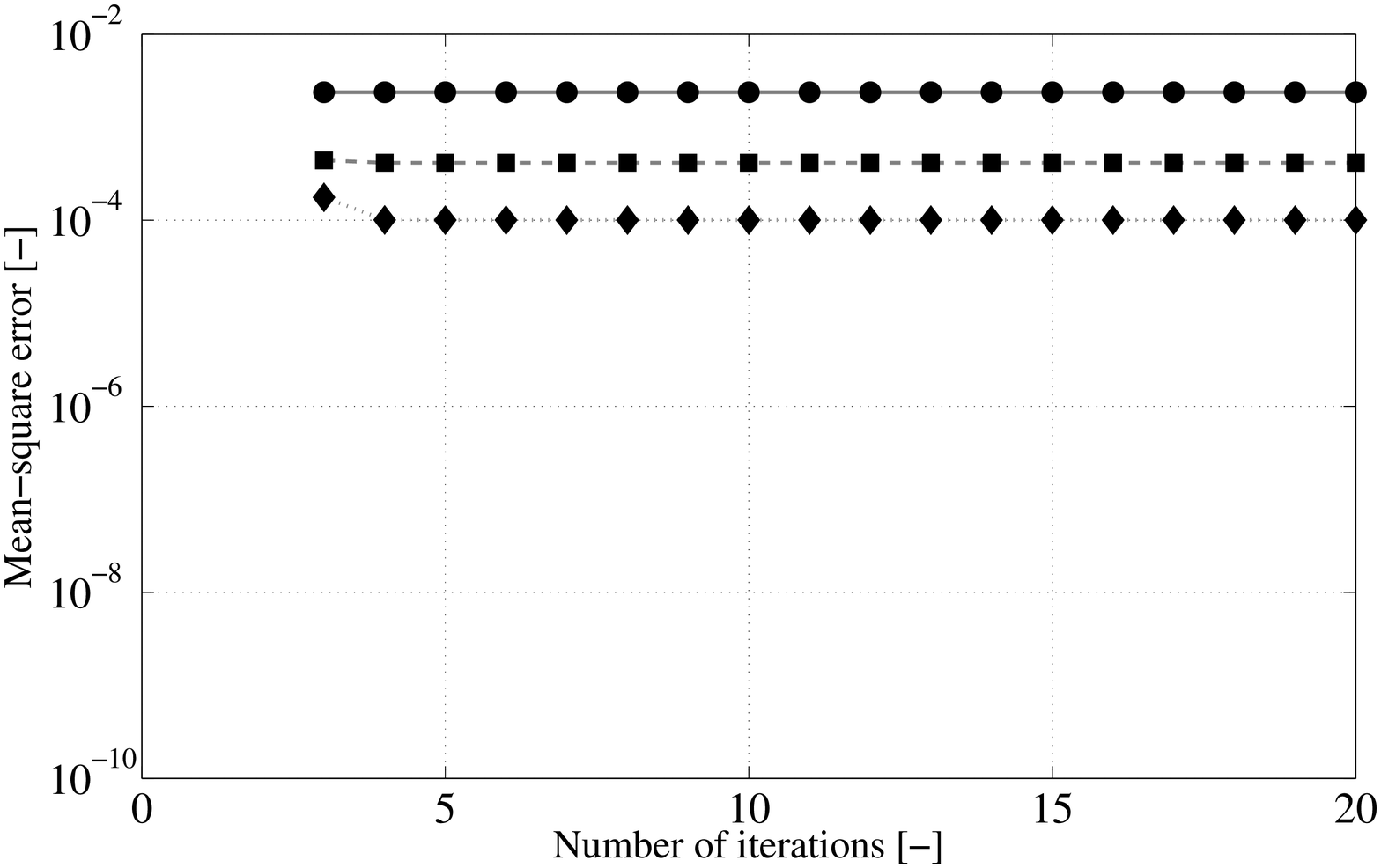}}
     \caption{Convergence analysis: (a) reduced dimension; and (b) and (c) mean-square distance between the successive approximations determined by the Monte Carlo and PC-based iterative methods for $\epsilon_{1}=0.20$~(circles), $\epsilon_{1}=0.05$~(squares), and $\epsilon_{1}=0.01$~(diamonds) and $\epsilon_{2}=0.0001$ as a function of the iteration.}\label{fig:figure13}
  \end{center}
\end{figure}

\begin{figure}[htp]
  \begin{center}
    \subfigure[Total degree $q$.]{\includegraphics[width=0.55\textwidth]{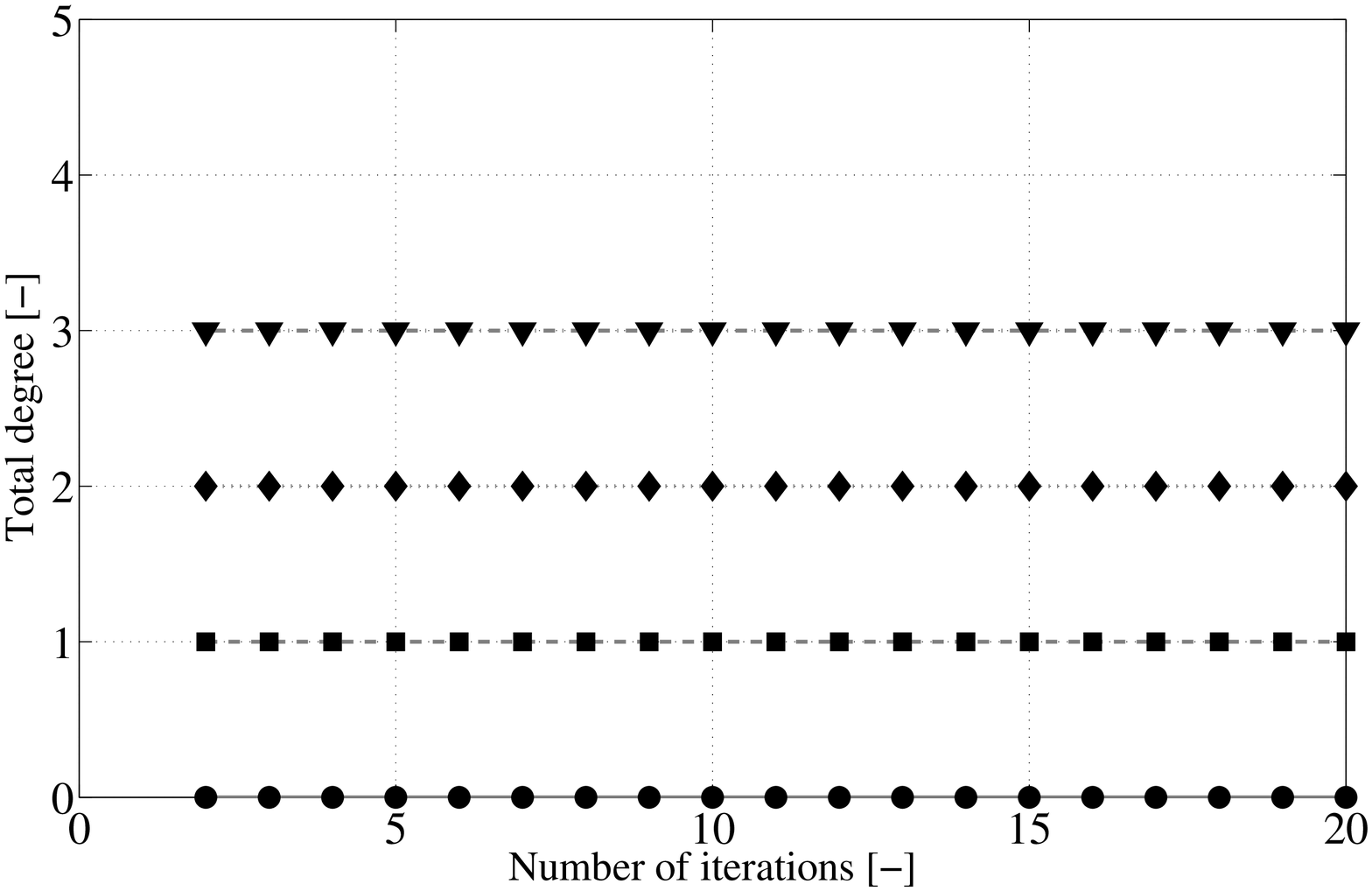}}
    \vfill
    \subfigure[$\ell\mapsto\sqrt{\frac{1}{MC}\sum_{k=1}^{MC}\|\boldsymbol{T}^{\ell}(\boldsymbol{\xi}_{k})-\widehat{\boldsymbol{T}}{}^{\ell,p}(\boldsymbol{\xi}_{k})\|_{\boldsymbol{W}}^{2}}\Big/\sqrt{\frac{1}{MC}\sum_{k=1}^{MC}\|\boldsymbol{T}^{\ell}(\boldsymbol{\xi}_{k})\|_{\boldsymbol{W}}^{2}}$.]{\includegraphics[width=0.55\textwidth]{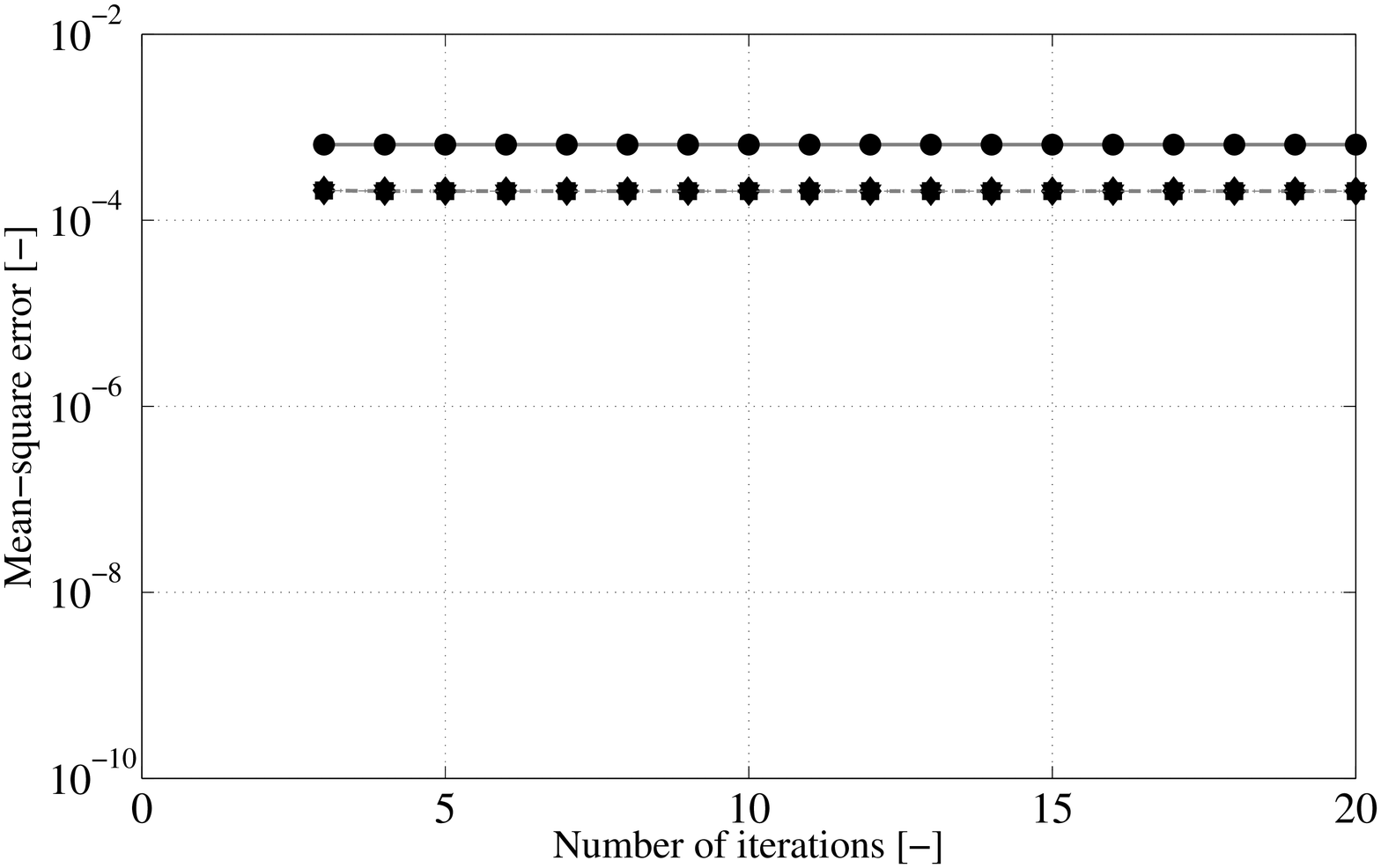}}
    \vfill
    \subfigure[$\ell\mapsto\sqrt{\frac{1}{MC}\sum_{k=1}^{MC}\|\boldsymbol{\Phi}^{\ell}(\boldsymbol{\xi}_{k})-\widehat{\boldsymbol{\Phi}}{}^{\ell,q}(\boldsymbol{\eta}^{\ell,p}(\boldsymbol{\xi}_{k}))\|_{\boldsymbol{W}}^{2}}\Big/\sqrt{\frac{1}{MC}\sum_{k=1}^{MC}\|\boldsymbol{\Phi}^{\ell}(\boldsymbol{\xi}_{k})\|_{\boldsymbol{W}}^{2}}$.]{\includegraphics[width=0.55\textwidth]{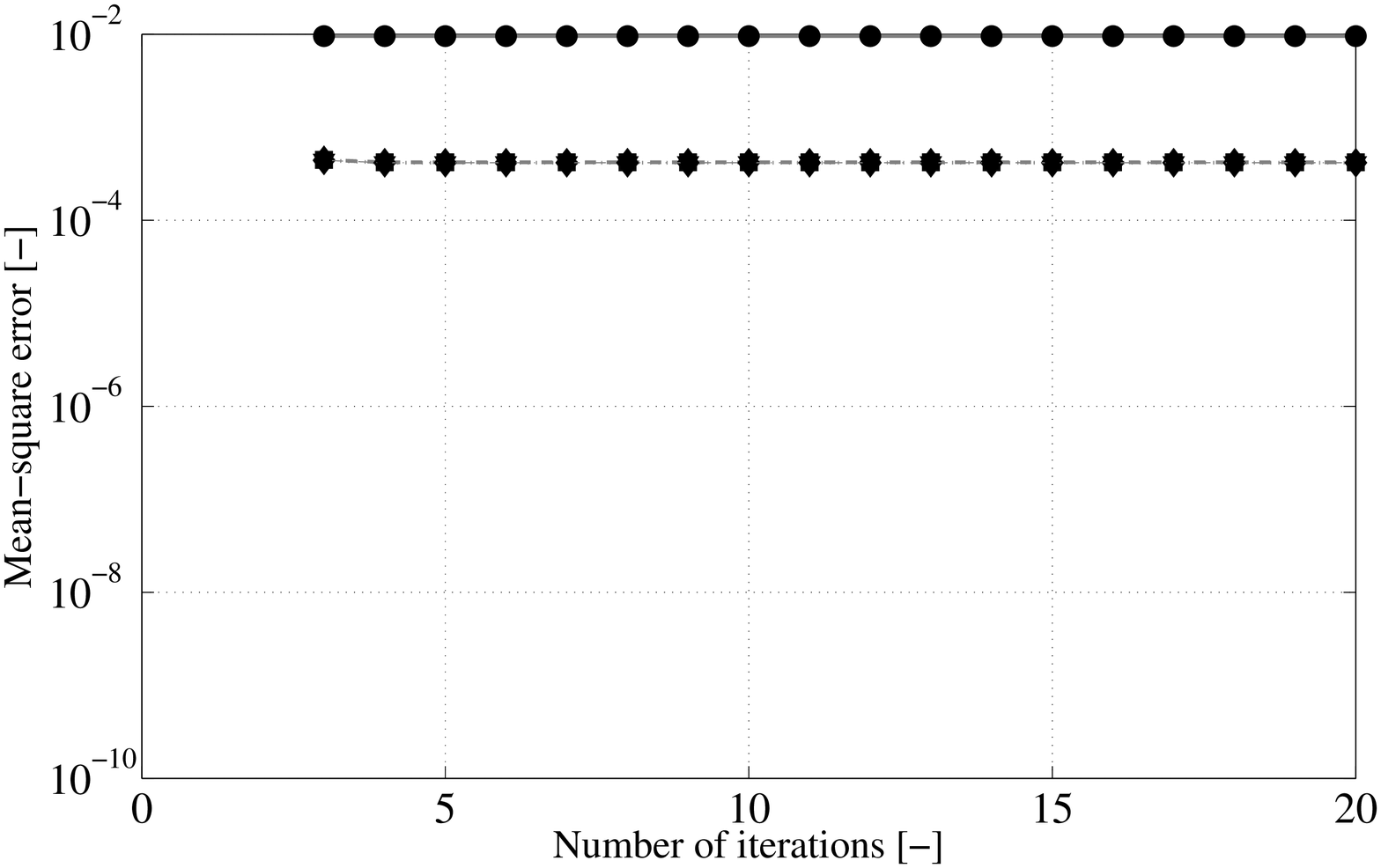}}
     \caption{Convergence analysis: (a) total degree; and (b) and (c) mean-square distance between the successive approximations determined by the Monte Carlo and PC-based iterative methods for $\epsilon_{1}=0.05$ and $\epsilon_{2}=1$~(circles), $\epsilon_{2}=0.01$~(squares), $\epsilon_{2}=0.0001$~(diamonds), and $\epsilon_{2}=0.000001$~(triangles) as a function of the iteration.}\label{fig:figure14}
  \end{center}
\end{figure}

We repeated the PC-based simulation for several values of the error tolerance levels; each of these values corresponded to specific accuracies that the KL decomposition of the random temperature and the PC expansion of the random neutron flux were required to maintain at each iteration. 
Figures~\ref{fig:figure13}(a) and~\ref{fig:figure14}(a) indicate that this KL decomposition and the PC expansion systematically retained more terms when higher accuracy was required.

Further, Figs.~\ref{fig:figure13}((b) and~(c)) and~\ref{fig:figure14}((b) and~(c)) indicate that the distance between the successive approximations determined by the Monte Carlo and PC-based iterative methods remained bounded as the iterations progressed and this distance can be reduced systematically by improving the accuracy of the KL decomposition of the random temperature and the PC expansion of the random neutron flux by decreasing the respective error tolerance levels.

\subsection{Concluding remarks}
The proposed methodology enabled the solution of the neutronics subproblem in a reduced-dimensional space because the KL decomposition was able to extract a low-dimensional representation of the random temperature as it passed from the heat subproblem to the neutronics subproblem.
While accuracy was maintained, the solution of the neutronics subproblem in a reduced-dimensional space resulted in computational gains, for two reasons.
First, it enabled the accurate representation of the random neutron flux by a PC expansion that had only a few terms.
Second, it enabled the computation of the coordinates in the PC expansion of the random neutron flux using a quadrature rule that had only a few nodes, thus requiring the solution of only a few samples of the neutronics subproblem per iteration. 

\section{Conclusion}
While most coupled models can be expected to be affected by a large number of sources of uncertainty, information exchanged between subproblems and iterations often resides in a considerably lower dimensional space than the sources themselves. 
In this work, we thus used a dimension-reduction technique to extract a low-dimensional representation of information as it passes from subproblem to subproblem and from iteration to iteration, and we proposed measure-transformation techniques that allows implementations to exploit this dimension reduction to achieve a computationally efficient solution of the subproblems in a reduced-dimensional space.
The effectiveness of the proposed methodology was demonstrated on a multiphysics problem relevant to nuclear reactors.

\section*{Acknowledgements}
This work was supported by DOE through an ASCR grant. 
The authors would also like to thank the anonymous reviewer for the valuable suggestions and Professor Christian Soize for relevant discussions during the final stages of the preparation of this paper. 

\bibliography{multiphysics3}

\end{document}